\documentclass{article}
\usepackage[utf8]{inputenc}
\usepackage{amsfonts,amsmath,amssymb,mathtools}
\usepackage[left=2cm,right=2cm,top=2cm,bottom=2cm,bindingoffset=0cm]{geometry}
\usepackage[russian,english]{babel}
\usepackage{graphicx}
\usepackage{svg}
\usepackage{authblk}
\usepackage{xcolor}
\usepackage{hyperref}
\usepackage{slashed}
\usepackage{accents}
\newcommand{\dbwidetilde}[1]{\widetilde{\widetilde{#1}}}
\usepackage{braket}
\usepackage{wrapfig}
\usepackage{multirow}
\usepackage{comment}
\usepackage{pb-diagram}
\usepackage{amsthm}
\usepackage{cancel}
\usepackage{tikz}\usetikzlibrary{tikzmark}
\usetikzlibrary{fadings} 
\usetikzlibrary{positioning}
\usetikzlibrary{backgrounds} 
\usetikzlibrary[backgrounds]
\usetikzlibrary{decorations.pathreplacing}
\usetikzlibrary{arrows}
\definecolor{airforceblue}{rgb}{0.36, 0.54, 0.66}	
\definecolor{beige}{rgb}{0.96, 0.96, 0.86}
\definecolor{bittersweet}{rgb}{1.0, 0.44, 0.37}
\definecolor{melon}{rgb}{0.99, 0.74, 0.71}
\definecolor{mustard}{rgb}{1.0, 0.86, 0.35}
\definecolor{lava}{rgb}{0.81, 0.06, 0.13}
\definecolor{magnolia}{rgb}{0.97, 0.96, 1.0}
\definecolor{lavendermist}{rgb}{0.9, 0.9, 0.98}
\definecolor{lavendergray}{rgb}{0.77, 0.76, 0.82}
\definecolor{palepink}{rgb}{0.98, 0.85, 0.87}
\definecolor{palesilver}{rgb}{0.79, 0.75, 0.73}
\definecolor{cadetgrey}{rgb}{0.57, 0.64, 0.69}
\definecolor{anti-flashwhite}{rgb}{0.95, 0.95, 0.96}
\colorlet{Light0anti-flashwhite}{anti-flashwhite!70!white}
\colorlet{Lightanti-flashwhite}{anti-flashwhite!50!white}
\colorlet{Light2anti-flashwhite}{anti-flashwhite!30!white}
\definecolor{linkcolor}{rgb}{0,0,1}
\definecolor{urlcolor}{rgb}{0,0,1}

\usepackage{tikz-cd}
\usetikzlibrary{decorations.markings}
\usepackage{array}
\usepackage{ytableau}
\usepackage{longtable}
\usepackage{empheq}
\usepackage{nicematrix}
\usepackage{multicol}
\usepackage{graphicx} 
\usepackage{subfig}
\usepackage{mathrsfs}
\usepackage{diagbox}
\usepackage[vcentermath]{youngtab}
\usepackage[natbib=true, backend = bibtex, style=numeric-comp, sorting=none,maxbibnames=99]{biblatex}
\hypersetup{pdfstartview=FitH,  linkcolor=linkcolor,urlcolor=urlcolor,citecolor=urlcolor, colorlinks=true}

\newcommand\bem{\begin{pmatrix}}
\newcommand\eem{\end{pmatrix}}
\newcommand\beq{\begin{equation}}
\newcommand\eeq{\end{equation}}
\newcommand\beqs{\begin{equation*}}
\newcommand\eeqs{\end{equation*}}

\newcommand{\sgn}{\,\text{sgn}}

\date{}

\usepackage{braids}
\usepackage{array}
\usepackage{ytableau}
\usepackage{longtable}
\usepackage{empheq}
\usepackage{multicol}
\usepackage{mathrsfs}
\usepackage{diagbox}

\definecolor{bright-gray}{rgb}{0.89, 0.90, 0.92}	

\def\be{\begin{eqnarray}}
\def\ee{\end{eqnarray}}

\def\bphi{\bar \phi}

\definecolor{red}{rgb}{1,0,0}
\definecolor{orange}{rgb}{1,0.5,0}
\definecolor{violet}{rgb}{0.7,0,1}

\newtheorem{theorem}{Theorem}[section]

\tikzstyle{every picture}+=[remember picture]

\begin{document}

\title{\bf Analogue of Goeritz matrices for computation\\  of bipartite HOMFLY-PT polynomials 
}

\author[2,3,5]{{\bf A. Anokhina}\thanks{\href{mailto:anokhina@itep.ru}{anokhina@itep.ru}}}
\author[1,4]{{\bf D. Korzun}\thanks{\href{mailto:korzun.dv@phystech.edu}{korzun.dv@phystech.edu}}}
\author[1,2,3,4,5]{{\bf E. Lanina}\thanks{\href{mailto:lanina.en@phystech.edu}{lanina.en@phystech.edu}}}
\author[1,2,3,5]{{\bf A. Morozov}\thanks{\href{mailto:morozov@itep.ru}{ morozov@itep.ru}}}

\affil[1]{Moscow Institute of Physics and Technology, 141700, Dolgoprudny, Russia}
\affil[2]{Institute for Information Transmission Problems, 127051, Moscow, Russia}
\affil[3]{NRC ”Kurchatov Institute”, 123182, Moscow, Russia}
\affil[4]{Saint Petersburg University, 199034, St. Petersburg, Russia}
\affil[5]{Institute for Theoretical and Experimental Physics, 117218, Moscow, Russia}
\renewcommand\Affilfont{\itshape\small}

\vspace{5cm}

\maketitle

\vspace{-8cm}

\begin{center}
	\hfill MIPT/TH-11/25\\
	\hfill ITEP/TH-14/25\\
	\hfill IITP/TH-12/25
\end{center}

\vspace{4.5cm}

\begin{abstract}

{
The Goeritz matrix is an alternative to the Kauffman bracket and, in addition, makes it possible to calculate the Jones polynomial faster with some minimal choice of a checkerboard surface of a link diagram. We introduce a modification of the Goeritz method that generalizes the Goeritz matrix for computing the HOMFLY--PT polynomials for any $N$ in the special case of bipartite links. Our method reduces to purely algebraic operations on matrices, and therefore, it can be easily implemented as a computer program. Bipartite links form a rather large family including a special class of Montesinos links constructed from the so-called rational tangles. We demonstrate how to obtain a bipartite diagram of such links and calculate the corresponding HOMFLY--PT polynomials using our developed generalized Goeritz method. 
}
\end{abstract}

\tableofcontents

\section{Introduction}
\setcounter{equation}{0}

Goeritz matrices are the connection between combinatorial and topological methods of research in knot theory. They first appeared in the work of Lebrecht Goeritz, where a connection between knot diagrams and quadratic forms was demonstrated~\cite{goeritz1933knoten}. The appearing quadratic forms and, consequently, symmetric matrices give various knot invariants — the knot signature~\cite{gordon1978signature}, the knot determinant~\cite{goeritz1933knoten}, and the Kauffman bracket~\cite{Kauff,boninger2023jones}. The relationship between the Goeritz matrix and the Kauffman bracket provides a method for computing the famous Jones polynomial $J(q)$~\cite{jones1985polynomial,jones1987hecke,boninger2023jones}, see Section~\ref{sec:JonesFromGoeritzMatrix}. This invariant is of particular interest, since its discovery revealed an equivalence between quantum knot polynomials and observables in topological quantum field theories~\cite{Witten1988hf}. Despite numerous studies, many questions regarding the Jones polynomial remain open. For instance, it is still unclear which topological information about a knot it contains.~\cite{jones2005jones}. It is possible that further investigation into the relationship between Jones polynomials and Goeritz matrices will improve our understanding of this matter.

There is a generalization of the Jones polynomial to the case of two variables --- the HOMFLY--PT polynomials\footnote{The HOMFLY–PT polynomials $H(A,q)$ and the Jones polynomials $J(q) = H(q^2, q)$ can be defined via the so-called skein relation:
	
\be\label{Hrelation}
\begin{picture}(300,25)(-70,-30)
	
	\put(0,-30){
		\put(-40,10){\mbox{$A\cdot H\Bigg( $}}
		\put(15,10){\mbox{$\Bigg) $}}
		\put(13,0){\vector(-1,1){24}}
		\put(-11,0){\line(1,1){10}}
		\put(3,14){\vector(1,1){10}}
		
		\put(64,0){
			\put(-37,10){\mbox{$= \ \, A^{-1}\cdot H\Bigg($}}
			\put(39,10){\mbox{$\Bigg)$}}
			
			\put(1,0){
				\put(13,0){\vector(1,1){24}}
				\put(37,0){\line(-1,1){10}}
				\put(23,14){\vector(-1,1){10}}
			}
		}
		
		\put(152,0){
			\put(-39,10){\mbox{$+ \, (q-q^{-1})\cdot H\Bigg($}}
			\put(46,10){\mbox{$\Bigg)$}}
			\put(30,0){\vector(0,1){24}}
			\put(40,0){\vector(0,1){24}}
		}
	}
	
\end{picture}
\ee
The idea of these relations is as follows. One takes a neighborhood of a single crossing and replaces only that crossing in the diagram with two other possible diagram fragments from~\eqref{Hrelation}. One then expresses the desired polynomial in terms of the polynomials of the two other diagrams from~\eqref{Hrelation} and repeats this operation, along with ambient isotopies, until all diagrams are transformed into diagrams of the trivial knot (a knot without crossings), which is conventionally normalized to $D_N = \{A\}/\{q\}$, where $\{x\} = x - x^{-1}$.
	
} $H(A = q^N,q)$ ~\cite{freyd1985new, przytycki1987kobe}. We can obtain it in the same way as the Jones polynomial — by taking the quantum trace of the $\mathfrak{R}$-matrix representation of a knot~\cite{reshetikhin1990ribbon,turaev1990yang,reshetikhin1991invariants}. The HOMFLY–PT polynomial corresponds to a representation of the universal enveloping algebra $U_q(\mathfrak{sl}_N)$~\cite{mironov2012character}, while the Jones polynomial, in turn, corresponds to a representation of $U_q(\mathfrak{sl}_2)$ and can be obtained from the HOMFLY–PT polynomial by reduction under the substitution $A \to q^2$. In this regard, it is natural to assume the existence of a matrix analogous to the Goeritz matrix, such that this new matrix would allow one to compute the HOMFLY–PT polynomials.

In this paper, we generalize the Goeritz matrix and present a method for computing the HOMFLY–PT polynomials (in the fundamental representation) for a special class of bipartite knots~\cite{BipKnots}. We will also refer to the generalized matrix as the quaternary Goeritz matrix, since its entries consist of four independent types of numbers, corresponding to different bipartite crossings and their colorings.

The quaternary Goeritz matrix for the HOMFLY–PT polynomial of a bipartite knot is related to the Goeritz matrix for the Jones polynomial of a precursor link, see Section~\ref{sec:BonNorm}. Since a single precursor with $n$ crossings corresponds to $2^n$ bipartite diagrams\footnote{Each positive crossing in the precursor link corresponds to two distinct bipartite crossings (similarly for negative crossings), see Fig.~\ref{fig:precursor-bip}.}, this generalization of the precursor Goeritz matrix to the bipartite case is not unique and depends on the choice of a particular bipartite diagram. At the same time, the inverse reduction is unique.

As our computation of bipartite HOMFLY–PT polynomials generalizes the algorithm of Boninger from~\cite{boninger2023jones}, we first explicitly explain that the Goeritz approach to computing the Jones polynomial is equivalent to the Kauffman bracket. The restriction to bipartite diagrams is motivated by the existence of a well-developed technique of planar decomposition for this class of knots~\cite{anokhina2024planar,anokhina2025bipartite,anokhina2025planar,anokhina2025khrbip}, which gives hope for the possibility of generalizing the Goeritz method to bipartite HOMFLY–PT polynomials. Despite equivalence to the Kauffman bracket, the Goeritz method is more convenient for practical computations, since it consists of purely algebraic (and easily programmable) matrix procedures and frees us from the need to work with knot diagrams and their planar decompositions.

In the subsections below, we will provide a brief introduction to the objects that we use in this text.

\subsection{Kauffman bracket and Jones polynomial: $N=2$}

In this text, we will use the symmetric form of the Kauffman bracket $\langle {\cal L} \rangle$ of a link $\cal L$ together with the condition that each closed cycle contributes a factor of $D_2 = q + q^{-1}$ (see Fig.~\ref{fig:Kauff}).

\begin{figure}[h!]
\begin{picture}(100,180)(-200,-125)

\put(0,20){
\put(-55,0){

\put(-60,-20){\line(1,1){18}}
\put(-20,-20){\line(-1,1){40}}
\put(-37.5,2){\line(1,1){18}}

\put(-68,22){\mbox{$i$}}  \put(-15,22){\mbox{$j$}}  \put(-68,-28){\mbox{$k$}}  \put(-15,-28){\mbox{$l$}}

\put(10,-2){\mbox{$=$}}
}

\put(-15,0){\mbox{$(-q)^{-1/2}$}}

\qbezier(30,20)(50,0)(70,20)
\qbezier(30,-20)(50,0)(70,-20)

\put(85,-2){\mbox{$+ \ \ \ (-q)^{+1/2} $}}

\put(25,0){
\qbezier(125,20)(145,0)(125,-20)
\qbezier(150,20)(130,0)(150,-20)
}
}

\put(-55,-50){

\put(-60,-20){\line(1,1){40}}
\put(-20,-20){\line(-1,1){18}}
\put(-60,20){\line(1,-1){18}}

\put(10,-2){\mbox{$=$}}

\put(-68,22){\mbox{$i$}}  \put(-15,22){\mbox{$j$}}  \put(-68,-28){\mbox{$k$}}  \put(-15,-28){\mbox{$l$}}
}

\put(0,-50){

\put(-15,0){\mbox{$(-q)^{+1/2}$}}

\qbezier(30,20)(50,0)(70,20)
\qbezier(30,-20)(50,0)(70,-20)

\put(85,-2){\mbox{$+ \ \ \ (-q)^{-1/2} $}}

\put(25,0){
\qbezier(125,20)(145,0)(125,-20)
\qbezier(150,20)(130,0)(150,-20)
}
}

\put(0,-40){

\put(-100,-70){\circle{30}}
\put(-75,-70){\mbox{$=\ \ \ D_2 \ = \ q+q^{-1}$}}
}

\end{picture}
\caption{\footnotesize The Kauffman bracket is the planar decomposition of an $\mathfrak{R}$-matrix vertex in the fundamental representation of $U_q(\mathfrak{sl}_2)$. In this case ($N = 2$), the fundamental representation is isomorphic to its conjugate and therefore the tangles in the diagram have no orientation.
}
\label{fig:Kauff}
\end{figure}

Reidemeister moves generate equivalence relations on knots. Diagrams of the same knot are isotopic with respect to Reidemeister moves (see Fig.\,\ref{fig:Reid-moves}). The Kauffman bracket is invariant under the second and third Reidemeister moves and changes under the first move.

\begin{figure}[h!]
    \centering
\begin{picture}(300,85)(65,-20)

\thicklines

\put(-33,0){\line(0,1){63}}

\put(-20,28){\mbox{$=$}}

\put(0,0){\line(0,1){25}}
\put(0,38){\line(0,1){25}}

\qbezier(0,25)(0,40)(13,40)
\qbezier(13,20)(8,19)(4,23)

\put(13,30){\oval(20,20)[r]}

\put(50,0){

\put(0,0){\line(0,1){20}}
\put(0,30){\line(0,1){33}}

\qbezier(13,20)(0,20)(0,30)
\qbezier(13,40)(5,40)(4,34)

\put(13,30){\oval(20,20)[r]}

\put(-18,28){\mbox{$=$}}

}

\put(0,-20){\mbox{I.R}}

\put(-30,0){

\put(170,0){

\put(-33,0){\line(0,1){63}}

\put(-20,0){\line(0,1){63}}

\put(-5,30){\mbox{$=$}}

\put(13,30){\oval(45,60)[r]}

\qbezier(50,60)(44,55)(38,50)
\qbezier(30,42)(20,32)(30,22)
\put(40,16){\line(1,-1){15}

\put(-25,-36){\mbox{II.R}}}

\put(57,30){\mbox{$=$}}

}

\put(270,0){

\put(13,30){\oval(45,60)[l]}

\put(-25,60){\line(1,-1){12}}

\put(10,-7){
\qbezier(-15,50)(-5,40)(-15,30)
}

\put(-29,1){\line(1,1){15}}

}
}

\put(240,-30){

\put(40,75){\line(1,0){57}}
    \put(113,75){\line(1,0){15}}
    \put(60,30){\line(1,1){60}}
    \put(50,90){\line(1,-1){10}}
    \put(70,70){\line(1,-1){10}}
    \put(90,50){\line(1,-1){20}}

\put(135,55){\mbox{$=$}}

\put(-40,0){
    \put(200,30){\line(1,1){60}}
    \put(225,45){\line(1,0){57}}
    \put(187,45){\line(1,0){18}}
    \put(215,90){\line(1,-1){18}}
    \put(242,63){\line(1,-1){13}}
    \put(265,40){\line(1,-1){10}}
    
\put(165,10){\mbox{III.R}}

}

}

\end{picture}
\caption{\footnotesize Reidemeister moves. There is only single realization of Reidemeister third move (III.R) shown in the figure.}
\label{fig:Reid-moves}
\end{figure}

Therefore, to obtain a knot invariant from the Kauffman bracket, it is necessary to normalize it by a factor that accounts for the first Reidemeister move. Specifically, the Jones polynomial for a link ${\cal L}$ is expressed in terms of the Kauffman bracket $\langle {\cal L} \rangle$ as follows:

\begin{equation}\label{JonesFromKauffman}
    J^{\cal L} = (-(-q)^{3/2})^{-{\rm W}({\cal D})}\cdot \langle {\cal L} \rangle\,,
\end{equation}
when ${\rm W}({\cal D})$ is the algebraic crossing number for the diagram  $\cal D$ of a link $\cal L$ also called the writhe number (being the difference of the number of positive and negative crossings, see Fig.\,\ref{fig:pos-neg-cr}).

\begin{figure}[h!]
	\centering
	\begin{picture}(300,65)(-135,-40)
		
		\put(-60,-20){\line(1,1){18}}
		\put(-20,-20){\vector(-1,1){40}}
		\put(-37.5,2.5){\vector(1,1){18}}
		
		\put(-45,-38){\mbox{\Large $+$}}
		
		\put(100,0){
			
			\put(-60,-20){\vector(1,1){40}}
			\put(-20,-20){\line(-1,1){18}}
			\put(-42.5,2.5){\vector(-1,1){18}}
			
			\put(-45,-38){\mbox{\Large $-$}}
		}
		
	\end{picture}
	\caption{\footnotesize Agreement on crossing signs. A positive crossing is indicated by a plus sign, a negative crossing is indicated by a minus sign.}
	\label{fig:pos-neg-cr}
\end{figure}

\subsection{Planar decomposition of a bipartite vertex for any $N$}\label{sec:BonNorm}

In~\cite{anokhina2024planar}, the Kauffman bracket was generalized to the case of an arbitrary $N$, but for a special class of bipartite knots consisting entirely of antiparallel\footnote{For $N\neq 2$, the orientation becomes significant.} {\it lock} elements (see Fig.~\ref{fig:pladeco}). Here, we again use the symmetric planar decomposition of a lock tangle.

\begin{figure}[h!]
\begin{picture}(100,175)(-150,-230)

\put(10,-85){
\put(-105,15){\vector(1,-1){12}} \put(-87,3){\vector(1,1){12}}
\put(-93,-3){\vector(-1,-1){12}} \put(-75,-15){\vector(-1,1){12}}
\put(-90,0){\circle*{6}}

\put(-60,-2){\mbox{$:=$}}
}

\put(10,-150){
\put(-105,15){\vector(1,-1){12}} \put(-87,3){\vector(1,1){12}}
\put(-93,-3){\vector(-1,-1){12}} \put(-75,-15){\vector(-1,1){12}}
\put(-90,0){\circle{6}}

\put(-60,-2){\mbox{$:=$}}
}

\put(0,-85){

\put(-10,0){

\put(-17,20){\line(1,-1){17}}\put(-17,20){\vector(1,-1){14}}   \put(0,3){\vector(1,1){17}}
\put(0,-3){\vector(-1,-1){17}}   \put(17,-20){\line(-1,1){17}} \put(17,-20){\vector(-1,1){14}}
\put(0,-3){\line(0,1){6}}

\put(30,-2){\mbox{$:=$}}

\qbezier(51,20)(75,-40)(97,-4)  \qbezier(102,4)(105,9)(110,20)
\qbezier(50,-20)(55,-9)(58,-4) \qbezier(63,4)(85,40)(110,-20)
\put(56,8){\vector(1,-2){2}} 
\put(90,-13){\vector(1,1){2}} 
\put(109,18){\vector(1,2){2}}
\put(105,-9){\vector(-1,2){2}} 
\put(70,13){\vector(-1,-1){2}} \put(51,-18){\vector(-1,-2){2}}

\put(130,-2){\mbox{$=$}}
}

\put(100,65){

\put(40,-67){\mbox{$\phi^{-1/2}$}}

\put(70,-60){\vector(1,0){40}}
\put(110,-70){\vector(-1,0){40}}
\put(120,-67){\mbox{$+\ \ \phi^{+1/2}$}}

\put(15,0){
\put(155,-50){\vector(0,-1){30}}
\put(165,-80){\vector(0,1){30}}
}}}


\put(0,-150){

\put(-10,0){

\put(-17,20){\line(1,-1){17}}
\put(-17,20){\vector(1,-1){14}}   \put(0,3){\vector(1,1){17}}
\put(0,-3){\vector(-1,-1){17}}   \put(17,-20){\line(-1,1){17}} 
\put(17,-20){\vector(-1,1){14}}
\put(-1,-3){\line(0,1){6}}  
\put(1,-3){\line(0,1){6}}

\put(30,-2){\mbox{$:=$}}

\qbezier(50,20)(55,9)(58,4) \qbezier(63,-4)(85,-40)(110,20)
\put(56,8){\vector(1,-2){2}} \put(90,-13){\vector(1,1){2}} \put(109,18){\vector(1,2){2}}
\qbezier(50,-20)(75,40)(97,4)  \qbezier(102,-4)(105,-9)(110,-20)
\put(104,-8){\vector(-1,2){2}} \put(70,13){\vector(-1,-1){2}} \put(51,-18){\vector(-1,-2){2}}

\put(130,-2){\mbox{$=$}}
}

\put(100,65){

\put(40,-67){\mbox{$\bphi^{-1/2}$}}

\put(70,-60){\vector(1,0){40}}
\put(110,-70){\vector(-1,0){40}}
\put(120,-67){\mbox{$+\ \ \bar \phi^{+1/2}$}}

\put(15,0){
\put(155,-50){\vector(0,-1){30}}
\put(165,-80){\vector(0,1){30}}
}}}

\put(20,-140){

\put(-100,-70){\circle{30}}
\put(-75,-70){\mbox{$=\ \ \ D_N \ = \ \frac{\{A\}}{\{q\}}$}}
}

\end{picture}
\caption{\footnotesize Vertical antiparallel lock from  \cite{anokhina2024planar}. Also drawn is the inverse lock made from inverse vertices. Here $\phi = A\{q\}$, $\bphi = -A^{-1}\{q\}$.
} \label{fig:pladeco}
\end{figure}

In topological fraiming HOMFLY-PT polynomial for a bipartite link ${\cal L}_{\rm bip}$ is expressed through the planar decomposition $\langle {\cal L}_{\rm bip} \rangle_{\rm pl}$:
\begin{equation}\label{HOMFLYNorm}
    H^{{\cal L}_{\rm bip}} = A^{-2(N_\bullet - N_\circ)} \cdot \phi^{N_\bullet/2}\bphi^{N_\circ/2}\langle {\cal L}_{\rm bip} \rangle_{\rm pl}\,,
\end{equation}
where $\phi = A\{q\}$, $\bphi = -A^{-1}\{q\}$, $N_\bullet$ and $N_\circ$ are the numbers of the upper and lower lock elements in Fig.\,\ref{fig:pladeco}, respectively, as represented in the bipartite diagram. We recall that throughout this text we use the notation $\{x\} = x - x^{-1}$.

It should be noted that under the substitution $\phi \rightarrow -q$, $\bphi \rightarrow -q^{-1}$, $D_N \rightarrow D_2$, the planar decomposition of lock tangles transforms into the Kauffman bracket shown in Fig.\,\ref{fig:Kauff}. We will refer to a {\it precursor link} as a link obtained by replacing the corresponding bipartite vertices with single crossings. Since the Kauffman bracket does not distinguish link orientation, a single crossing can be represented either as a vertically oriented lock or as a horizontally oriented inverse lock, as illustrated in Fig.\,\ref{fig:precursor-bip}.

\begin{figure}[h!]
	\centering
	\begin{picture}(300,180)(-60,-85)
		
		\put(-80,75){\mbox{Precursor diagram}}
		
		\put(110,75){\mbox{Bipartite diagram}}
		
		\put(-77,55){\mbox{Jones polynomial}}
		
		\put(93,55){\mbox{HOMFLY-PT polynomial}}
		
		\put(0,20){
			
			\put(-60,-20){\line(1,1){18}}
			\put(-20,-20){\line(-1,1){40}}
			\put(-37.5,2){\line(1,1){18}}
			
		}
		
		\put(0,-50){
			
			\put(-60,-20){\line(1,1){40}}
			\put(-20,-20){\line(-1,1){18}}
			\put(-60,20){\line(1,-1){18}}
			
		}
		
		\put(10,20){\mbox{\Large $\Longleftrightarrow$}}
		
		\put(10,-50){\mbox{\Large $\Longleftrightarrow$}}
		
		\put(90,20){
			
			\put(-17,20){\line(1,-1){17}}\put(-17,20){\vector(1,-1){14}}   \put(0,3){\vector(1,1){17}}
			\put(0,-3){\vector(-1,-1){17}}   \put(17,-20){\line(-1,1){17}} \put(17,-20){\vector(-1,1){14}}
			\put(0,-3){\line(0,1){6}}
			
		}
		
		\put(90,-50){
			
			\put(-17,20){\line(1,-1){17}}\put(-17,20){\vector(1,-1){14}}   \put(0,3){\vector(1,1){17}}
			\put(0,-3){\vector(-1,-1){17}}   \put(17,-20){\line(-1,1){17}} \put(17,-20){\vector(-1,1){14}}
			\put(-1,-3){\line(0,1){6}}
			\put(1,-3){\line(0,1){6}}
			
		}
		
		\put(145,20){\mbox{or}}
		
		\put(145,-50){\mbox{or}}
		
		\put(200,20){
			
			\put(-20,17){\line(1,-1){17}}\put(-20,17){\vector(1,-1){14}}   \put(3,0){\vector(1,1){17}}
			\put(-3,0){\vector(-1,-1){17}}   \put(20,-17){\line(-1,1){17}} \put(20,-17){\vector(-1,1){14}}
			\put(-3,1){\line(1,0){6}}
			\put(-3,-1){\line(1,0){6}}
			
		}
		
		\put(200,-50){
			
			\put(-20,17){\line(1,-1){17}}\put(-20,17){\vector(1,-1){14}}   \put(3,0){\vector(1,1){17}}
			\put(-3,0){\vector(-1,-1){17}}   \put(20,-17){\line(-1,1){17}} \put(20,-17){\vector(-1,1){14}}
			\put(-3,0){\line(1,0){6}}
			
		}
		
	\end{picture}
	\caption{\footnotesize If in the expansion of the HOMFLY--PT polynomial in variables $\phi$, $\bphi$, and $D_N$ we apply the substitution $\phi \rightarrow -q$, $\bphi \rightarrow -q^{-1}$, $D_N \rightarrow D_2$, then the HOMFLY--PT polynomial reduces to the Jones polynomial of a precursor link. This follows from the equivalence between the planar decomposition of a bipartite vertex (see Fig.\,\ref{fig:pladeco}) and the Kauffman bracket for a single crossing (see Fig.\,\ref{fig:Kauff}). We will refer to bipartite crossings as positive if they correspond to a positive precursor crossing. Furthermore, we call lock tangles vertical or horizontal depending on the orientation of their connecting bar.}
	\label{fig:precursor-bip}
\end{figure}

\paragraph{Text plan.} In this section, we have introduced the main objects that we work with throughout the text. The main result of our work is an algorithm for computing HOMFLY--PT polynomials using a quaternary Goeritz matrix proposed herein. However, it do not appear until Section~\ref{sec:HOMFLYfromGoeritz}; before that we must clarify a number of concepts that are used along the way. Our algorithm can be applied to an arbitrary bipartite link. Therefore, in Section~\ref{sec:bip-diag-obtain}, we describe how bipartite link diagrams can be obtained. In Section~\ref{sec:Goeritz-mat}, we introduce the Goeritz matrix and in Section~\ref{sec:JonesFromGoeritzMatrix}, we show that the action of the Kauffman bracket on link diagrams can be represented as transformations of Goeritz matrices. The method of obtaining Jones polynomials from a Goeritz matrix follows the approach of~\cite{boninger2023jones}. However, the logical order is reversed: we begin with the Kauffman bracket and transfer its action onto Goeritz matrices. Then, in Section~\ref{sec:Jones-ex} we present examples of computing Jones polynomials from a Goeritz matrix. Based on the developed understanding of the relationship between the Kauffman bracket and the algorithm on Goeritz matrices, we generalize the Goeritz matrix in Section~\ref{sec:HOMFLYfromGoeritz} so that it yields HOMFLY--PT polynomials for a distinguished class of bipartite links. Finally, in Section~\ref{sec:HOMFLY-ex}, we provide examples of computing HOMFLY--PT polynomials using the proposed quaternary Goeritz matrix.

\section{Obtaining bipartite diagrams}
\label{sec:bip-diag-obtain}
\setcounter{equation}{0}

At present, there is no general algorithm for constructing a bipartite diagram of a knot from an arbitrary knot diagram. Moreover, there is no algorithm for determining whether such a realization exists\footnote{Note that not all knots are bipartite. For example, the knots $9_{35}$, $9_{46}$, and $10_{74}$ are not bipartite. Additional examples and justification can be found in~\cite{lewark2016new}.}.  In this work, however, we focus on specific diagrams and classes of knots on which we illustrate our method: two-strand knots\footnote{By Alexander theorem, any link can be represented as the closure of a braid. Here, we refer to links that admit a representation as two-strand braids.} (also known as torus knots $T[2,p]$), twist knots ${\rm Tw}_m$, rational knots $R\!\left(\frac{p}{q}\right)$, and certain pretzel knots $P(t_1,t_2,\dots,t_n)$ and Montesinos knots $K\!\left(\frac{p_1}{q_1}, \frac{p_2}{q_2}, \dots, \frac{p_n}{q_n}\right)$. For the examples considered, we provide an explicit algorithm for obtaining bipartite diagrams which are used in subsequent sections to compute bipartite HOMFLY--PT polynomials.

\subsection{Algorithm for  two-strand knots}\label{sec:TwoBraidsKnotsAlgorithm}

Let us demonstrate that any two-strand torus link\footnote{Diagrams with a single crossing are trivially bipartized (i.e. made bipartite) by a first Reidemeister move of an appropriate orientation; of course, this transforms the diagram into that of the unknot.} $T[2,p]$ has a bipartite diagram. 

\begin{figure}[h!]
		\centering	
		\includegraphics[width =0.6\linewidth]{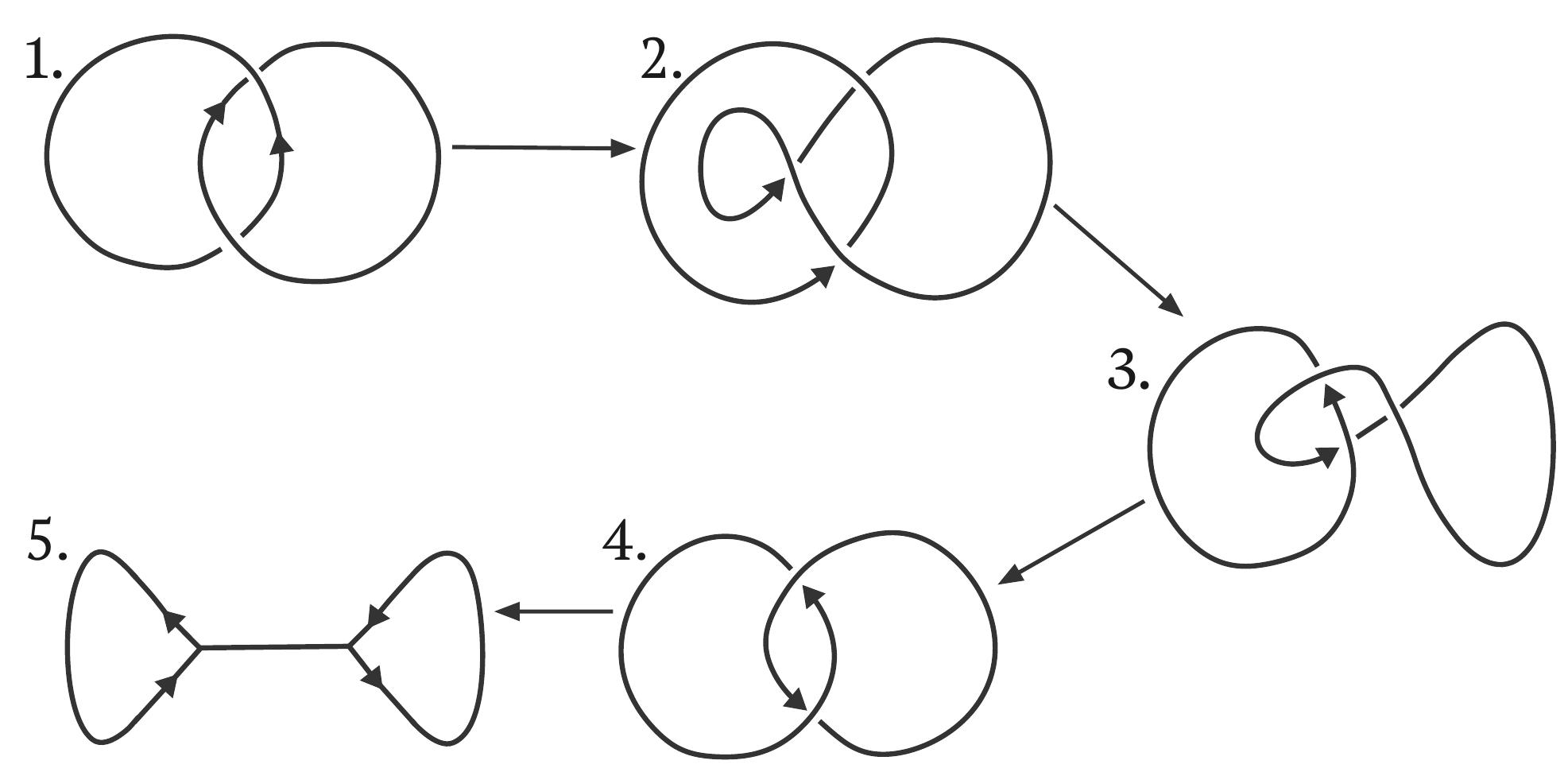}
        \caption{\footnotesize{Construction of a bipartite diagram for the Hopf link $T[2,2]$.}}
        \label{fig:BipHopf}
	\end{figure}

\begin{wrapfigure}[18]{r}{0.72\textwidth}
	\centering
	\vspace{-0.0cm}
	\includegraphics[width =1.07\linewidth]{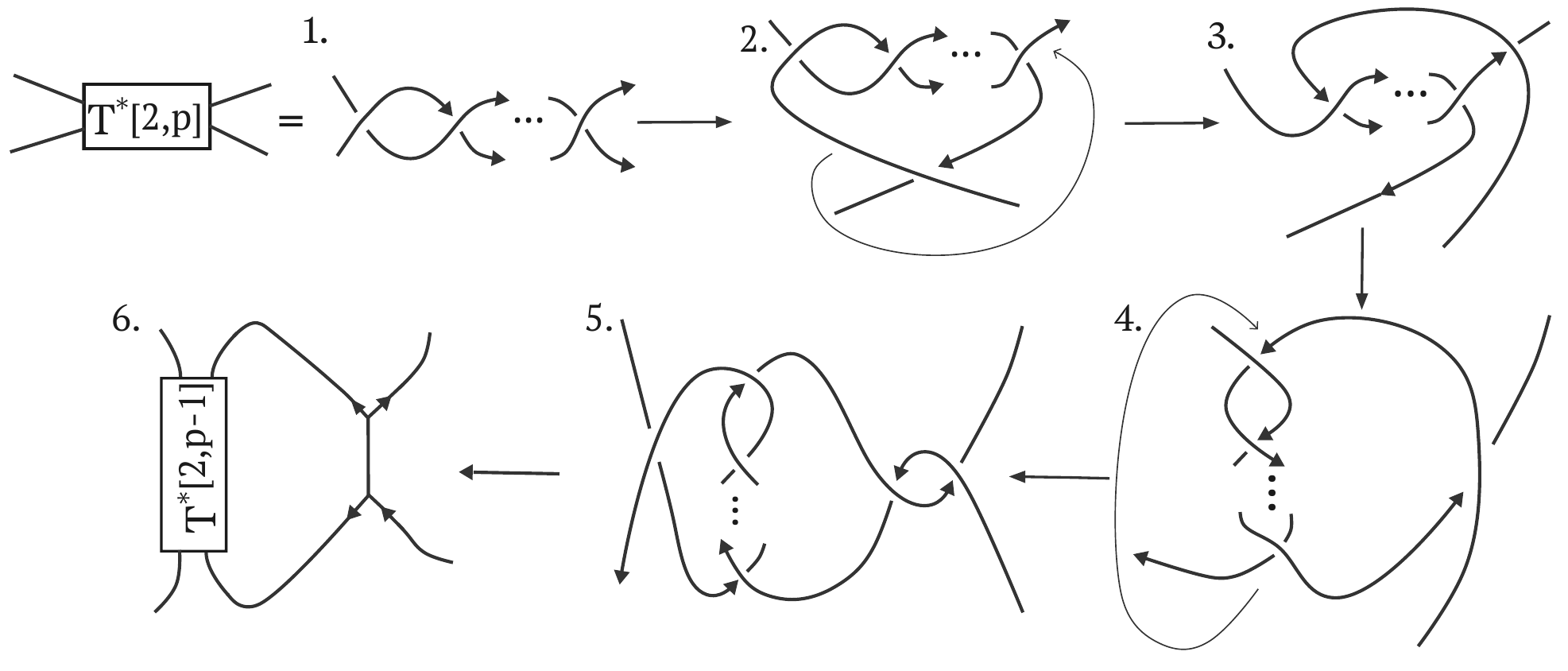}
	\caption{\footnotesize{Iteration of the algorithm for constructing a bipartite realization for $T^*[2,p]$.}}
	\label{fig:BipTp}
\end{wrapfigure}

To begin, we consider the simple examples of the Hopf link and the trefoil, and then show that the bipartization\footnote{By bipartization we mean a sequence of transformations that produces an explicitly bipartite diagram from an ordinary diagram.} of an arbitrary $T[2,p]$ link reduces to the bipartization of a 4-tangle for the trefoil $T^*[2,3]$, as illustrated in Fig.\,\ref{fig:BipTp}. For the Hopf link, there are parallel and antiparallel orientations of the components relative to one another. The antiparallel orientation is bipartite by definition (see Fig.\,\ref{fig:pladeco}). The parallel orientation can, however, be converted into the antiparallel one by a few elementary moves, as shown in Fig.\,\ref{fig:BipHopf}.

\begin{figure}[h!]
	\centering	
	\includegraphics[width =0.6\linewidth]{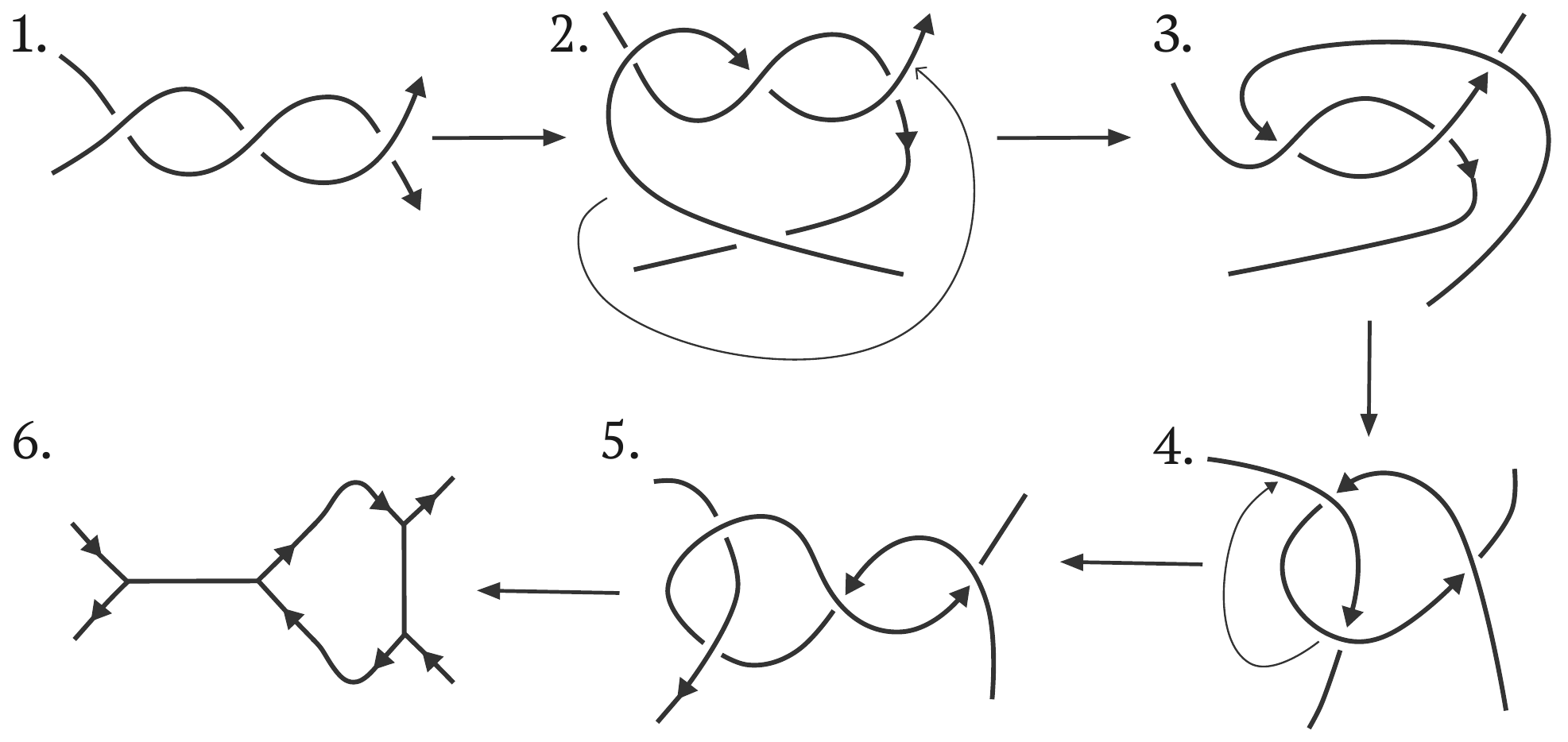}
	\caption{\footnotesize{Construction of a bipartite diagram for 4-tangle $T^*[2,3]$.}}
	\label{fig:BipT3}
\end{figure}

A bipartite realization of the trefoil is well known, it is the twist knot ${\rm Tw}_{2}$. Among other places, it can be found in~\cite{mironov2012character}. It can be obtained directly without difficulty, as we illustrate in the example of the tangle $T^*[2,3]$ in Fig.\,\ref{fig:BipT3}. The trefoil knot is obtained by pairwise closure of the upper and lower endpoints.

To explain the bipartization algorithm for the 4-tangle $T^*[2,p]$ in the case of an arbitrary $p$, we provide a sequence of transformations that replace a single ordinary crossing in the tangle $T^*[2,p]$ with a single bipartite crossing (see Fig.\,\ref{fig:BipTp}). As in the previous example, a link is obtained from a tangle by pairwise closure of the upper and lower endpoints.

In the first step $(1 \to 2)$, we pull one lower endpoint over the other so that the resulting crossing is inverse to the crossings in the original tangle (at the level of links, this corresponds to the first Reidemeister move). In the second step $(2 \to 3)$, we carry the lower endpoint positioned on top across all crossings of the tangle to the right edge of the diagram (the tangle is now decomposed into the 4-tangle $T^*[2,p-1]$ and a single crossing). In the third step $(3 \to 5)$, we rotate the 4-tangle $T^*[2,p-1]$ by $\pi$ in the plane perpendicular to the diagram plane. This adds one crossing on each side of $T^*[2,p-1]$. As a result, a bipartite crossing (with antiparallel orientation) appears on the right, while on the left, we obtain a $T^*[2,p-1]$ tangle rotated by $\frac{\pi}{2}$, to which we may apply the same sequence of transformations.

\begin{figure}[h!]
	\centering	
	\includegraphics[width =0.85\linewidth]{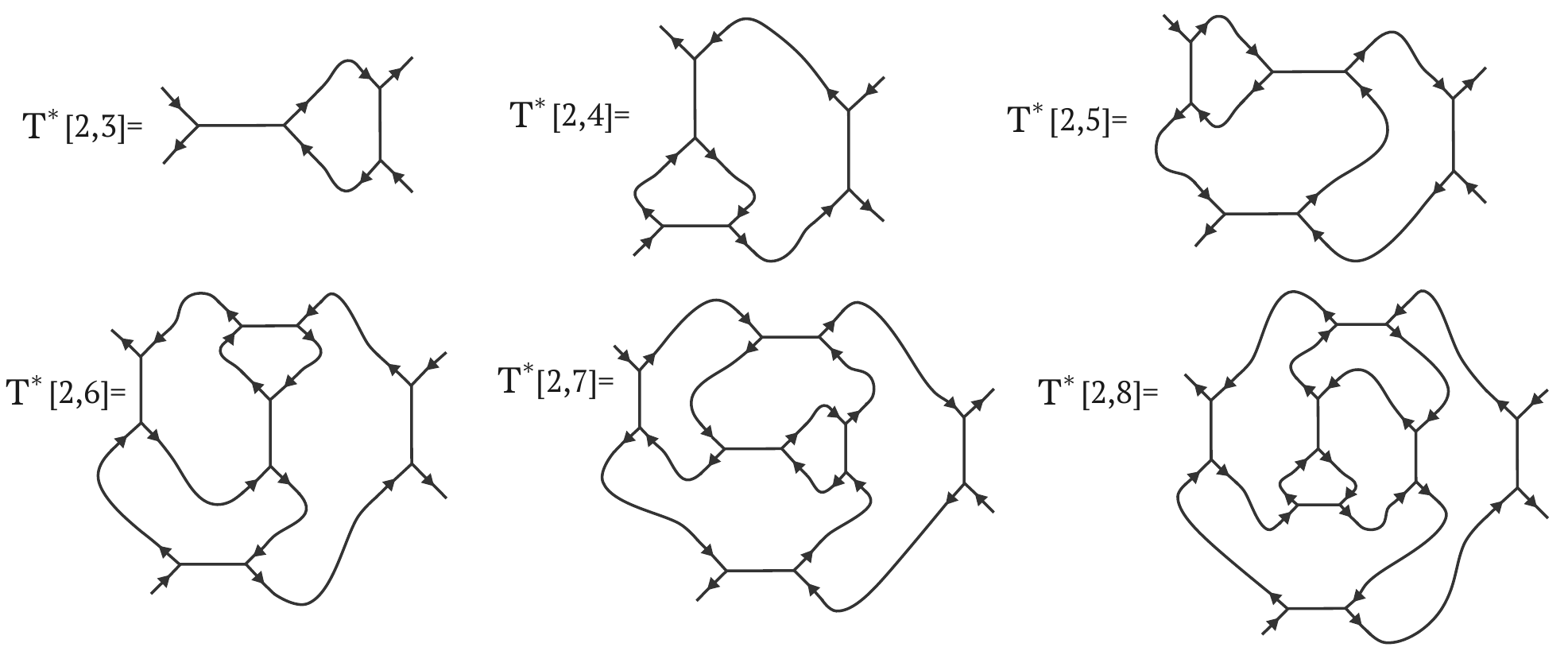}
	\caption{\footnotesize{Bipartite diagrams of tangles $T^*[2,p]$.}}
	\label{fig:BipTorTangle}
\end{figure}

We may proceed in this way until the original diagram is reduced to a combination of lock elements and the 4-tangle of the trefoil, whose bipartite diagram is known (see Fig.~\ref{fig:BipT3}). In this manner, we obtain a recipe for constructing a bipartite diagram of an arbitrary two-strand link (see examples in Fig.~\ref{fig:BipTorTangle}): we take a bipartite diagram of a two-strand link corresponding to a tangle $T^*[2,p]$, rotate it by $\frac{\pi}{2}$ clockwise and glue it on the right to a lock element which yields the bipartite diagram of the tangle $T^*[2,p+1]$.

\subsection{Remark about twist knots}\label{sec:TwistKnots}

\begin{wrapfigure}[10]{r}{0.3\textwidth}
	\centering
	\vspace{-1.0cm}
	\includegraphics[width =1.0\linewidth]{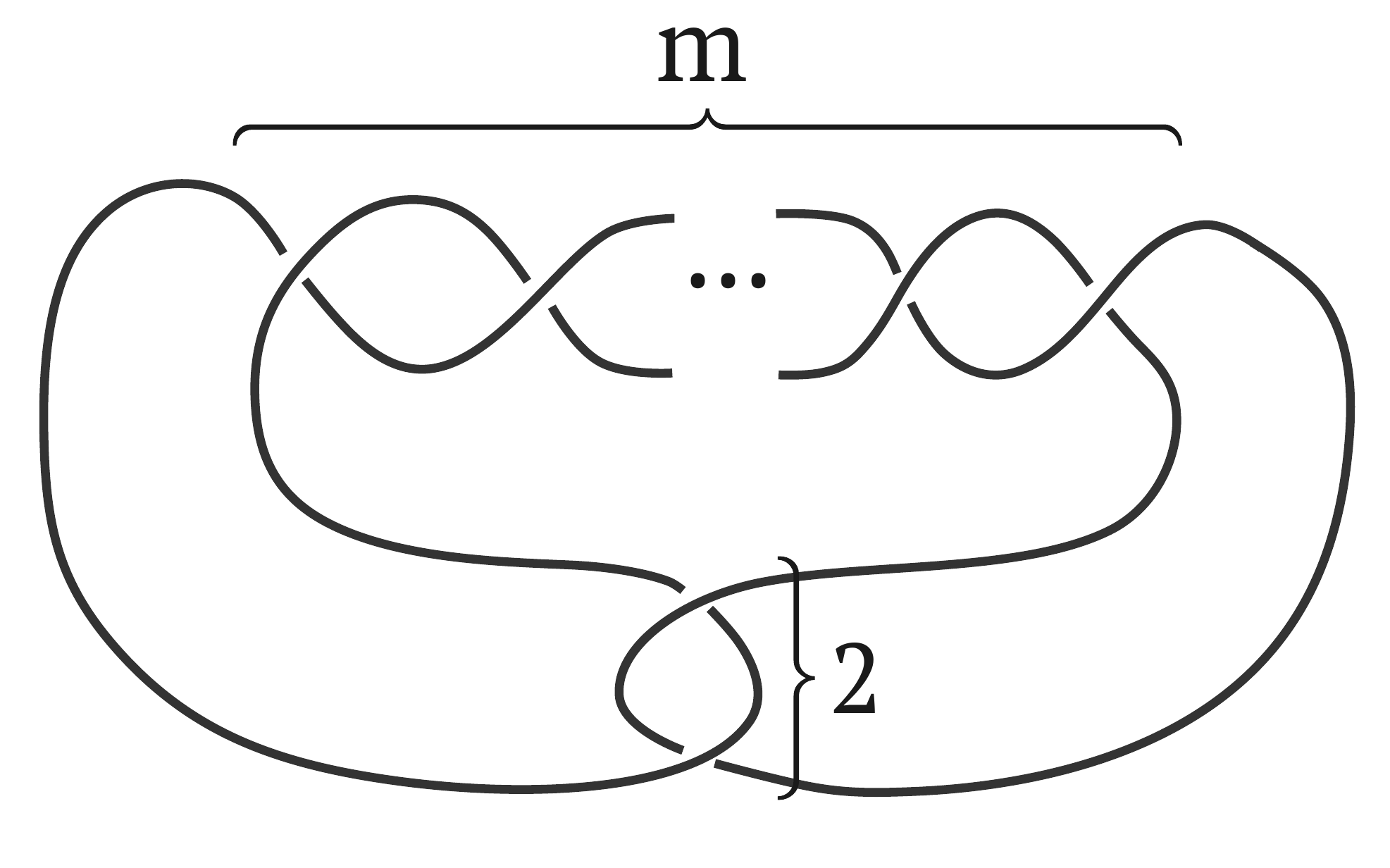}
	\caption{\footnotesize{Example of a twist knot Tw$_m$ with $m+2$ crossings.}}
	\label{fig:TwistExample}
\end{wrapfigure}

Twist knots ${\rm Tw}_{m}$ form a one-parameter family of knots with $m + 2$ crossings, where $m$ and $2$ correspond to the usual two-strand tangles, as illustrated in Fig.~\ref{fig:TwistExample}.

If $m = 2k$ is even, then the diagram naturally decomposes into lock elements, whose orientations are antiparallel by construction. Thus, the diagram is immediately bipartite. If $m = 2k + 1$ is odd, the upper tangle with $m$ crossings also has an antiparallel orientation while the lower one becomes parallel. In this case we need to reverse the orientation of the lower tangle and add one crossing to the upper one. Figure~\ref{fig:TwistToBip} demonstrates that this can be achieved through simple diagrammatic transformations:

\begin{figure}[h!]
	\centering	
	\includegraphics[width =0.85\linewidth]{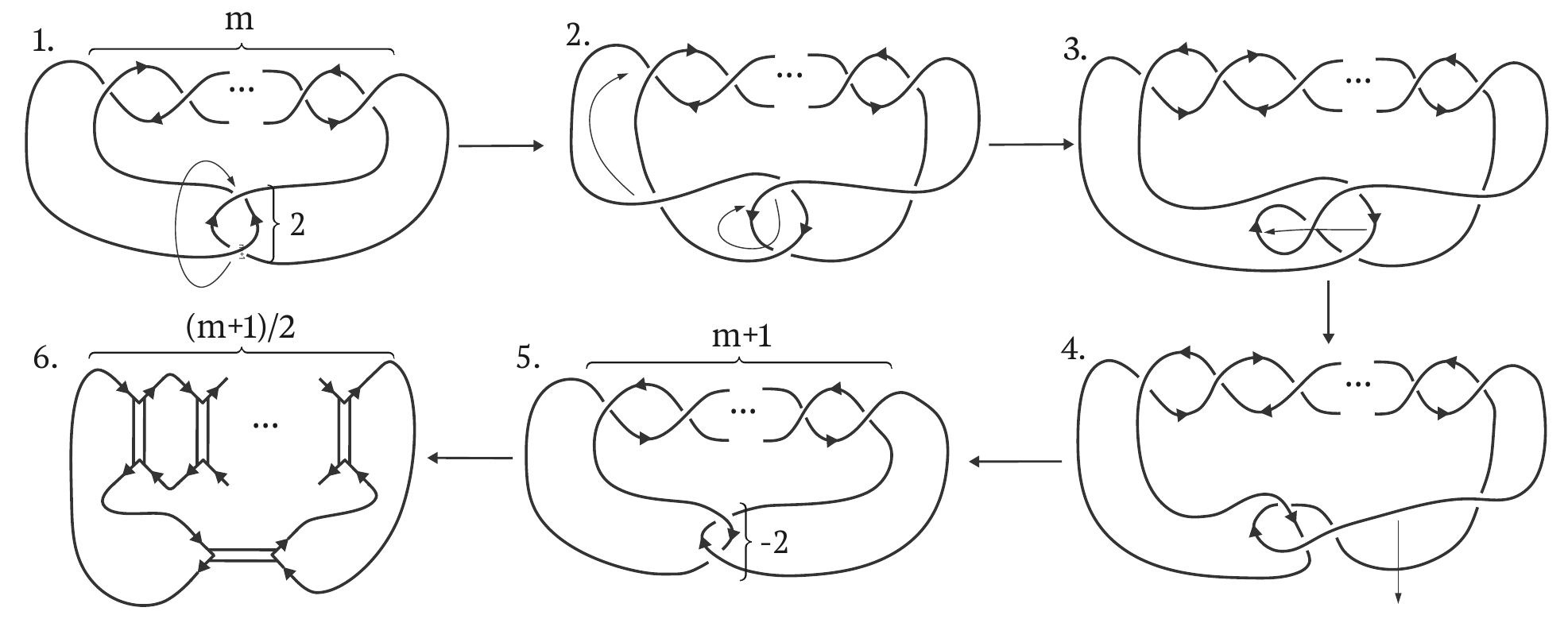}
	\caption{\footnotesize{Bipartization of ${\rm Tw}_m$ with odd $m$.}}
	\label{fig:TwistToBip}
\end{figure}

In the first step $(1 \to 2)$, we twist the lower tangle by $\pi$, which results in two crossings of opposite signs appearing on its sides. In the next step $(2 \to 3)$, one of these new crossings is moved into the upper tangle so that the crossings in the upper tangle are grouped into $\frac{m+1}{2}$ lock elements. In the following steps $(2 \to 4)$, we use the remaining crossing to reverse the orientation of the lower tangle similarly to the orientation change in the Hopf link (see Fig.~\ref{fig:BipHopf}). As a result, we obtain a twist knot with an even number $m = 2k + 2$ of crossings in the upper tangle and with the lower tangle being the mirror image of the initial one.

\subsection{Bipartiteness of pretzel knots
}\label{sec:BipForPretzelKnots}

\begin{wrapfigure}[10]{r}{0.225\textwidth}
	\centering
	\vspace{-2.2cm}
	\includegraphics[width =1.0\linewidth]{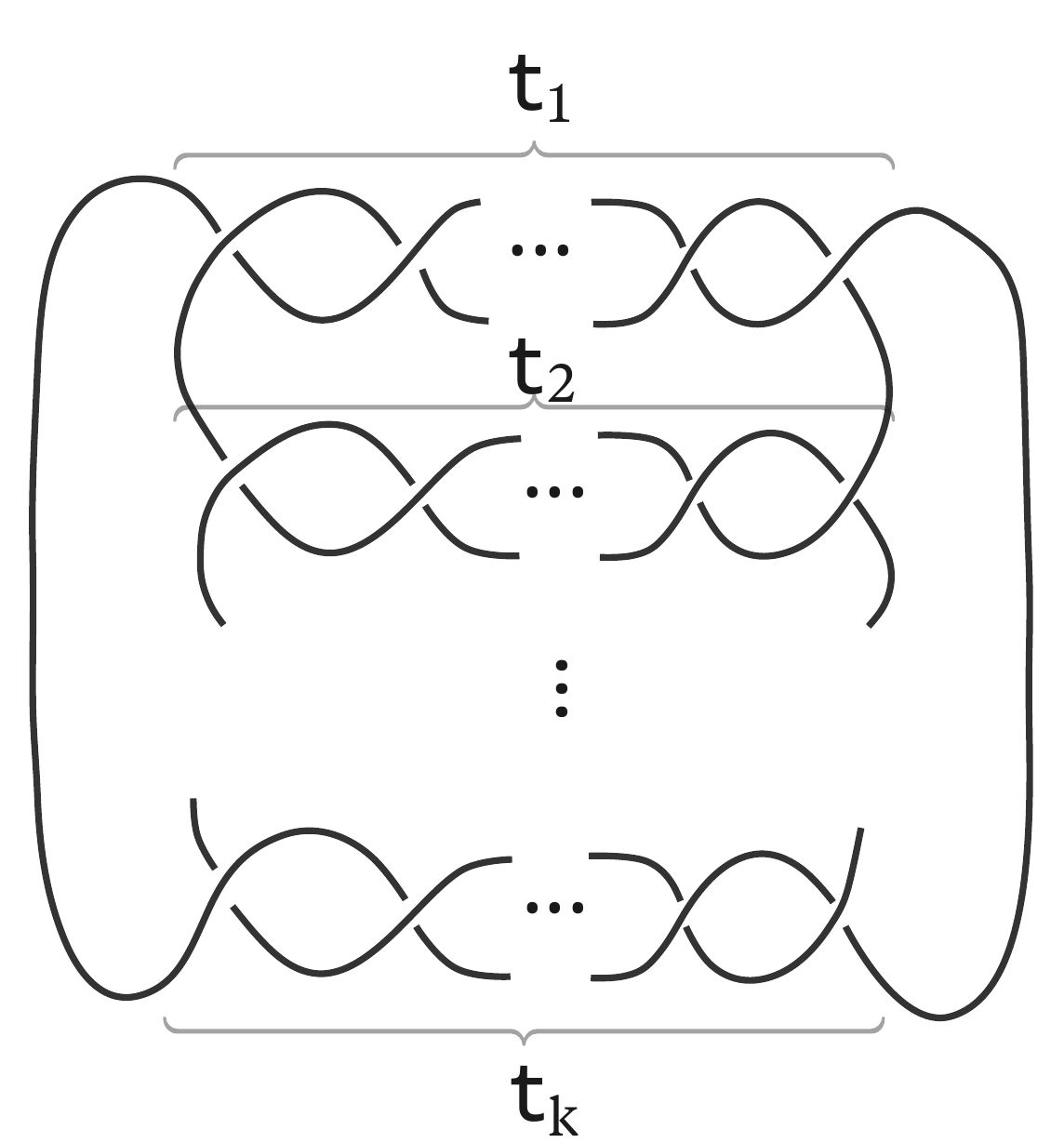}
	\caption{\footnotesize{Pretzel link $P(t_1, \dots, t_k)$ consisting of $k$ tangles.}}
	\label{fig:PretzelExample}
\end{wrapfigure}

Pretzel knots are the generalization of two-strand links, in which the usual crossings are replaced by two-strand tangles inserted perpendicularly to the initial braid (see Fig.~\ref{fig:PretzelExample}). This class is convenient for analysis and, at the same time, more general than the two-strand or twist families (the latter two are subsets of the pretzel family). 

It is of particular interest that for pretzel knots, there exist constraints on the parameters that guarantee the presence or absence of bipartiteness in their diagrammatic representations. Obviously, the most straightforward condition ensuring bipartiteness is the evenness of all $t_i$. However, we are interested in bipartiteness of more advanced configurations. The following theorem serves as an example.

\begin{theorem}\label{3Pretzel}{\bf {\cite{lewark2016new}}}
	If all the numbers of crossings in the two-strand tangles $t_1, t_2, t_3$ for a pretzel link $P(t_1, t_2, t_3)$ are odd, and the greatest common divisor satisfies $\gcd(t_1, t_2, t_3) > 1$, then the knot is non-bipartite. If at least one of the numbers $t_1, t_2, t_3$ is even, the link admits a bipartite realization.
\end{theorem}

\begin{figure}[h!]
	\centering	
	\includegraphics[width =0.7\linewidth]{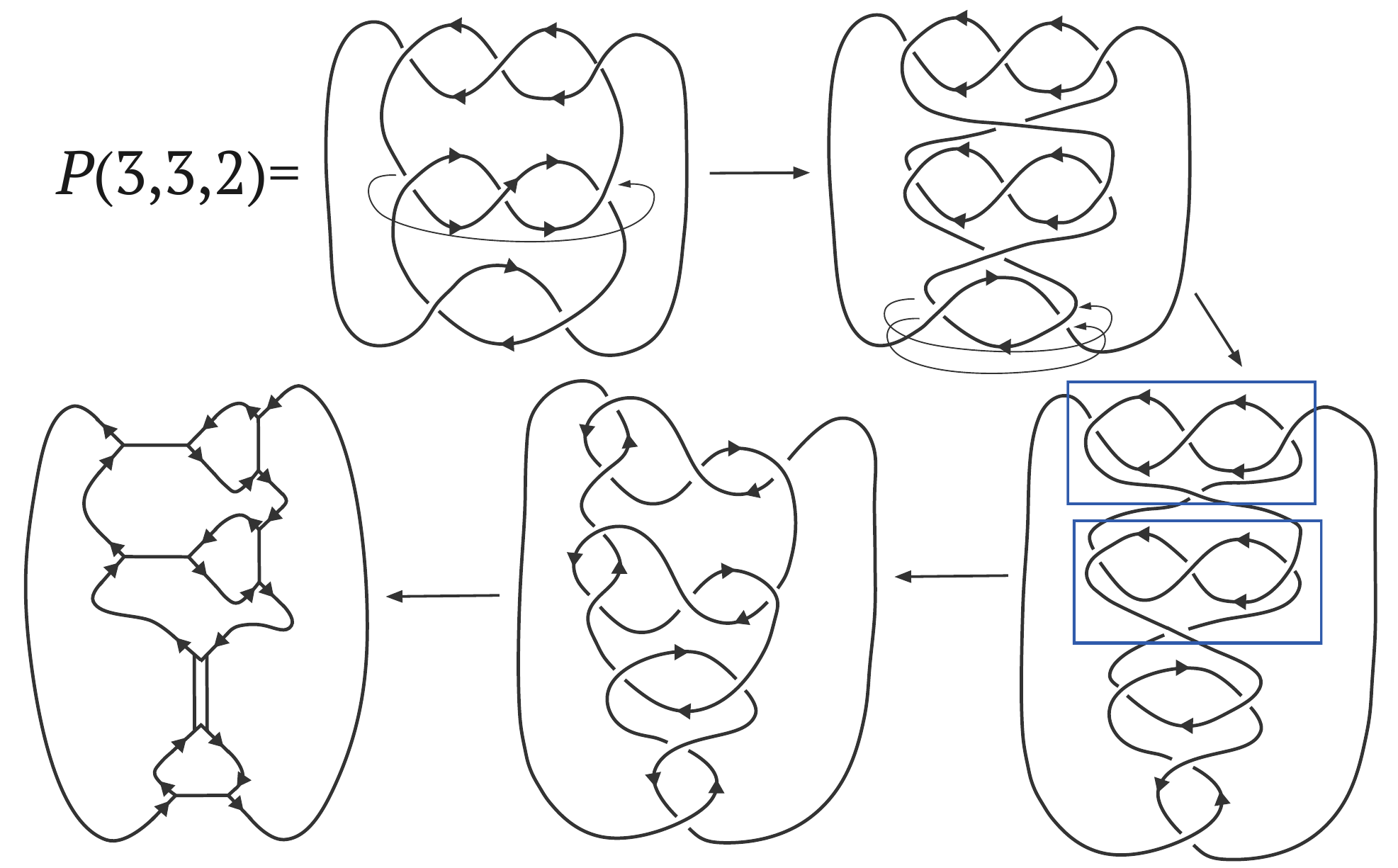}
	\caption{\footnotesize{Bipartization for knot $8_{5}$.}}
	\label{fig:Knot85ToBip}
\end{figure}

An illustration of Theorem~\ref{3Pretzel} can be provided by the previously mentioned knot $9_{46}$, which is the pretzel knot $P(3,3,-3)$ and, according to the theorem, is non-bipartite. As an opposite example, consider the knot $8_{5}$ (in pretzel notation, it is $P(3,3,2)$). This example is interesting not only because it admits a bipartite diagram, which we obtain in Fig.~\ref{fig:Knot85ToBip}, but also because its bipartite diagram contains lock elements of both types allowing us to consider, on its example, a more general quaternary Goeritz matrix.

\subsection{Rational knots}\label{sec:ratioknots}

\begin{figure}[h!]
	\centering	
	\includegraphics[width =0.9\linewidth]{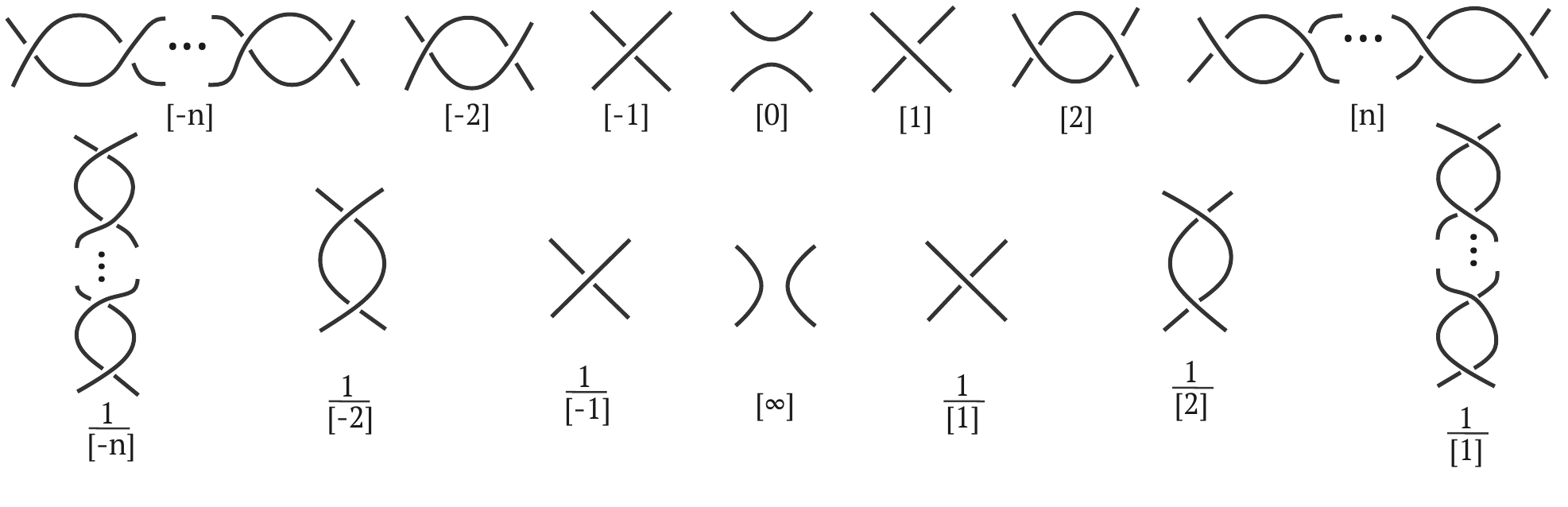}
	\caption{\footnotesize{The notation of composite tangles for rational knots.}}
	\label{fig:RTangle}
\end{figure}

As a family with bipartite tangles, the class of rational knots $R\!\left(\frac{p}{q}\right)$~\cite{kauffman2004classification,goldman1997rational} is particularly illustrative, since all of its members are bipartite~\cite{duzhin2010formula}. Rational links consist of two–strand tangles, are parameterized by a rational fraction $\frac{p}{q}$, and are constructed using its continued fraction expansion. Let us introduce a notation for denoting the two-strand tangles that compose $R\!\left(\frac{p}{q}\right)$, see Fig.~\ref{fig:RTangle}.

\begin{wrapfigure}[9]{r}{0.29\textwidth}
	\centering
	\vspace{-0.75cm}
	\includegraphics[width =0.75\linewidth]{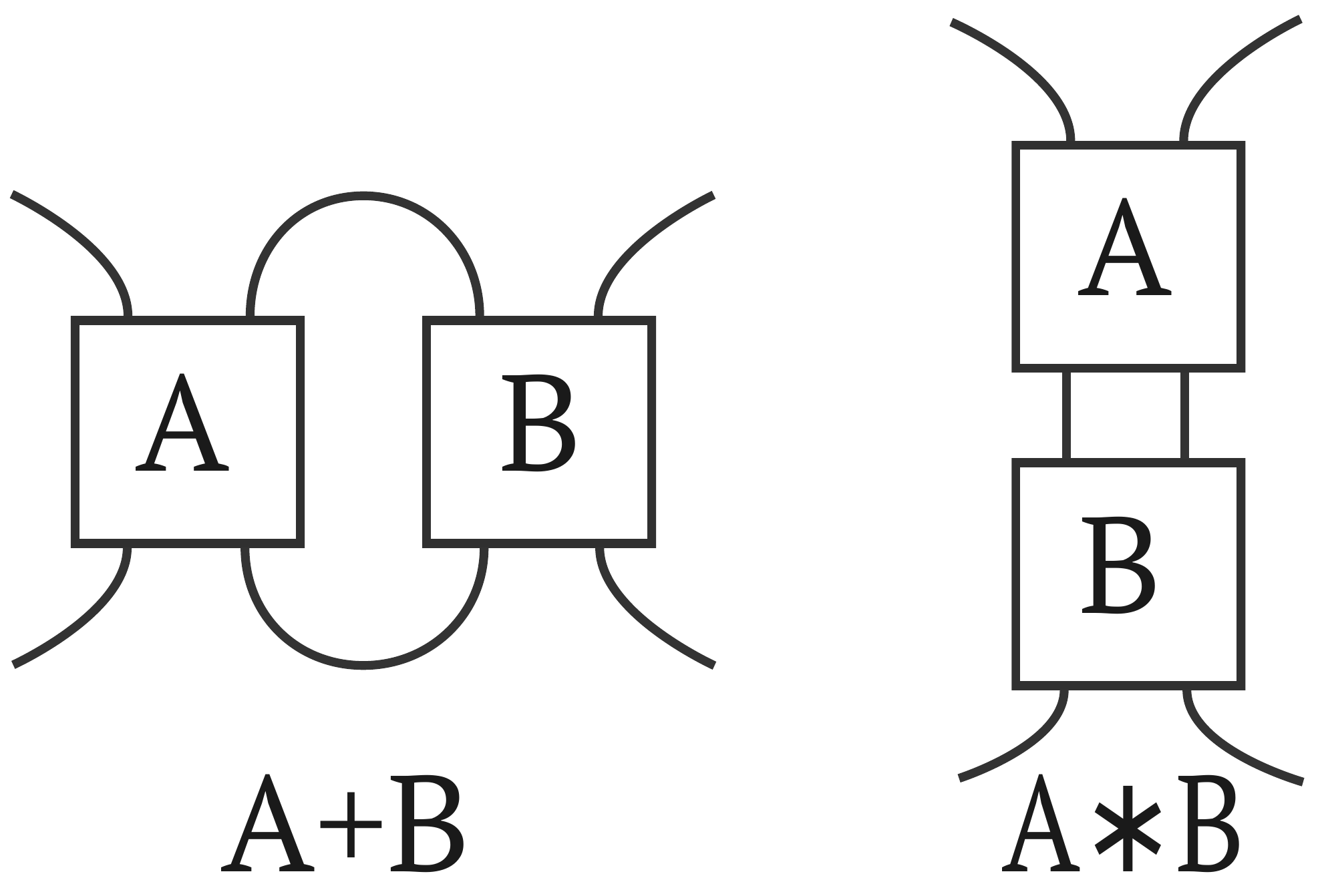}
	\caption{\scriptsize{The addition and multiplication operations on rational knots.}}
	\label{fig:RTangle2}
\end{wrapfigure}

We denote horizontal tangles by a number $[n]$, and vertical tangles by a fraction $\frac{1}{[n]}$. We require that these numbers encode the information about crossings in two-strand tangles: crossings shown in Fig.\,\ref{fig:RTangle} in the left columns are counted with a negative sign, while their mirror crossings in the right columns are counted with a positive sign. We also introduce the tangles $[0]$ and $[\infty]$ which represent horizontal and vertical parallel strands, respectively.

\begin{wrapfigure}[13]{r}{0.3\textwidth}
	\centering
	\vspace{-0.5cm}
	\includegraphics[width =1.0\linewidth]{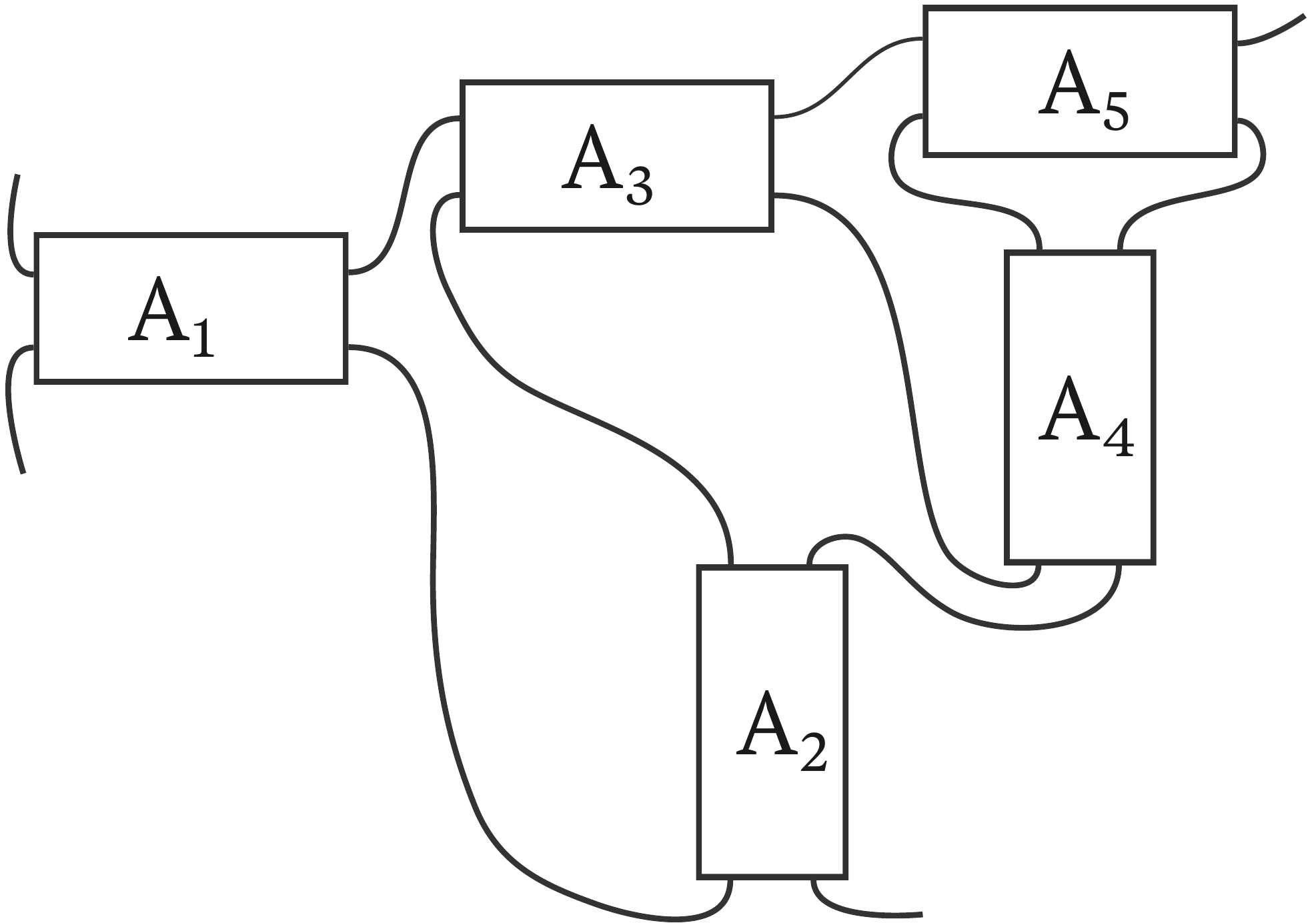}
	\caption{\footnotesize{Example of a rational tangle $R^*\left(\frac{p}{q}\right)$in standart.}}
	\label{fig:RatKnotA5}
\end{wrapfigure}

Rational knots can be represented in the standard form as in the example shown in Fig.\,\ref{fig:RatKnotA5}. This form is defined as follows. We start with a tangle $[0]$ or $[\infty]$. We introduce the operations of addition ($+$) and multiplication ($*$), see Fig.~\ref{fig:RTangle2}. Addition attaches a tangle to the right or left, while multiplication attaches a tangle to the top or bottom:

\begin{equation}
    [A] + [B] = [A+B]\,, \ \ \ \frac{1}{[A]}*\frac{1}{[B]} = \frac{1}{[A+B]}\,.
\end{equation}
By applying these operations in sequence we construct a rational tangle. In the standard form, we restrict ourselves to attaching crossing tangles only to the bottom-left\footnote{Alternatively, one can use only the top-right variant.}. Thus, every rational tangle consists of an alternation of horizontal and vertical two-strand tangles. The resulting diagram can be encoded by an ordered list of the crossing numbers  $A_i$ in such two-strand tangles as $[A_1,A_2,\dots,A_n]$:

\begin{equation}\label{GeneralRatTangle}
    [A_1,A_2,\dots,A_n]=[A_1] + \left(\frac{1}{[A_2]}*\dots*\left([A_{n-2}]+\left(\frac{1}{[A_{n-1}]}*[A_n]\right)\right)\right)\,. 
\end{equation}

\noindent This expression has a representation in the form of a continued fraction:
\begin{equation}\label{RationalFraction}
    \frac{p}{q} = A_1 + \cfrac{1}{A_2 + \cfrac{1}{\cdots + \cfrac{1}{A_{n-1} + \cfrac{1}{A_n}}}}\,.
\end{equation}
Due to the integer coefficients $A_i$, the entire fraction $\frac{p}{q}$ is rational, and it uniquely determines the tangle\footnote{In the notation of rational tangles used here, one must close the right and left ends, respectively, to obtain links.}. However, the continued fraction corresponding to $\frac{p}{q}$ is not unique.

\begin{theorem}\label{RatTheorem}{\bf {\cite{duzhin2010formula}}}
	If either $p$ or $q$ is even, then the fraction $\frac{p}{q} = [A_1,A_2,\dots,A_k]$ has a unique representation with all coefficients $A_i$ being even. In this case, the rational tangle $R\left(\frac{p}{q}\right)$ can be written in a bipartite form.
\end{theorem}

Rational knots with odd $p$ and $q$ also admit bipartite diagrams, since links $R\left(\frac{p}{q}\right)$ and $R\left(\frac{p'}{q}\right)$ $(p' = p - q, \, p>0,$ and $\, p' = p + q, \, p<0)$ are isotopic to each other \cite{duzhin2010formula}. Hence, by Theorem \ref{RatTheorem}, $R\left(\frac{p'}{q}\right)$, and consequently $R\left(\frac{p}{q}\right)$, are bipartite. Thus, we have the following statement.

\begin{theorem}\label{RatKnotTheorem}{\bf {\cite{duzhin2010formula}}}
    All rational knots are bipartite.
\end{theorem}

\begin{figure}[h!]
		\centering	
		\includegraphics[width =0.75\linewidth]{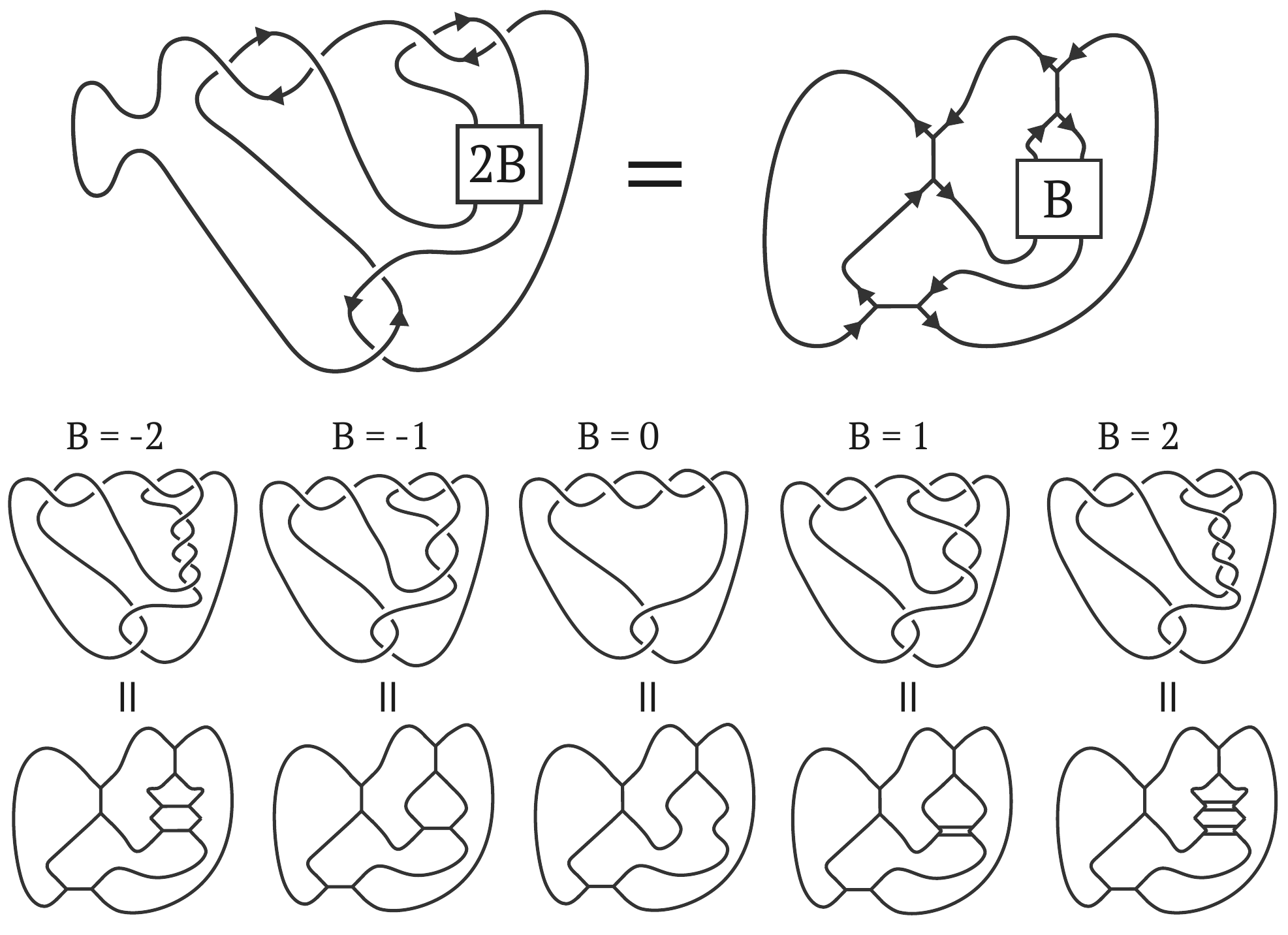}
        \caption{\footnotesize{Examples of knots from a one-parametric family of rational knots.}}
        \label{fig:RatKnotBExample}
	\end{figure}

As an example of writing a bipartite rational knot by a rational fraction let us consider the following one-parameter family of continued fractions which is given by a continuous fraction:

\begin{equation}\label{RatTangleExample}
    R\left(\frac{p}{q}\right) =  \frac{8B+4}{-12B-7} = \cfrac{1}{-2 + \cfrac{1}{2 + \cfrac{1}{2B + \cfrac{1}{2}}}}\,.
\end{equation}
We distinguish the positive values of the parameter $B$, since in this case, both types of lock elements appear in diagrams (see Fig.\,\ref{fig:RatKnotBExample}). The fraction~\eqref{RatTangleExample} corresponds to a tangle consisting of five elements: $A_1 = 0,\ A_2 = -\frac{1}{2},\ A_3 = 2,\ A_4 = \frac{1}{2B},\ A_5 = 2$, as in the example shown in Fig.\,\ref{fig:RatKnotA5}. For small values of the parameter $B \in \{-2, -1, 0, 1, 2\}$, we provide the correspondence with knots from Rolfsen table~\cite{rolfsen2003knots}, which can be established by computing the HOMFLY–PT polynomials from the bipartite diagrams (see Section~\ref{sec:HOMFLYForRatioKnotsFamily}). For accuracy, we note that the results obtained using our formula~\eqref{Hrelation} correspond to the mirror images of the knots listed in Rolfsen table.

\begin{table}[h]
	\centering
	\begin{tabular}{|c|c|c|c|c|c|}
		\hline
		$B$& $-2$ & $-1$ & 0 & 1 & 2 \\
		\hline
		${\text{Knot}}$  $R\left(\frac{p}{q}\right)$ & $ \overline{7}_5$ & $ \overline{5}_1$ &$\overline{5}_2$&$\overline{7}_6$ & $\overline{9}_8$  \\ 
		\hline
		\raisebox{-0.1cm}{$\frac{p}{q}$} & \raisebox{-0.1cm}{$-\frac{12}{17}$} & \raisebox{-0.1cm}{$-\frac{4}{5}$}  & \raisebox{-0.1cm}{$-\frac{4}{7}$} & \raisebox{-0.1cm}{$-\frac{12}{19}$} &  \raisebox{-0.1cm}{$-\frac{20}{31}$} \\ [1.5ex]
		\hline
	
	\end{tabular}
	\caption{\footnotesize{Some rational knots from the family~\eqref{RatTangleExample}. Their bipartite diagrams are shown in Fig.\,\ref{fig:RatKnotBExample}.}}	
	\label{tab:char}
\end{table}

\newpage

\subsection{Montesinos knots}

\begin{wrapfigure}[10]{r}{0.35\textwidth}
	\centering
	\vspace{-0.0cm}
	\includegraphics[width =\linewidth]{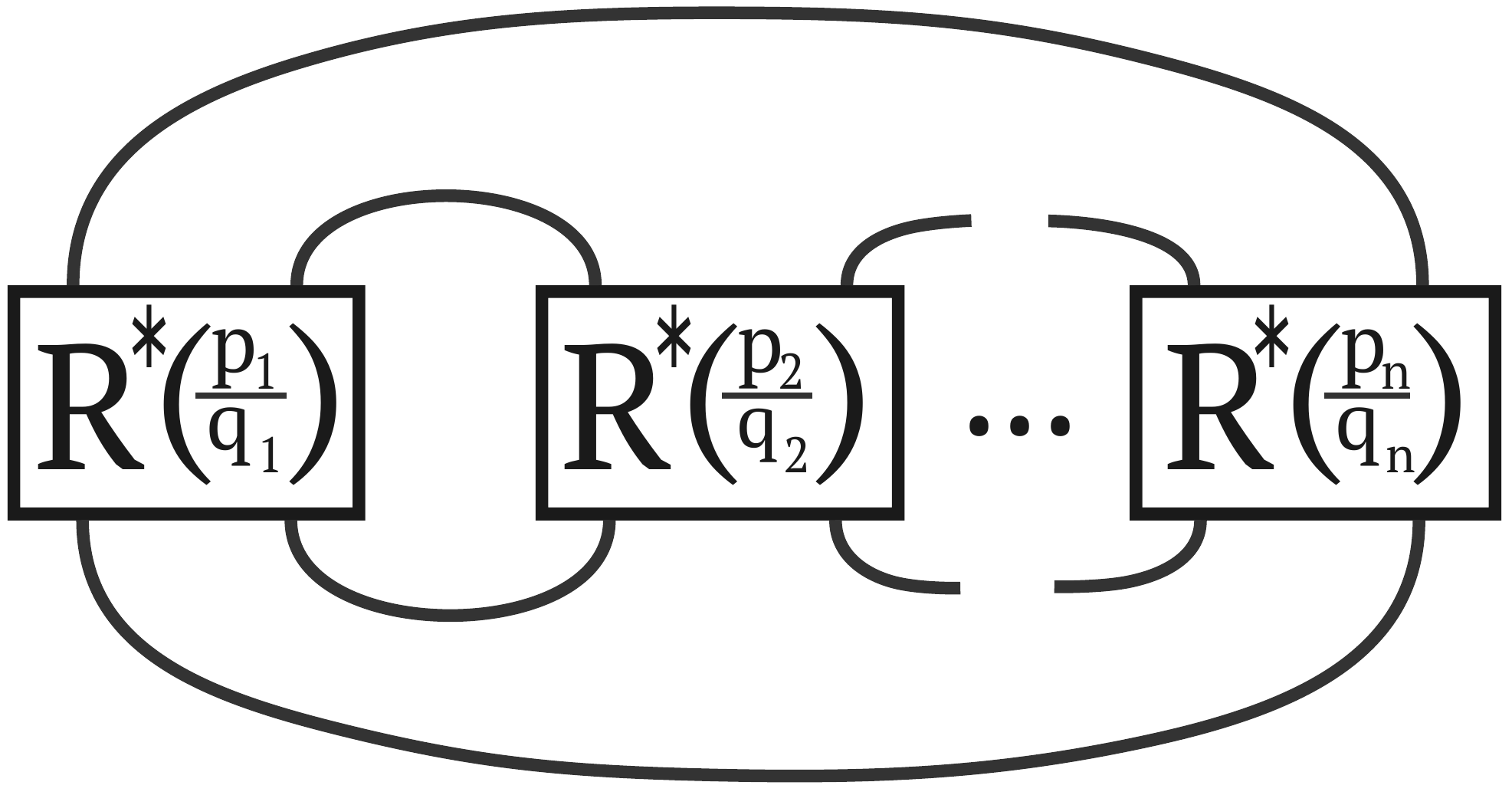}
	\caption{\footnotesize{Montesinos knot $K\left(\frac{p_1}{q_1}, \dots , \frac{p_n}{q_n}\right)$.}}
	\label{fig:Montesinos}
\end{wrapfigure}

\begin{figure}[h!]
		\centering	
		\includegraphics[width =0.75\linewidth]{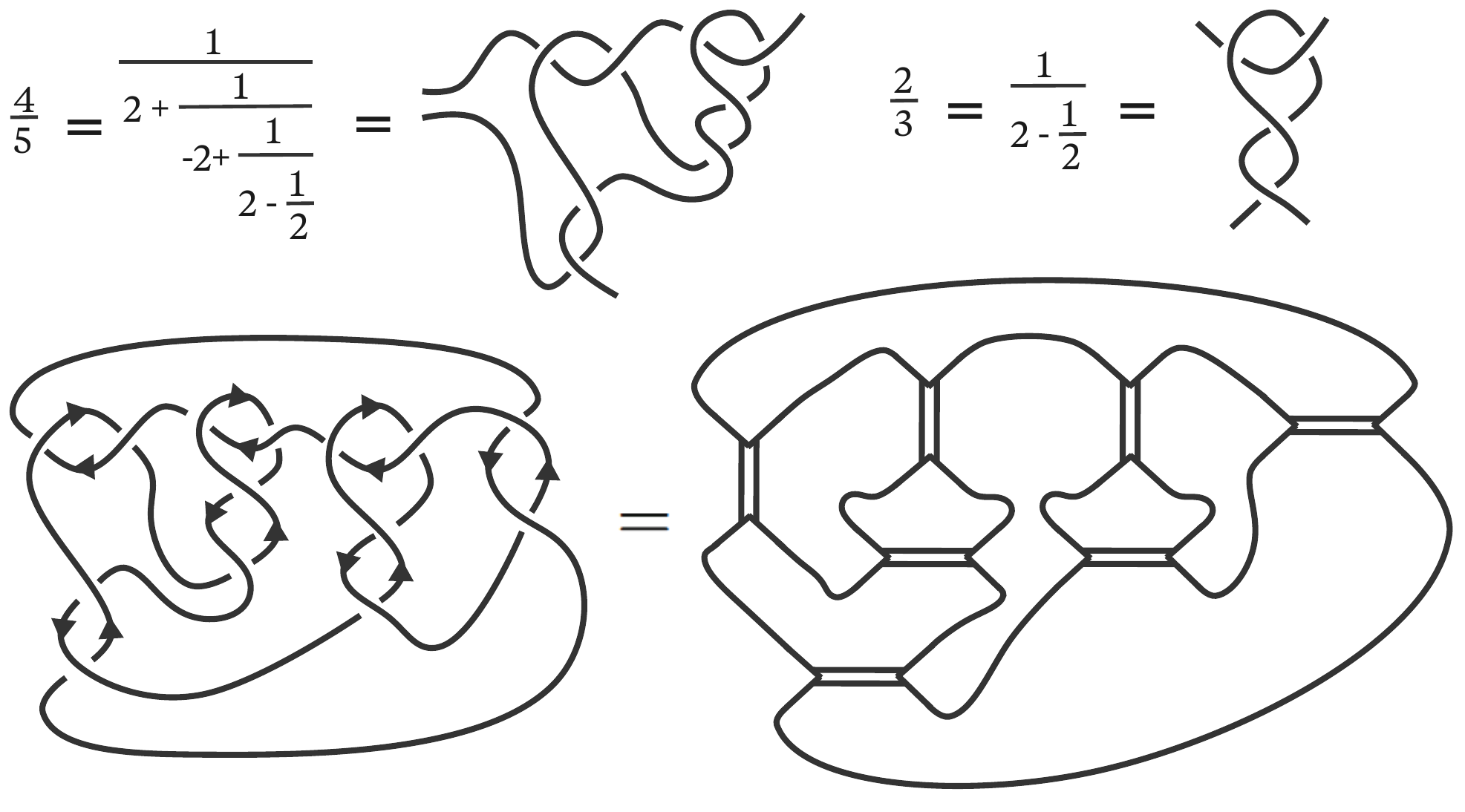}
        \caption{\footnotesize{Composite tangles and diagrams for the bipartite Montesinos knot $10_{49}$.}}
        \label{fig:MontesinosKnot1049}
	\end{figure}

There exists another family of links that is particularly interesting from the viewpoint of bipartite representation, which generalize pretzel links — the Montesinos links $K\left(\frac{p_1}{q_1}, \dots , \frac{p_n}{q_n}\right)$~\cite{lewark2016new}. In Montesinos links, rational tangles are placed instead of two-strand tangles used in pretzel links (see Fig.\,\ref{fig:Montesinos}). A list of knots from Rolfsen table in Montesinos notation can be found in~\cite{dunfield2001table}. In an obvious way, pretzel links $P(t_1, t_2, \dots, t_n)$ can be rewritten in Montesinos notation as $K\left(\frac{1}{t_1}, \frac{1}{t_2}, \dots, \frac{1}{t_n}\right)$; the previously discussed knot $8_5$ has the form $K\left(\frac{1}{3},\,\frac{1}{3},\,\frac{1}{2}\right)$ in this notation.\footnote{For Montesinos knots we use a closure rotated by $\frac{\pi}{2}$ relative to the closure of pretzel knots (see Fig.\,\ref{fig:PretzelExample}), while the orientation of the tangles remains unchanged. For this reason, the relabeling $t_i \to \frac{1}{t_i}$ occurs when switching from pretzel to Montesinos notation in order to keep consistency with the conventions used in~\cite{lewark2016new, dunfield2001table}.}  
As a less trivial example, we can consider the knot $10_{49} = K\left(\frac{4}{5},\,\frac{2}{3},\,\frac{1}{2}\right)$ shown in Fig.\,\ref{fig:MontesinosKnot1049}. There also exist parameter constraints that ensure the existence of a bipartite diagram for Montesinos links.

\begin{theorem}\label{MontesinosTheorem}{\bf {\cite{lewark2016new}}}
If a link $K\left(\frac{p_1}{q_1}, \dots , \frac{p_k}{q_k}\right)$ has several components, then this link is bipartite. If $K\left(\frac{p_1}{q_1}, \dots , \frac{p_k}{q_k}\right)$ is a knot and at least one of the denominators $q_i$ is even, then this knot is bipartite.
\end{theorem}

The knot $10_{49}$ satisfies the condition of the theorem, since one of the fractions has the denominator $2$. Therefore, $10_{49}$ admits a bipartite diagram, which can be explicitly demonstrated by expressing the rational fractions that parametrize the knot as continued fractions with all even coefficients (see Fig.\,\ref{fig:MontesinosKnot1049}).

\section{Goeritz matrix}
\label{sec:Goeritz-mat}
\setcounter{equation}{0}
The Goeritz matrix is an object constructed from the checkerboard coloring of a link diagram, carrying both topological and combinatorial information about a knot diagram. Due to this, the Goeritz matrix gives rise to several invariants at once. Among them, there is the Jones polynomial, which we will discuss in detail in Section~\ref{sec:JonesFromGoeritzMatrix}. In addition, the absolute value of the determinant of the Goeritz matrix yields the knot determinant ${\rm Det}^{K}$ discovered simultaneously with the Goeritz matrix itself~\cite{goeritz1933knoten}.  
The topological nature of the Goeritz matrix is also revealed through the Gordon–Litherland form~\cite{gordon1978signature}, which allows one to compute the knot signature, often in a more convenient way than through the Seifert matrix, since the Goeritz matrix typically has a smaller size and a simpler construction algorithm, which we present below.
\begin{enumerate}

\item To each knot diagram we associate a surface consisting of two types of regions (white and black\footnote{In the figures in this paper, the regions are shaded in gray so that the contour of the knot remains visible.} — shaded according to the rule of checkerboard coloring). All white regions are then numbered.

\item We define the symmetric unreduced Goeritz matrix $\widetilde{G} = (\widetilde{g}_{ij})$ for a given surface according to the following rule:  

\begin{equation}\label{IrGoerMatr}
\centering
\widetilde{g}_{ij} = \left\{\begin{split}
-\sum_{c_{i,j}}\sigma(c_{i,j}) &, \ \ \ \   i\ne j\,, \\
-\sum_{k \ne i} \widetilde{g}_{kj}  &, \ \ \ \  i = j\,, \\ 
\end{split}\right.
\end{equation}
where $c_{i,j}$ are crossings between regions $i$ and $j$. Depending on a position of shaded regions, each crossing $c$ contributes the value $\sigma(c) = \pm 1$ to corresponding matrix elements (see Fig.~\ref{fig:sign-convention}).

\begin{figure}[h!]
\begin{picture}(300,60)(-180,-15)

\thicklines

\put(0,0){\line(1,1){40}}
\put(0,40){\line(1,-1){17}}
\put(23,17){\line(1,-1){17}}

\put(0,41){\tikz \coordinate (point1);}
\put(40,41){\tikz \coordinate (point2);}
\put(20,21){\tikz \coordinate (point3);}
\put(0,0){\tikz[overlay] \fill[bright-gray] (point1) -- (point2) -- (point3) -- cycle;}

\put(0,-1){\tikz \coordinate (point4);}
\put(40,-1){\tikz \coordinate (point5);}
\put(20,19){\tikz \coordinate (point6);}
\put(0,0){\tikz[overlay] \fill[bright-gray] (point4) -- (point5) -- (point6) -- cycle;}

\put(0,-15){\mbox{$\sigma(c)=+1$}}

\put(100,40){\line(1,-1){40}}
\put(100,0){\line(1,1){17}}
\put(123,23){\line(1,1){17}}

\put(100,0){

\put(0,41){\tikz \coordinate (point1);}
\put(40,41){\tikz \coordinate (point2);}
\put(20,21){\tikz \coordinate (point3);}
\put(0,0){\tikz[overlay] \fill[bright-gray] (point1) -- (point2) -- (point3) -- cycle;}

\put(0,-1){\tikz \coordinate (point4);}
\put(40,-1){\tikz \coordinate (point5);}
\put(20,19){\tikz \coordinate (point6);}
\put(0,0){\tikz[overlay] \fill[bright-gray] (point4) -- (point5) -- (point6) -- cycle;}

\put(0,-15){\mbox{$\sigma(c)=-1$}}
}
         
\end{picture}
    \caption{\footnotesize Checkerboard sign convention.}
    \label{fig:sign-convention}
\end{figure}

\item To obtain the (reduced) Goeritz matrix $G = g_{ij}$, we remove the $k$-th row and column $(\forall k)$. The remaining elements then form the desired matrix $g_{ij} = \widetilde{g}_{ij},\; \forall\, i\ne k,\, j\ne k$. Thus, if an initial unreduced Goeritz matrix is of size $n \times n$, a reduced matrix is of size $(n-1) \times (n-1)$, where $n$ is a number of unshaded (white) regions.

\end{enumerate}

\begin{figure}[h!]
		\centering	
		\includegraphics[width =0.7\linewidth]{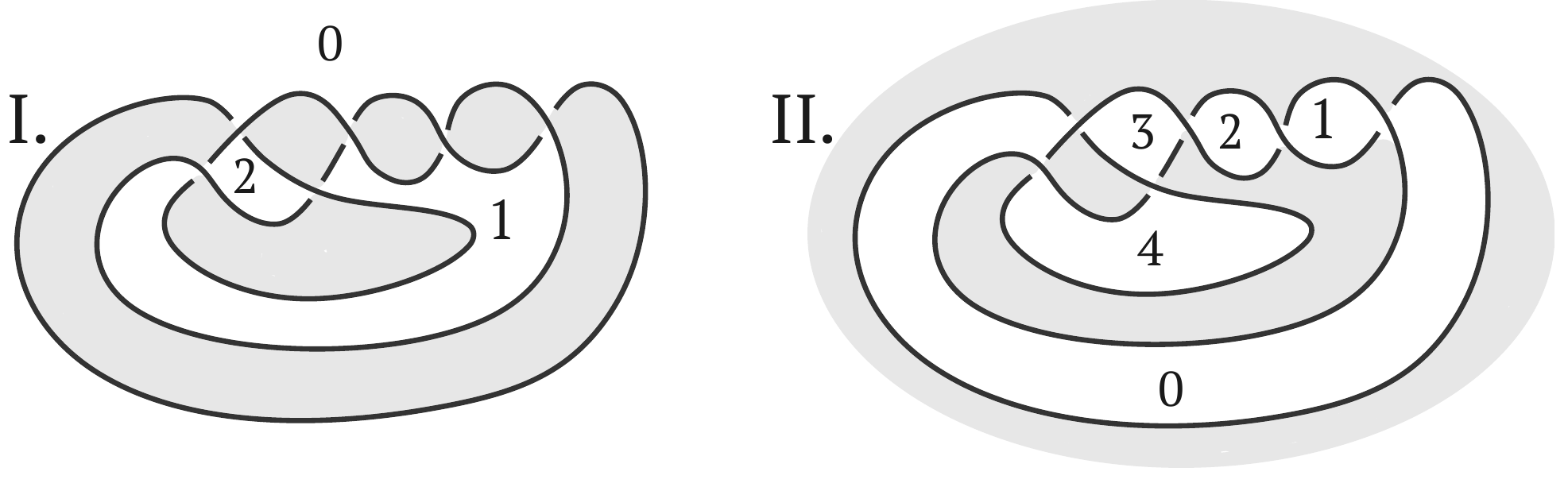}
        \caption{\footnotesize The two checkerboard colorings of the knot $5_2$.}\label{pic:Knot52}
	\end{figure}

As an example we present the Goeritz matrices and checkerboard colorings for the trefoil knot $3_1$ and the knot $5_2$ (see Figs.~\ref{fig:trefoil}--\ref{pic:Knot52}).

\begin{figure}[h!]
    \centering
\begin{tikzpicture}[scale=0.5]

\node at (0.5, 1.8) {1};
\node at (0.5, -0.5) {2};
\node at (0.5, -2.5) {3};
\node at (0.5, -4.8) {1};

\fill[bright-gray] (0,0) -- (0,1) -- (0.5,0.5) -- cycle; 
\fill[bright-gray] (1,1) -- (1,0) -- (0.5,0.5) -- cycle; 

\fill[bright-gray] (0,-2) -- (0,-1) -- (0.5,-1.5) -- cycle; 
\fill[bright-gray] (1,-1) -- (1,-2) -- (0.5,-1.5) -- cycle; 

\fill[bright-gray] (0,-4) -- (0,-3) -- (0.5,-3.5) -- cycle; 
\fill[bright-gray] (1,-3) -- (1,-4) -- (0.5,-3.5) -- cycle; 

\draw[thick] (0,0) -- (1,1); 
\draw[thick] (1,0) -- (0.6,0.4);
\draw[thick] (0,1) -- (0.4,0.6);

\draw[thick] (0,-2) -- (1,-1);
\draw[thick] (1,-2) -- (0.6,-1.6);
\draw[thick] (0,-1) -- (0.4,-1.4);

\draw[thick] (0,-4) -- (1,-3);
\draw[thick] (1,-4) -- (0.6,-3.6);
\draw[thick] (0,-3) -- (0.4,-3.4);

\fill[bright-gray] (1,1) .. controls (3,3) and (3,-6) .. (1,-4);
\draw[thick] (1,1) .. controls (3,3) and (3,-6) .. (1,-4);
\fill[white] (1,0) .. controls (1.5,-0.5) .. (1,-1);
\draw[thick] (1,0) .. controls (1.5,-0.5) .. (1,-1);
\fill[white] (1,-2) .. controls (1.5,-2.5) .. (1,-3);
\draw[thick] (1,-2) .. controls (1.5,-2.5) .. (1,-3);

\fill[bright-gray] (0,1) .. controls (-2,3) and (-2,-6) .. (0,-4);
\draw[thick] (0,1) .. controls (-2,3) and (-2,-6) .. (0,-4);
\fill[white] (0,0) .. controls (-0.5,-0.5) .. (0,-1);
\draw[thick] (0,0) .. controls (-0.5,-0.5) .. (0,-1);
\fill[white] (0,-2) .. controls (-0.5,-2.5) .. (0,-3);
\draw[thick] (0,-2) .. controls (-0.5,-2.5) .. (0,-3);

\fill[bright-gray] (10.5,-1.5) circle (4);

\fill[white] (10,0) -- (10,1) -- (10.5,0.5) -- cycle; 
\fill[white] (11,1) -- (11,0) -- (10.5,0.5) -- cycle; 

\fill[white] (10,-2) -- (10,-1) -- (10.5,-1.5) -- cycle; 
\fill[white] (11,-1) -- (11,-2) -- (10.5,-1.5) -- cycle; 

\fill[white] (10,-4) -- (10,-3) -- (10.5,-3.5) -- cycle; 
\fill[white] (11,-3) -- (11,-4) -- (10.5,-3.5) -- cycle;

\draw[thick] (10,0) -- (11,1); 
\draw[thick] (11,0) -- (10.6,0.4);
\draw[thick] (10,1) -- (10.4,0.6);

\draw[thick] (10,-2) -- (11,-1);
\draw[thick] (11,-2) -- (10.6,-1.6);
\draw[thick] (10,-1) -- (10.4,-1.4);

\draw[thick] (10,-4) -- (11,-3);
\draw[thick] (11,-4) -- (10.6,-3.6);
\draw[thick] (10,-3) -- (10.4,-3.4);

\fill[white] (11,1) .. controls (13,3) and (13,-6) .. (11,-4);
\draw[thick] (11,1) .. controls (13,3) and (13,-6) .. (11,-4);

\fill[bright-gray] (11,0) .. controls (11.5,-0.5) .. (11,-1);
\draw[thick] (11,0) .. controls (11.5,-0.5) .. (11,-1);

\fill[bright-gray] (11,-2) .. controls (11.5,-2.5) .. (11,-3);
\draw[thick] (11,-2) .. controls (11.5,-2.5) .. (11,-3);

\fill[white] (10,1) .. controls (8,3) and (8,-6) .. (10,-4);
\draw[thick] (10,1) .. controls (8,3) and (8,-6) .. (10,-4);

\fill[bright-gray] (10,0) .. controls (9.5,-0.5) .. (10,-1);
\draw[thick] (10,0) .. controls (9.5,-0.5) .. (10,-1);

\fill[bright-gray] (10,-2) .. controls (9.5,-2.5) .. (10,-3);
\draw[thick] (10,-2) .. controls (9.5,-2.5) .. (10,-3);

\node at (9.5, -1.5) {1};
\node at (11.5, -1.5) {2};

\end{tikzpicture}
    \caption{\footnotesize The trefoil knot and its two checkerboard colorings.}
    \label{fig:trefoil}
\end{figure}

\noindent For the right-hand coloring of the trefoil the unreduced and reduced Goeritz matrices are as follows:
\begin{equation}
    \widetilde{G} = \begin{pmatrix}\label{G-mat-3-1-1}
        3 & -3 \\
        -3 & 3
    \end{pmatrix} \ \ \mapsto \ \ G = 3\,.
\end{equation}

\noindent For the left-hand coloring:
\begin{equation}\label{G-mat-3-1-2}
    \widetilde{G} = \begin{pmatrix}
        -2 & 1 & 1 \\
        1 & -2 & 1 \\
        1 & 1 & -2
    \end{pmatrix}\ \ \mapsto \ \ G = \begin{pmatrix}
        -2 & 1 \\
        1 & -2
    \end{pmatrix}\,.
\end{equation}
Note that the determinant of the Goeritz matrix for the trefoil knot $3_1$ is the same for both colorings: $\det(G_1^{3_1}) = 3$.

While for the Goeritz matrices of the trefoil the reduction does not depend on the choice of the removed row and column, for the knot $5_2$ and both of its colorings (Fig.~\ref{pic:Knot52}), there are several distinct reduced Goeritz matrices. For the left-hand coloring:

\begin{equation}\label{GoeritzForKnot52N1}
\widetilde{G}{}^{5_2}_{I}  =
\begin{pmatrix}
 \begin{array}{ccc}  2&-3&1  \\  -3&1&2\\1&2&-3       
 \end{array}
 \end{pmatrix} \ \ \mapsto \ \ G_1^{5_2}  =
\begin{pmatrix}
 \begin{array}{cc}   1&2\\2&-3       
 \end{array}
 \end{pmatrix}, \ \ G_2^{5_2}  =
\begin{pmatrix}
 \begin{array}{cc}   2&-3\\-3&1       
 \end{array}
 \end{pmatrix}  \ \ \text{or} \ \ G_3^{5_2}  =
\begin{pmatrix}
 \begin{array}{cc}   2&1\\1&-3       
 \end{array}
 \end{pmatrix}.
\end{equation}
For the right-hand coloring the number of white regions and, consequently, the number of reduced matrices is larger, so we present only one Goeritz matrix:

\begin{equation}\label{GoeritzForKnot52N1}
\dbwidetilde{G}{}^{5_2}_{II}  =
\begin{pmatrix}
 \begin{array}{ccccc}  1&2&0&-1&-1  \\  1&-2&1&0&0\\0&1&-2&1&0\\-1&0&1&1&-1\\
 -1&0&0&-1&2
 \end{array}
 \end{pmatrix} \ \ \mapsto \ \ G_4^{5_2}  =
\begin{pmatrix}
 \begin{array}{cccc}  1&1&0&-1  \\  1&-2&1&0\\0&1&-2&1\\-1&0&1&1
 \end{array}
 \end{pmatrix}, 
\end{equation}
Using these matrices, one can verify the invariance of the absolute value of the knot determinant for different realizations of the Goeritz matrix:  
${\rm Det}^{5_2} = |\det(G_1^{5_2})| = |\det(G_2^{5_2})| = |\det(G_3^{5_2})| = |\det(G_4^{5_2})| = 7$.  
Note that this invariant, in particular, distinguishes the knots $3_1$ and $5_2$.

\section{Method for obtaining Jones polynomials from the Goeritz matrix}\label{sec:JonesFromGoeritzMatrix}
\setcounter{equation}{0}

Let us choose the Goeritz matrix corresponding to a given knot diagram. To obtain the Jones polynomial, we seek a function $\mu[G]$ of the Goeritz matrix $G$. It turns out that this function corresponds to the Kauffman bracket, as was rigorously proven in~\cite{boninger2023jones}. Therefore, we can obtain the Jones polynomial simply by applying the proper normalization.\footnote{For the Kauffman bracket, the normalization factor depends only on the algebraic sum of crossings ${\rm W}^{\cal L}$. In the Goeritz matrix approach, however, the same normalization factor does not always yield a correct topological result, since the Goeritz matrices are insensitive to self-crossings of regions, while ${\rm W}^{\cal L}$ changes. To account for this in the normalization, one needs to retouch the coloring by counting all self-crossings of the regions (see Section~\ref{sec:RetMu}).} The algorithm is described in Sections~\ref{sec:GoeritzMatTransf}--\ref{sec:ExceptionalCases}, where we also demonstrate its relation to the Kauffman bracket. Before that in Section~\ref{sec:before-ex}, we examine a preliminary example to build intuition about the transformations of Goeritz matrices and about the coefficients appearing in the Boninger algorithmic step which should correspond to the coefficients in the Kauffman bracket.

\subsection{Preliminary example}
\label{sec:before-ex}

As a starting point let us illustrate the process using the example of the trefoil knot and its corresponding Goeritz matrices~\eqref{G-mat-3-1-2}, showing how the Kauffman bracket resolutions transform the diagram, the checkerboard coloring, and the Goeritz matrix (see Fig.\,\ref{fig:ExampelTrefoilKauffmanGoeritz}). At each step, we explicitly draw the diagrams corresponding to the different Goeritz matrices in order to visualize the respective matrix transformations. The coefficients appearing in the application of the Kauffman bracket are indicated above the arrows.

\begin{figure}[h!]
\centering
\begin{picture}(300,350)(50,-10)

\put(0,143){

\put(170,110){\vector(2,-1){60}}

\put(55,100){\mbox{\small $(-q)^{+1/2}$}}

\put(200,100){\mbox{\small $(-q)^{-1/2}=P_1(q)$}}

\put(115,110){\vector(-2,-1){60}}

\put(172,157){\mbox{\footnotesize $\widetilde{G} = \begin{pmatrix}
        -2 & \boxed{1} & 1 \\
        1 & -2 & 1 \\
        1 & 1 & -2
\end{pmatrix}$}}

\put(172,124){\mbox{\footnotesize $G=\begin{pmatrix}
    -2 & \boxed{1} \\
    1 & -2
\end{pmatrix}$}}

\put(57,65){\mbox{\footnotesize $\widetilde{G}=\begin{pmatrix}
    -1 & 0 & 1 \\
    0 & -1 & \boxed{1} \\
    1 & 1 & 2
\end{pmatrix}$}}

\put(57,30){\mbox{\footnotesize $G=\begin{pmatrix}
    -1 & 0 \\
    0 & \boxed{-1}
\end{pmatrix}$}}

\put(285,65){\mbox{\footnotesize $\widetilde{G}=\begin{pmatrix}
    -2 & \boxed{2} \\
    2 & -2 
\end{pmatrix}$}}

\put(285,30){\mbox{\footnotesize $G=-2$}}

\put(40,17){\vector(0,-1){30}}

\put(17,17){\vector(0,-1){30}}

\put(47,-2){\mbox{\small $(-q)^{-1/2}$}}

\put(-35,-2){\mbox{\small $(-q)^{+1/2}D_2 \ \ \ \ + $}}

\put(57,-35){\mbox{\footnotesize $\widetilde{G}=\begin{pmatrix}
    -1 & \boxed{1} \\
    1 & -1 
\end{pmatrix}$}}

\put(57,-65){\mbox{\footnotesize $G=-1$}}

\put(57,-130){\mbox{\footnotesize $\widetilde{G}=0$}}

\put(57,-145){\mbox{\footnotesize $G=\varnothing$}}

\put(225,105){

\put(57,-130){\mbox{\footnotesize $\widetilde{G}=0$}}

\put(57,-145){\mbox{\footnotesize $G=\varnothing$}}
}

\put(0,-107){

\put(20,17){\vector(0,-1){30}}

\put(35,17){\vector(0,-1){30}}

\put(47,-2){\mbox{\small $(-q)^{-1/2}= (-q)^{+1/2}D_2 + P_1(q)$}}

\put(-35,-2){\mbox{\small $(-q)^{+1/2}D_2 \ \ \ \ + $}}

}

\put(265,17){\vector(0,-1){30}}

\put(248,17){\vector(0,-1){30}}

\put(220,0){\mbox{\small $-q D_2 \ \ \; \, + $}}

\put(273,0){\mbox{\small $2 - q^{-1} D_2 = (-q)^{+1/2\cdot 2} D_2 + P_2(q)$}}

}

\begin{tikzpicture}[scale=0.4]

\fill[bright-gray] (-9.5,-26) circle (0.7);
\draw[thick] (-9.5,-26) circle (0.7);

\fill[bright-gray] (10.5,-16.5) circle (0.7);
\draw[thick] (10.5,-16.5) circle (0.7);

\node at (0.5, 1.8) {1};
\node at (0.5, -0.5) {2};
\node at (0.5, -2.5) {3};
\node at (0.5, -4.8) {1};

\fill[bright-gray] (0,0) -- (0,1) -- (0.5,0.5) -- cycle; 
\fill[bright-gray] (1,1) -- (1,0) -- (0.5,0.5) -- cycle; 

\fill[bright-gray] (0,-2) -- (0,-1) -- (0.5,-1.5) -- cycle; 
\fill[bright-gray] (1,-1) -- (1,-2) -- (0.5,-1.5) -- cycle; 

\fill[bright-gray] (0,-4) -- (0,-3) -- (0.5,-3.5) -- cycle; 
\fill[bright-gray] (1,-3) -- (1,-4) -- (0.5,-3.5) -- cycle; 

\draw[thick] (0,0) -- (1,1); 
\draw[thick] (1,0) -- (0.6,0.4);
\draw[thick] (0,1) -- (0.4,0.6);

\draw[thick] (0,-2) -- (1,-1);
\draw[thick] (1,-2) -- (0.6,-1.6);
\draw[thick] (0,-1) -- (0.4,-1.4);

\draw[thick] (0,-4) -- (1,-3);
\draw[thick] (1,-4) -- (0.6,-3.6);
\draw[thick] (0,-3) -- (0.4,-3.4);

\fill[bright-gray] (1,1) .. controls (3,3) and (3,-6) .. (1,-4);
\draw[thick] (1,1) .. controls (3,3) and (3,-6) .. (1,-4);
\fill[white] (1,0) .. controls (1.5,-0.5) .. (1,-1);
\draw[thick] (1,0) .. controls (1.5,-0.5) .. (1,-1);
\fill[white] (1,-2) .. controls (1.5,-2.5) .. (1,-3);
\draw[thick] (1,-2) .. controls (1.5,-2.5) .. (1,-3);

\fill[bright-gray] (0,1) .. controls (-2,3) and (-2,-6) .. (0,-4);
\draw[thick] (0,1) .. controls (-2,3) and (-2,-6) .. (0,-4);
\fill[white] (0,0) .. controls (-0.5,-0.5) .. (0,-1);
\draw[thick] (0,0) .. controls (-0.5,-0.5) .. (0,-1);
\fill[white] (0,-2) .. controls (-0.5,-2.5) .. (0,-3);
\draw[thick] (0,-2) .. controls (-0.5,-2.5) .. (0,-3);

\begin{scope}[shift={(-10,-17)}]
\node at (0.5, 1.8) {1};
\node at (0.5, -2.5) {23};
\node at (0.5, -4.8) {1};

\fill[bright-gray] (0,1) rectangle (1,-2);

\fill[white] (0,1) .. controls (0.2,0.7) and (0.8,0.7) .. (1,1);
\draw[thick] (0,1) .. controls (0.2,0.7) and (0.8,0.7) .. (1,1);

\fill[white] (0,-2) .. controls (0.2,-1.8) and (0.8,-1.8) .. (1,-2);
\draw[thick] (0,-2) .. controls (0.2,-1.8) and (0.8,-1.8) .. (1,-2);

\fill[bright-gray] (0,-4) -- (0,-3) -- (0.5,-3.5) -- cycle; 
\fill[bright-gray] (1,-3) -- (1,-4) -- (0.5,-3.5) -- cycle; 

\draw[thick] (0,-4) -- (1,-3);
\draw[thick] (1,-4) -- (0.6,-3.6);
\draw[thick] (0,-3) -- (0.4,-3.4);

\fill[bright-gray] (1,1) .. controls (3,3) and (3,-6) .. (1,-4);
\draw[thick] (1,1) .. controls (3,3) and (3,-6) .. (1,-4);
\fill[white] (1,-2) .. controls (1.5,-2.5) .. (1,-3);
\draw[thick] (1,-2) .. controls (1.5,-2.5) .. (1,-3);

\fill[bright-gray] (0,1) .. controls (-2,3) and (-2,-6) .. (0,-4);
\draw[thick] (0,1) .. controls (-2,3) and (-2,-6) .. (0,-4);
\fill[white] (0,-2) .. controls (-0.5,-2.5) .. (0,-3);
\draw[thick] (0,-2) .. controls (-0.5,-2.5) .. (0,-3);
\end{scope}

\begin{scope}[shift={(-10,-8)}]
\node at (0.5, 1.8) {1};
\node at (0.5, -0.5) {2};
\node at (0.5, -2.5) {3};
\node at (0.5, -4.8) {1};

\fill[bright-gray] (0,0) rectangle (1,1); 

\fill[bright-gray] (0,-2) -- (0,-1) -- (0.5,-1.5) -- cycle; 
\fill[bright-gray] (1,-1) -- (1,-2) -- (0.5,-1.5) -- cycle; 

\fill[bright-gray] (0,-4) -- (0,-3) -- (0.5,-3.5) -- cycle; 
\fill[bright-gray] (1,-3) -- (1,-4) -- (0.5,-3.5) -- cycle; 

\fill[white] (0,1) .. controls (0.2,0.7) and (0.8,0.7) .. (1,1);
\draw[thick] (0,1) .. controls (0.2,0.7) and (0.8,0.7) .. (1,1);

\fill[white] (0,0) .. controls (0.2,0.2) and (0.8,0.2) .. (1,0);
\draw[thick] (0,0) .. controls (0.2,0.2) and (0.8,0.2) .. (1,0);

\draw[thick] (0,-2) -- (1,-1);
\draw[thick] (1,-2) -- (0.6,-1.6);
\draw[thick] (0,-1) -- (0.4,-1.4);

\draw[thick] (0,-4) -- (1,-3);
\draw[thick] (1,-4) -- (0.6,-3.6);
\draw[thick] (0,-3) -- (0.4,-3.4);

\fill[bright-gray] (1,1) .. controls (3,3) and (3,-6) .. (1,-4);
\draw[thick] (1,1) .. controls (3,3) and (3,-6) .. (1,-4);
\fill[white] (1,0) .. controls (1.5,-0.5) .. (1,-1);
\draw[thick] (1,0) .. controls (1.5,-0.5) .. (1,-1);
\fill[white] (1,-2) .. controls (1.5,-2.5) .. (1,-3);
\draw[thick] (1,-2) .. controls (1.5,-2.5) .. (1,-3);

\fill[bright-gray] (0,1) .. controls (-2,3) and (-2,-6) .. (0,-4);
\draw[thick] (0,1) .. controls (-2,3) and (-2,-6) .. (0,-4);
\fill[white] (0,0) .. controls (-0.5,-0.5) .. (0,-1);
\draw[thick] (0,0) .. controls (-0.5,-0.5) .. (0,-1);
\fill[white] (0,-2) .. controls (-0.5,-2.5) .. (0,-3);
\draw[thick] (0,-2) .. controls (-0.5,-2.5) .. (0,-3);
\end{scope}

\begin{scope}[shift={(10,-8)}]

\node at (0.5, -0.5) {12};
\node at (0.5, -2.5) {3};
\node at (0.5, -4.8) {12};

\fill[bright-gray] (0,-2) -- (0,-1) -- (0.5,-1.5) -- cycle; 
\fill[bright-gray] (1,-1) -- (1,-2) -- (0.5,-1.5) -- cycle; 

\fill[bright-gray] (0,-4) -- (0,-3) -- (0.5,-3.5) -- cycle; 
\fill[bright-gray] (1,-3) -- (1,-4) -- (0.5,-3.5) -- cycle;

\fill[bright-gray] (0,1) .. controls (0.2,0.8) and (0.2,0.2) .. (0,0);
\draw[thick] (0,1) .. controls (0.2,0.8) and (0.2,0.2) .. (0,0);

\fill[bright-gray] (1,1) .. controls (0.8,0.8) and (0.8,0.2) .. (1,0);
\draw[thick] (1,1) .. controls (0.8,0.8) and (0.8,0.2) .. (1,0);

\draw[thick] (0,-2) -- (1,-1);
\draw[thick] (1,-2) -- (0.6,-1.6);
\draw[thick] (0,-1) -- (0.4,-1.4);

\draw[thick] (0,-4) -- (1,-3);
\draw[thick] (1,-4) -- (0.6,-3.6);
\draw[thick] (0,-3) -- (0.4,-3.4);

\fill[bright-gray] (1,1) .. controls (3,3) and (3,-6) .. (1,-4);
\draw[thick] (1,1) .. controls (3,3) and (3,-6) .. (1,-4);
\fill[white] (1,0) .. controls (1.5,-0.5) .. (1,-1);
\draw[thick] (1,0) .. controls (1.5,-0.5) .. (1,-1);
\fill[white] (1,-2) .. controls (1.5,-2.5) .. (1,-3);
\draw[thick] (1,-2) .. controls (1.5,-2.5) .. (1,-3);

\fill[bright-gray] (0,1) .. controls (-2,3) and (-2,-6) .. (0,-4);
\draw[thick] (0,1) .. controls (-2,3) and (-2,-6) .. (0,-4);
\fill[white] (0,0) .. controls (-0.5,-0.5) .. (0,-1);
\draw[thick] (0,0) .. controls (-0.5,-0.5) .. (0,-1);
\fill[white] (0,-2) .. controls (-0.5,-2.5) .. (0,-3);
\draw[thick] (0,-2) .. controls (-0.5,-2.5) .. (0,-3);
  \end{scope}

\end{tikzpicture}
\end{picture}
\caption{\footnotesize Resolution of the trefoil knot using the Kauffman bracket and the corresponding reduced and unreduced Goeritz matrices. Above the arrows (according to the Kauffman bracket, see Fig.\,\ref{fig:Kauff}), there are the coefficients with which the trefoil resolutions must be weighted in order to obtain the Jones polynomial. At the second and third steps, the horizontal and vertical Kauffman bracket resolutions yield identical diagrams up to a cycle, which we factor out as the multiplier $D_2$. Therefore, in the figure instead of showing two identical diagrams we display only one. The coefficient for this diagram is simply the sum of the coefficients corresponding to the horizontal and vertical resolutions. The functions $P_i(q)$ are obtained in the same way as in formula~\eqref{Step3}. In the Goeritz matrices, we highlight the crossing numbers that are to be resolved at the corresponding next step.}
\label{fig:ExampelTrefoilKauffmanGoeritz}
\end{figure}

Let us list the observations made from the example of the two-strand trefoil knot (see Fig.\,\ref{fig:ExampelTrefoilKauffmanGoeritz}), which should be reflected in the general algorithm for computing the Jones polynomial in Sections~\ref{sec:GoeritzMatTransf} and~\ref{sec:muInv}.

\paragraph{Step 1.} At the first step, we resolve a single crossing between regions 1 and 2. Note that the horizontal resolution does not change the number of regions and therefore does not affect the size of the Goeritz matrix. In this case, the number of crossings between regions 1 and 2 becomes zero, which corresponds to the matrix transformation (I) described in Section~\ref{sec:GoeritzMatTransf}. The vertical resolution, on the other hand, merges regions 1 and 2 into a single region 12. Hence, the size of the Goeritz matrix decreases by one, and the number of crossings between regions 12 and 3 becomes equal to the sum of the crossings between regions 1 and 3 and between 2 and 3, which fully corresponds to transformation (II) for the Goeritz matrices.

\paragraph{Steps 2 and 3.} At these steps, each region crosses only one other region\footnote{For example, in the right diagram after the first step, region 12 crosses only with region 3.}, and the Goeritz matrix becomes diagonal. In this case, both the horizontal and vertical Kauffman resolutions produce the same knot (if the isolated cycle is replaced by the factor $D_2$). The corresponding transformation on the matrices acts on the diagonal elements and removes the respective row and column after resolving the crossings, just as in transformation (III) described in Section~\ref{sec:GoeritzMatTransf}.

\medskip
 
\noindent Thus, we see that the natural transformation of the unreduced Goeritz matrices acts at each step on a single off-diagonal element of the matrix, resolving the corresponding crossings in the knot. In the final steps when each region crosses only one other region, we can act on diagonal elements (instead of off-diagonal ones), which allows us to work with the reduced Goeritz matrix instead of the unreduced one. The operations on the Goeritz matrices are in full agreement with the general method described in Section~\ref{sec:GoeritzMatTransf}, while the coefficients from the Kauffman bracket correspond to formulas~\eqref{Step3}--\eqref{Step0}.

\subsection{Matrix transformations}
\label{sec:GoeritzMatTransf}

As we have already illustrated in the example from the previous section, at each step of the algorithm, all crossings between adjacent regions in a knot diagram are resolved. Therefore, each step of the algorithm should produce not only diagrams of new links but also transformations of the Goeritz matrices with respect to a matrix of an initial diagram. All such transformations can be classified into three types.

Once again, we derive these transformations starting from the Kauffman bracket. As we have already noted, elements of the Goeritz matrix are numbers of crossings between adjacent regions and at each step of the algorithm we resolve such crossings between two regions. Consider the simplest case (see Fig.\,\ref{fig:Kauff-braid}), where the regions $i$ and $j$ cross along a single braid. By applying the Kauffman bracket sequentially to each crossing in the braid, we ultimately obtain its vertical and horizontal resolutions. From this, the first two transformations of the Goeritz matrices naturally follow.

\bigskip

\noindent {\bf (I)} The first transformation $G \;\mapsto\; G^{(I)}$ corresponds to the vertical resolution of a braid. It is clear that under vertical resolution the number of crossings $g_{ij}$ becomes zero, which is also reflected correspondingly in the two diagonal elements $g_{ii}$ and $g_{jj}$:

\begin{equation}\label{GoeritzMovie1}
G=\begin{pmatrix}
\ddots & \\
& g_{ii} & \dots & \boxed{g_{ij}} \\
& \vdots & \ddots & \vdots \\
& g_{ji} & \dots & g_{jj} \\
& & & & \ddots \\
\end{pmatrix}\quad \mapsto \quad
G^{(I)}=\begin{pmatrix}
\ddots & \\
& g_{ii} + g_{ij} & \dots & 0 \\
& \vdots & \ddots & \vdots \\
& 0 & \dots & g_{jj} + g_{ij} \\
& & & & \ddots \\
\end{pmatrix}.   
\end{equation}
The explicit transformation of the matrix elements is given by:
\begin{equation}
\begin{aligned}
g_{i i} \quad &\mapsto \quad g_{i i}+g_{i j}\,, \\
g_{j j} \quad &\mapsto \quad g_{j j}+g_{i j}\,, \\
g_{i j},\, g_{j i} \quad &\mapsto \quad 0\,.
\end{aligned}
\end{equation}

\noindent {\bf (II)} The second transformation $G \;\mapsto\; G^{(II)}$ corresponds to the horizontal resolution of a braid and reduces the size of the matrices, since the regions $i$ and $j$ are merged into a single region $ij$. It is clear that if the initial regions had crossings with other regions, then the crossing numbers of the unified region $ij$ with those regions are simply added: $g_{ij,k} = g_{ik} + g_{jk}$. Accordingly, the diagonal elements also change, as they must be equal to the sum of all crossings with a fixed region taken with the opposite sign:
\begin{equation}\label{GoeritzMovie2}
\begin{picture}(300,100)(90,-45)
    \put(0,0){\mbox{$G=\begin{pmatrix}
\ddots & & & & \vdots \\
& g_{ii} & g_{ik} & \dots & \boxed{g_{ij}} \\
& g_{ki} & \ddots & & g_{kj} \\
& \vdots & & \ddots & \vdots \\
\dots & g_{ji} & g_{jk} & \dots & g_{jj} & \dots \\
& & & & \vdots &  \ddots \\
\end{pmatrix} \quad \mapsto \quad G^{(II)} =
\begin{pmatrix}
\ddots & & & & {\vdots} \\
& g_{ii} + g_{jj} + 2g_{ij} & g_{ik} + g_{jk} & \dots & g_{ij} \\
& g_{ki} + g_{kj} & \ddots & & g_{kj} \\
& \vdots & & \ddots & \vdots \\
{\dots} & g_{ji} & g_{jk} & \dots & g_{jj} & {\dots} \\
& & & & {\vdots} &  \ddots \\
\end{pmatrix}. $}}

\thicklines

{\color{red}

\put(420,-53){\line(0,1){110}}
\put(237,-25){\line(1,0){225}}

}

\end{picture}
\end{equation}

The explicit transformation of the matrix elements:
\begin{equation}
\begin{aligned}
& g_{i i} \quad \mapsto \quad g_{i i}+g_{j j}+2 g_{i j}\,, \\
& g_{i k} \quad \mapsto \quad g_{i k}+g_{j k}, \; \text { for all } \; k \neq i\,, \\
& g_{k i} \quad \mapsto \quad g_{k i}+g_{k j}, \; \text { for all } \; k \neq i\,.
\end{aligned}
\end{equation}

\noindent {\bf (III)} According to our trefoil example from the previous section, there must exist a third type of transformation that acts on diagonal Goeritz matrices. Such matrices correspond to the case in which all regions, except one, cross only with a single other region, whose corresponding row and column were removed during the reduction of the Goeritz matrix. An example of this situation is shown in Fig.\,\ref{fig:diag-Goeritz-ex}. It can be seen that applying the Kauffman bracket to any braid yields the same knot diagram up to a cycle that can be replaced by the factor $D_2$. Thus, the third transformation $G \;\mapsto\; G^{(III)}$ acts on a diagonal matrix by removing the $i$-th row and column since the crossings with a single region are eliminated together with the corresponding region itself:

\begin{equation}
    G = {\rm diag}(\dots,\, g_{i-1},\, \boxed{g_i}\,,\, g_{i+1},\, \dots) \quad \mapsto \quad G^{(III)} = {\rm diag}(\dots ,\, g_{i-1},\, \cancel{g_i},\, g_{i+1},\, \dots)\,. 
\end{equation}

\subsection{The function $\mu$ on Goeritz matrices and Jones polynomial}
\label{sec:muInv}

\begin{figure}[h!]
\centering
\begin{picture}(300,150)(-30,5)

\put(95,66){\vector(2,-1){50}}
\put(80,66){\vector(-2,-1){50}}

\put(20,55){\mbox{$(-q)^{g_{ij}/2}$}}

\put(125,55){\mbox{$P_{g_{ij}}(q)$}}


\begin{tikzpicture}[scale=0.4]

\fill[bright-gray] (10,0) rectangle (11,-1);

\fill[bright-gray] (10,1) -- (11,1) -- (10.5,0.5) -- cycle; 

\fill[bright-gray] (10,0) -- (11,0) -- (10.5,0.5) -- cycle; 

\fill[bright-gray] (10,-1) -- (11,-1) -- (10.5,-1.5) -- cycle; 

\fill[bright-gray] (10,-2) -- (11,-2) -- (10.5,-1.5) -- cycle;

\draw[thick] (10,0) -- (11,1); 
\draw[thick] (11,0) -- (10.65,0.35);
\draw[thick] (10,1) -- (10.35,0.65);

\draw[thick] (10,-2) -- (11,-1);
\draw[thick] (11,-2) -- (10.65,-1.65);
\draw[thick] (10,-1) -- (10.35,-1.35);

\fill[bright-gray] (11,0) .. controls (11.5,-0.5) .. (11,-1);
\draw[thick] (11,0) .. controls (11.5,-0.5) .. (11,-1);

\fill[bright-gray] (10,0) .. controls (9.5,-0.5) .. (10,-1);
\draw[thick] (10,0) .. controls (9.5,-0.5) .. (10,-1);

\begin{scope}[shift={(0,-2.5)}]

\fill[bright-gray] (10,-2) rectangle (11,-3);

\fill[bright-gray] (10,-3) -- (11,-3) -- (10.5,-3.5) -- cycle; 

\fill[bright-gray] (10,-4) -- (11,-4) -- (10.5,-3.5) -- cycle; 

\fill[bright-gray] (10,-1) -- (11,-1) -- (10.5,-1.5) -- cycle; 

\fill[bright-gray] (10,-2) -- (11,-2) -- (10.5,-1.5) -- cycle;

\draw[thick] (10,-2) -- (11,-1);
\draw[thick] (11,-2) -- (10.65,-1.65);
\draw[thick] (10,-1) -- (10.35,-1.35);

\draw[thick] (10,-4) -- (11,-3);
\draw[thick] (11,-4) -- (10.65,-3.65);
\draw[thick] (10,-3) -- (10.35,-3.35);

\fill[bright-gray] (11,-2) .. controls (11.5,-2.5) .. (11,-3);
\draw[thick] (11,-2) .. controls (11.5,-2.5) .. (11,-3);

\fill[bright-gray] (10,-2) .. controls (9.5,-2.5) .. (10,-3);
\draw[thick] (10,-2) .. controls (9.5,-2.5) .. (10,-3);
\end{scope}

\fill[bright-gray] (17.5,1) rectangle (18.5,-1);

\fill[bright-gray] (17.5,-6.5) rectangle (18.5,-4.5);

\draw[thick] (17,-4.5) rectangle (19,-1);

\draw[thick] (17.5,1) -- (17.5,-1);
\draw[thick] (18.5,1) -- (18.5,-1);

\draw[thick] (17.5,-6.5) -- (17.5,-4.5);
\draw[thick] (18.5,-6.5) -- (18.5,-4.5);

\node at (15.8, -2.5) {$i$};
\node at (20, -2.5) {$j$};

\node at (9, -2.5) {$i$};
\node at (12, -2.5) {$j$};
\node at (10.5, -2.5) {$\vdots$};
\node at (14, -2.5) {$\equiv$};
\node at (18, -2.5) {$g_{ij}$};

\fill[bright-gray] (8,-10) rectangle (9,-13);

\draw[thick] (8,-10) -- (8,-13);
\draw[thick] (9,-10) -- (9,-13);

\node at (6.9, -11.5) {$i$};
\node at (10.1, -11.5) {$j$};

\fill[bright-gray] (18,-10) .. controls (19,-11) and (20,-11) .. (21,-10);

\draw[thick] (18,-10) .. controls (19,-11) and (20,-11) .. (21,-10);

\fill[bright-gray] (18,-13) .. controls (19,-12) and (20,-12) .. (21,-13);

\draw[thick] (18,-13) .. controls (19,-12) and (20,-12) .. (21,-13);

\node at (19.5,-11.5) {$ij$};

\end{tikzpicture}
\end{picture}
\caption{\footnotesize Resolution of a two-strand braid using the Kauffman bracket. The coefficients associated with each resolution are shown above the arrows.} 
\label{fig:Kauff-braid}
\end{figure}

Our goal now is to introduce a function $\mu[g_{ij}]$ on Goeritz matrices that reproduces the Kauffman bracket. We will be guided by our preliminary example from Section~\ref{sec:before-ex} and by the result of applying the Kauffman bracket. Initially, we deal with a non-diagonal Goeritz matrix. To determine how the function $\mu$ acts on the non-diagonal elements, we compute the coefficients of the two-strand braid resolution using the Kauffman bracket (see Fig.\,\ref{fig:Kauff-braid}). For the vertical resolution of a braid with $g_{ij}$ crossings, corresponding to the matrix $G^{(I)}$ (as explained in Section~\ref{sec:GoeritzMatTransf}), the coefficient is simply the factor $(-q)^{g_{ij}/2}$, while for the horizontal resolution (corresponding to the matrix $G^{(II)}$), it is given by the combinatorial factor $P_{g_{ij}}(q)$ (see Fig.\,\ref{fig:Kauff-braid}):
{\small \begin{equation}\label{Step3} 
P_{g_{ij}}(q) = \sum^{|g_{ij}|}_{k=1} \begin{pmatrix}
|g_{ij}| \\
k 
\end{pmatrix} D_2^{k-1} (-q)^{{\rm sgn}(g_{ij})(|g_{ij}|-2k)/2} = (-q)^{g_{ij}/2}\;\frac{\left(1+(-q)^{-{\rm sgn}(g_{ij})}D_2 \right)^{|g_{ij}|}-1}{D_2} = (-q)^{g_{ij}/2}\;\frac{(-1)^{g_{ij}} q^{-2g_{ij}}-1}{D_2}\,.
\end{equation}}

\noindent This factor can be easily computed inductively by resolving the crossings in a two-strand braid. Thus, the action of the Kauffman bracket on non-diagonal elements $g_{ij}$, $(i \neq j)$ of the Goeritz matrix is given by:
\begin{equation}\label{Step1}
\boxed{\mu [g_{ij}] = (-q)^{g_{ij}/2} \mu [g_{ij}^{(I)}] + P_{g_{ij}}(q)\mu [g_{ij}^{(II)}]\,.}
\end{equation}

\begin{figure}[h!]
    \centering
\begin{picture}(300,105)(0,5)
\begin{tikzpicture}[scale=0.4]

\fill[bright-gray] (9.5,-1) ellipse [x radius=12, y radius=5];

\fill[white] (0,0) rectangle (4,-2);
\fill[white] (0,-0.5) .. controls (-3,2.5) and (7,2.5) .. (4,-0.5);
\fill[white] (0,-1.5) .. controls (-4,-5) and (23,-5) .. (19,-1.5);

\draw[thick] (0,0) rectangle (4,-2);
\draw[thick] (0,-0.5) .. controls (-3,2.5) and (7,2.5) .. (4,-0.5);

\draw[thick] (4,-1.5) -- (6.5,-1.5);

\draw[thick] (10.5,-1.5) -- (11.5,-1.5);

\draw[thick] (13.5,-1.5) -- (15,-1.5);

\draw[thick] (0,-1.5) .. controls (-4,-5) and (23,-5) .. (19,-1.5);

\node at (12.5,-1.5) {$\ldots$}; 

\begin{scope}[shift={(6.5,0)}]

\fill[white] (0,0) rectangle (4,-2);
\fill[white] (0,-0.5) .. controls (-3,2.5) and (7,2.5) .. (4,-0.5);

\draw[thick] (0,0) rectangle (4,-2);
\draw[thick] (0,-0.5) .. controls (-3,2.5) and (7,2.5) .. (4,-0.5);
\end{scope}

\begin{scope}[shift={(15,0)}]

\fill[white] (0,0) rectangle (4,-2);
\fill[white] (0,-0.5) .. controls (-3,2.5) and (7,2.5) .. (4,-0.5);

\draw[thick] (0,0) rectangle (4,-2);
\draw[thick] (0,-0.5) .. controls (-3,2.5) and (7,2.5) .. (4,-0.5);
\end{scope}

\node at (2,-1) {$g_{1n}$};
\node at (8.5,-1) {$g_{2n}$};
\node at (17,-1) {$g_{n-1,n}$};
\end{tikzpicture}
\end{picture}
    \caption{\footnotesize An example of a link that arises in the final steps of the Goeritz method for computing the Jones polynomial, when the matrix becomes diagonal. It can be seen that both the vertical and horizontal resolutions of any braid within rectangle yield the same diagram, up to a single cycle.}
    \label{fig:diag-Goeritz-ex}
\end{figure}

By successively applying this procedure to the Goeritz matrices we arrive at diagonal matrices. For corresponding links the Kauffman bracket acts more simply, since both resolutions of the braid produce the same knot diagram (up to a cycle) but with different coefficients. From Fig.\,\ref{fig:diag-Goeritz-ex} it is easy to see that the function $\mu$ for a diagonal element takes the following form:
\begin{equation}\label{Step2}
\boxed{\mu [g_{ii}] = \left((-q)^{-g_{ii}/2}D_2 + P_{-g_{ii}}(q)\right) \mu [g_{i}^{(III)}]\,.}
\end{equation}
At the final step we obtain a cycle, which contributes the factor $D_2 = q + q^{-1}$ to the Jones polynomial. The corresponding Goeritz matrix is $G = \varnothing$, and therefore we set
\begin{equation}\label{Step0}
\mu[\varnothing]=D_2=q+q^{-1}\,.
\end{equation}
We now explicitly state the following proposition that we have demonstrated. A rigorous mathematical proof that the Goeritz matrix invariant $\mu[g_{ij}]$ corresponds to the Kauffman bracket can be found in~\cite{boninger2023jones}.

\paragraph{Proposition.} The polynomial $\mu [G(\mathcal{D})]$ constructed from the Goeritz matrix of a link diagram $\mathcal{D}$ corresponding to a link $\mathcal{L}$ represents the result of resolving the diagram $\mathcal{D}$ using the Kauffman bracket. Thus, the Jones polynomial can be obtained from it as follows:
\begin{equation}\label{Jones} J^{\cal  {L}}(q) = (-(-q)^{3/2})^{-\rm W({\cal D})+\rm W_r({\cal D})}\mu [G]\,.
\end{equation}
Recall that $\rm W$ is the algebraic sum of the crossings in the diagram. The quantity $\rm W_r(\mathcal{D})$ counts the crossings corresponding to the first Reidemeister move which are not distinguished by the Goeritz matrix. In other words, $\rm W_r(\mathcal{D})$ represents the number of self-crossings of regions. The details of this normalization are discussed further in Section~\ref{sec:RetMu}.

\paragraph{Inaccuracy in Boninger’s normalization.}  
In~\cite{boninger2023jones}, the following normalization of the Jones polynomial is presented\footnote{To avoid unnecessary confusion, we reproduce Boninger normalization here using our notation.}:

\begin{equation}\label{normJones} J^{\cal  {L}}(q) = (-(-q)^{1/2})^{3\left(\sum_{i\le j}g_{ij}-e(S,\cal{L})\right)}\mu [G]\,,
\end{equation}
where $e(S,\mathcal{L})$ is a quantity that depends on the checkerboard coloring $S$ and the link $\mathcal{L}$; for knots --- $e(S,\mathcal{K}) = 0$.

This normalization is incorrect, because the sum of the upper-triangular elements of the Goeritz matrix does not equal the algebraic sum of crossings for the diagram. We illustrate this on the example of the knot $5_2$ (see Fig.\,\ref{pic:Knot52}). It is easy to observe that the sum of the upper-triangular elements of any Goeritz matrix \eqref{GoeritzForKnot52N1} for the knot $5_2$ vanishes $\sum_{i\le j} g_{ij}=0$, whereas ${\rm W}=4$. The correct normalization in this case is given by the algebraic sum of crossings, as can be seen explicitly:

\begin{equation}\label{ExampleTorKnot2}
    J^{5_2} =(-(-q)^{3/2})^{-4}\left( (-q)^{-1/2}\mu\begin{bmatrix}\begin{array}{cc}  -3 & 0\\0 & 2     
 \end{array}
     
 \end{bmatrix} + P_{-1}\mu[-1]\right) = \frac{-1 + q^2 - q^4 + 2 q^6 - q^8 + q^{10}}{q^{12}}\,.
\end{equation}

\subsection{Exceptional cases}\label{sec:ExceptionalCases}

In this section, we describe cases where two regions in the Goeritz matrix are separated by more than a single braid and we also consider how the first Reidemeister move is reflected in the Goeritz algorithm. In the context of the first Reidemeister move, we will highlight some interesting features of the Goeritz matrix and the corresponding approach to computing Jones polynomials.

\subsubsection{Retouching a polynomial $\mu[G]$}\label{sec:RetMu}

\begin{figure}[h!]
		\centering	
		\includegraphics[width =0.5\linewidth]{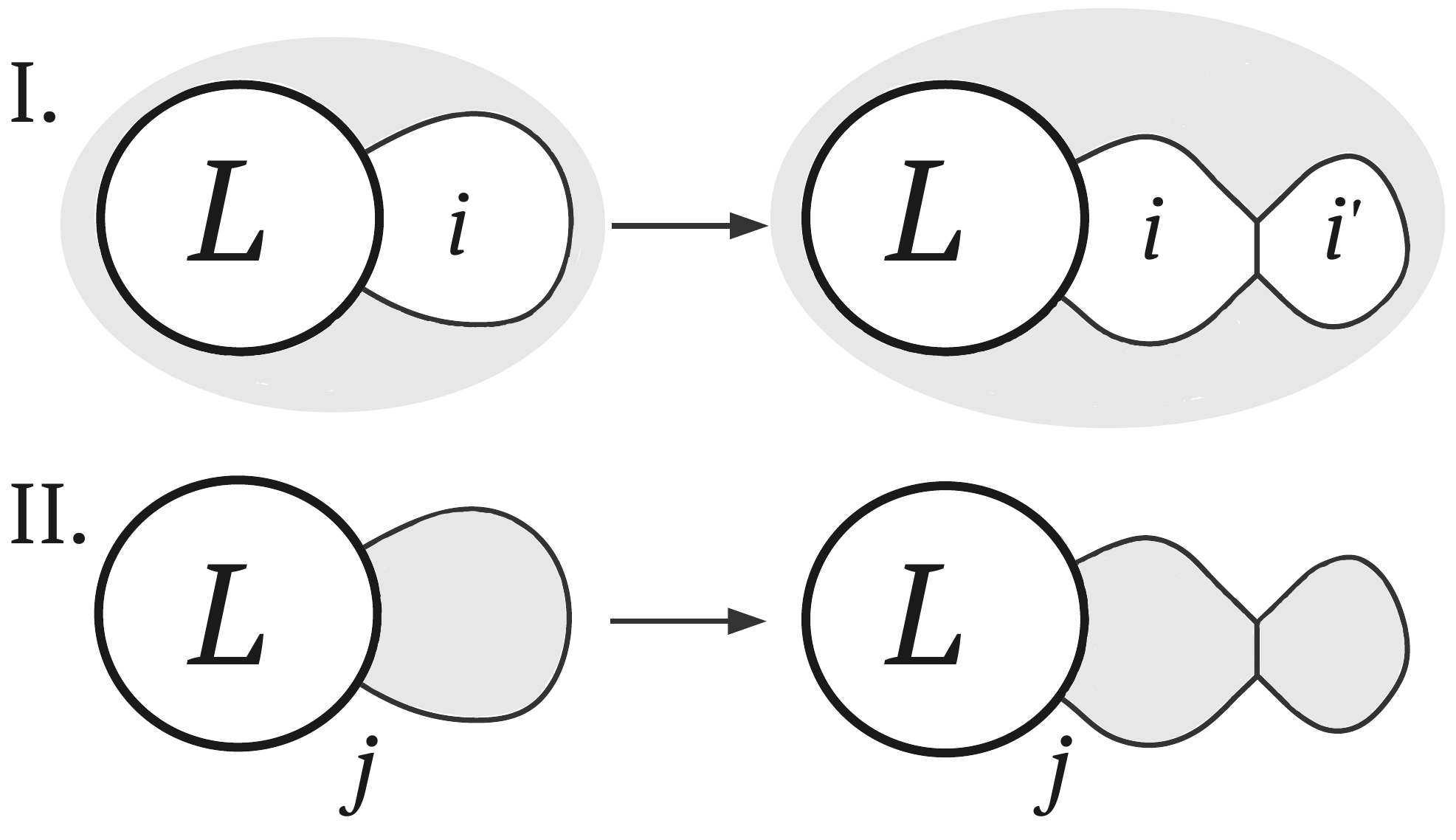}
        \caption{\footnotesize The first Reidemeister move on the diagram and the corresponding checkerboard coloring.}
        \label{fig:Retoucher1}
	\end{figure}

Let us now consider how the checkerboard coloring of a diagram changes under the first Reidemeister move (I.R) (see Fig.\,\ref{fig:Retoucher1}). In the first coloring option, the the first Reidemeister move adds a new white region $i'$, which intersects only with the region $i$. The Goeritz matrix therefore changes (its dimension increases by one), and consequently, the expression for $\mu[G]$ acquires an additional factor $(-(-q)^{\pm3/2})$, when computing the Jones polynomial, this factor is compensated by the change in the algebraic crossing number $\rm W' = \rm W \pm 1$.
	
The second coloring is of particular interest. As illustrated in Fig.\,\ref{fig:Retoucher1}, the number of white regions in new checkerboard coloring remains unchanged. Moreover, an additional crossing does not affect the Goeritz matrix, since by construction only crossings between \emph{distinct} regions are taken into account, while self-crossings of a single region are ignored. Because the new crossing leaves the matrix $G$ unchanged, the polynomial $\mu[G]$ likewise remains invariant under this coloring\footnote{In many research contexts, such a subtlety does not even arise, as knots are typically constructed without such self-crossings in practical applications.}. 
Therefore, in order to obtain the correct Jones polynomial via the Goeritz approach with normalization $(-(-q)^{3/2})^{-{\rm W}({\cal D})}$, it is necessary to remove the extra factors that corresponds to invisible for the Goeritz matrix the first Reidemeister moves and to change of ${\rm W}$. In other words, one must count the number of self-crossings of the regions.

The first approach to adjust the normalization is through retouching the polynomial $\mu[G]$. That is, for the original diagram $\cal D$, one counts the number of loops $\rm W^+_r(\cal D)$ that are removed by a the first Reidemeister move (see Fig.\,\ref{fig:Retoucher1}) and correspond to positive crossings. Likewise, one defines the number of negative loops ${\rm W}^-_r(\cal D)$\footnote{It is important to note that there may exist loops that are intercrossed multiple times but are still removed by the first Reidemeister move. In such cases the contribution of these loops to $\rm W_r(\cal D)$ equals the number of crossings counted with the appropriate sign.}. 
Then, the quantity $\rm W_r(\cal D) = \rm W^+_r(\cal{D}) - \rm W^-_r(\cal D)$ will be a counter-normalization factor \eqref{Jones}. A simpler method is to compute the algebraic sum of self-crossings over all regions.

In the second approach one may use the factor $(-(-q)^{3/2})^{-{\rm W}({\cal D})}$, but the knot diagram must be redrawn so that no regions contain self-crossings which is always possible.

\subsubsection{Intersection of regions by several separated braids}\label{sec:SeveralAreas}

Suppose that in the diagram $\cal D$ under the checkerboard coloring $S$ there are $k$ disjoint two-strand tangles between regions $i$ and $j$, as shown in Fig.\,\ref{fig:sum-par}.  
Then the element of the Goeritz matrix is given by $g_{ij} = n_1 + n_2 + \dots + n_k$ and the function $\mu$ applied to this element resolves all these two-strand crossings simultaneously, so that the Jones polynomial depends only on the sum of the crossings in these tangles. In this section, we will show that this does not contradict the Kauffman bracket.

\begin{figure}[h!]
    \centering
\begin{picture}(300,65)(0,0)
\begin{tikzpicture}[scale=0.4]

\fill[bright-gray] (-2,1.5) rectangle (0,0.5);
\draw[thick] (-2,1.5) -- (0,1.5);
\draw[thick] (-2,0.5) -- (0,0.5);

\fill[bright-gray] (3,1.5) rectangle (5.05,0.5);
\draw[thick] (3,1.5) -- (5.1,1.5);
\draw[thick] (3,0.5) -- (5.1,0.5);

\draw[thick] (0,0) rectangle (3,2);
\node at (1.5,1) {$n_1$}; 

\fill[bright-gray] (6.95,1.5) rectangle (9,0.5);
\draw[thick] (6.9,1.5) -- (9,1.5);
\draw[thick] (6.9,0.5) -- (9,0.5);

\draw[thick] (6,1) circle (1); 
\node at (6.1,1) {${\cal L}_1$};

\fill[bright-gray] (12,1.5) rectangle (14.05,0.5);
\draw[thick] (12,1.5) -- (14.1,1.5);
\draw[thick] (12,0.5) -- (14.1,0.5);

\draw[thick] (9,0) rectangle (12,2);
\node at (10.5,1) {$n_2$}; 

\fill[bright-gray] (15.95,1.5) rectangle (17.5,0.5);
\draw[thick] (15.9,1.5) -- (17.5,1.5);
\draw[thick] (15.9,0.5) -- (17.5,0.5);

\draw[thick] (15,1) circle (1); 
\node at (15,1) {${\cal L}_2$};

\node at (18.5,1) {$\ldots$};

\fill[bright-gray] (19.6,1.5) rectangle (21.5,0.5);
\draw[thick] (19.6,1.5) -- (21.5,1.5);
\draw[thick] (19.6,0.5) -- (21.5,0.5);

\fill[bright-gray] (24.5,1.5) rectangle (26.5,0.5);
\draw[thick] (24.5,1.5) -- (26.5,1.5);
\draw[thick] (24.5,0.5) -- (26.5,0.5);

\draw[thick] (21.5,0) rectangle (24.5,2);
\node at (23,1) {$n_k$};

\node at (12,-1) {$i$};
\node at (12,3) {$j$};

\end{tikzpicture}    
\end{picture}
    \caption{\footnotesize Two regions intersect along disjoint by ${\cal L}_m$ two-strand braids. Both the Kauffman bracket and the Goeritz method imply that the Jones polynomial for links containing such tangles depends only on the sum of the parameters in the braids.}
    \label{fig:sum-par}
\end{figure}

\begin{figure}[h!]
    \centering
\begin{picture}(300,180)(100,0)
\begin{tikzpicture}[scale=0.4] 

\node at (1,5.5) {${\cal L}_1$};

\node at (-1.4,5.5) {$i$};
\node at (3.4,5.5) {$j$};

\fill[bright-gray] (0.5,12.5) rectangle (1.5,11);
\draw[thick] (0.5,12.5) -- (0.5,11);
\draw[thick] (1.5,12.5) -- (1.5,11);

\fill[bright-gray] (0.5,8) rectangle (1.5,6.4);
\draw[thick] (0.5,8) -- (0.5,6.35);
\draw[thick] (1.5,8) -- (1.5,6.35);

\node at (1,9.5) {$n$};
\draw[thick] (0,8) rectangle (2,11);

\fill[bright-gray] (0.5,4.6) rectangle (1.5,3);
\draw[thick] (0.5,4.65) -- (0.5,3);
\draw[thick] (1.5,4.65) -- (1.5,3);

\draw[thick] (1,5.5) circle (1);

\fill[bright-gray] (0.5,0) rectangle (1.5,-1.5);
\draw[thick] (0.5,0) -- (0.5,-1.5);
\draw[thick] (1.5,0) -- (1.5,-1.5);

\node at (1,1.5) {$m$};
\draw[thick] (0,0) rectangle (2,3);

\node at (5,5.5) {$=$};

\node at (1,-3) {$J^{{\cal L}_{n,m}}$};

\begin{scope}[shift={(11.5,0)}]

\node at (-3,5.5) {\footnotesize $(-q)^{(n+m)/2}$};


\fill[bright-gray] (0.5,10) rectangle (1.5,6.4);
\draw[thick] (0.5,10) -- (0.5,6.35);
\draw[thick] (1.5,10) -- (1.5,6.35);

\fill[bright-gray] (0.5,4.6) rectangle (1.5,1);
\draw[thick] (0.5,4.65) -- (0.5,1);
\draw[thick] (1.5,4.65) -- (1.5,1);

\node at (1,5.5) {${\cal L}_1$};
\draw[thick] (1,5.5) circle (1);

\node at (3.5,5.5) {$+$};

\node at (1,-0.5) {$J^{{\cal L}_0}$};

\end{scope}

\begin{scope}[shift={(21.5,0)}]

\node at (-3,5.5) {\small $(-q)^{m/2} P_n$};

\fill[bright-gray] (0.5,4.6) rectangle (1.5,1);
\draw[thick] (0.5,4.65) -- (0.5,1);
\draw[thick] (1.5,4.65) -- (1.5,1);

\fill[bright-gray] (0.5,6.4) .. controls (0.7,7.5) and (1.3,7.5) .. (1.5,6.4);
\draw[thick] (0.5,6.4) .. controls (0.7,7.5) and (1.3,7.5) .. (1.5,6.4);

\node at (1,5.5) {${\cal L}_1$};
\draw[thick] (1,5.5) circle (1);

\fill[bright-gray] (0.5,9) .. controls (0.7,8) and (1.3,8) .. (1.5,9);
\draw[thick] (0.5,9) .. controls (0.7,8) and (1.3,8) .. (1.5,9);

\node at (3.5,5.5) {$+$};

\node at (1,-0.5) {$\frac{1}{D_2}\, J^{{\cal L}_1}J^{{\cal L}_2}$};

\end{scope}

\begin{scope}[shift={(31.5,0)}]

\node at (-3,5.5) {\small $(-q)^{n/2}P_m$};


\fill[bright-gray] (0.5,10) rectangle (1.5,6.4);
\draw[thick] (0.5,10) -- (0.5,6.35);
\draw[thick] (1.5,10) -- (1.5,6.35);

\fill[bright-gray] (0.5,4.6) .. controls (0.7,3.5) and (1.3,3.5) .. (1.5,4.6);
\draw[thick] (0.5,4.6) .. controls (0.7,3.5) and (1.3,3.5) .. (1.5,4.6);

\node at (1,5.5) {${\cal L}_1$};
\draw[thick] (1,5.5) circle (1);

\fill[bright-gray] (0.5,2) .. controls (0.7,3) and (1.3,3) .. (1.5,2);
\draw[thick] (0.5,2) .. controls (0.7,3) and (1.3,3) .. (1.5,2);

\node at (3.5,5.5) {$+$};

\node at (1,0) {$\frac{1}{D_2}\, J^{{\cal L}_1}J^{{\cal L}_2}$};

\end{scope}

\begin{scope}[shift={(39.5,0)}]

\node at (-2,5.5) {\small $P_n P_m$};

\fill[bright-gray] (0.5,4.6) .. controls (0.7,3.5) and (1.3,3.5) .. (1.5,4.6);
\draw[thick] (0.5,4.6) .. controls (0.7,3.5) and (1.3,3.5) .. (1.5,4.6);

\fill[bright-gray] (0.5,6.4) .. controls (0.7,7.5) and (1.3,7.5) .. (1.5,6.4);
\draw[thick] (0.5,6.4) .. controls (0.7,7.5) and (1.3,7.5) .. (1.5,6.4);

\fill[bright-gray] (0.5,9) .. controls (0.7,8) and (1.3,8) .. (1.5,9);
\draw[thick] (0.5,9) .. controls (0.7,8) and (1.3,8) .. (1.5,9);

\node at (1,5.5) {${\cal L}_1$};
\draw[thick] (1,5.5) circle (1);

\fill[bright-gray] (0.5,2) .. controls (0.7,3) and (1.3,3) .. (1.5,2);
\draw[thick] (0.5,2) .. controls (0.7,3) and (1.3,3) .. (1.5,2);

\node at (1,0) {$J^{{\cal L}_1}J^{{\cal L}_2}$};
    
\end{scope}

\end{tikzpicture}
\end{picture}
    \caption{\footnotesize On the left, only the two-strand tangle of a link ${\cal L}_{n,m}$ is shown. The two explicitly depicted braids are separated by some intermediate link ${\cal L}_1$. On the right-hand side, of the equality the Kauffman bracket is applied to these two braids (see Fig.\,\ref{fig:Kauff-braid}). In the second and third terms, one can see that the same composite links (the connected sum of two links) ${\cal L}_1 \# {\cal L}_2$ appear. In the fourth term, the disjoint union of these links ${\cal L}_1 \cup {\cal L}_2$ is obtained. The Jones polynomials for all the resulting links are indicated below.}
    \label{fig:sum-par-2}
\end{figure}

\noindent Let us first consider the case $k=2$, see Fig.\,\ref{fig:sum-par-2}. After resolving the braids with $n$ and $m$ crossings, we obtain (from left to right) some link ${\cal L}_0$, the connected sum ${\cal L}_1 \# {\cal L}_2$ and the disjoint union ${\cal L}_1 \cup {\cal L}_2$. Using the formulas below for the Jones polynomials of the connected and disjoint sums of links,
\begin{equation}
\begin{aligned}
    J^{{\cal L}_1 \# {\cal L}_2} &= \frac{1}{D_2}\, J^{{\cal L}_1}J^{{\cal L}_2}\,, \\
    J^{{\cal L}_1 \cup {\cal L}_2} &= J^{{\cal L}_1}J^{{\cal L}_2}\,,
\end{aligned}
\end{equation}
then, for the Jones polynomial of the initial link we obtain:
\begin{equation}
    J^{{\cal L}_{n,m}} = (-q)^{(n+m)/2} J^{{\cal L}_0} + \frac{1}{D_2}\, J^{{\cal L}_1}J^{{\cal L}_2} \underbrace{\left( (-q)^{m/2}P_n + (-q)^{n/2}P_m + P_n P_m D_2 \right)}_{P_{n+m}}  = (-q)^{(n+m)/2} J^{{\cal L}_0} + \frac{1}{D_2}\, J^{{\cal L}_1}J^{{\cal L}_2} P_{n+m}\,,
\end{equation}
which is achieved thanks to an important property
\begin{equation}\label{PSumFormula}
    \boxed{(-q)^{m/2}P_n + (-q)^{n/2}P_m + P_n P_m D_2 = P_{n+m}\,,}
\end{equation}
which can be easily verified either from the known expression for $P_n$~\eqref{Step3} or derived directly from the Kauffman bracket for the two-strand braid with $n+m$ crossings. Indeed, in this case, the Jones polynomial depends only on the sum $n+m$.

We can similarly verify that now for $k=3$ we obtain
\begin{equation}
\begin{aligned}
    J^{{\cal L}_{n,m,p}} &= (-q)^{(n+m+p)/2} J^{{\cal L}_0} + \frac{1}{D_2^2} \, J^{{\cal L}_1}J^{{\cal L}_2}J^{{\cal L}_3}\Big(\underbrace{(-q)^{n/2} P_m P_p D_2 + (-q)^{(n+p)/2} P_m + (-q)^{(m+n)/2} P_p}_{(-q)^{n/2}P_{m+p}} + (-q)^{(m+p)/2} P_n +\\ 
    &+ \underbrace{(-q)^{m/2} P_n P_p D_2 + (-q)^{p/2} P_n P_m D_2 + P_n P_m P_p D_2^2}_{P_{m+p}P_n D_2} \Big) = (-q)^{(n+m+p)/2} J^{{\cal L}_0} + \frac{1}{D_2^2} \, J^{{\cal L}_1}J^{{\cal L}_2}J^{{\cal L}_3} P_{n+m+p}\,,
\end{aligned}
\end{equation}
and for an arbitary $k$, it can be proven by induction that
\begin{equation}
    J^{{\cal L}_{n_1,n_2,\dots,n_k}} = (-q)^{\sum_j n_j/2} J^{{\cal L}_0} + \frac{1}{D_2^{k-1}} \left(\prod_{l=1}^k J^{{\cal L}_l}\right) P_{n_1+\dots+ n_k}\,.
\end{equation}
Thus, we have found that, in the general case as well, the Jones polynomial depends only on the sum of the braid parameters, and the resulting coefficients $(-q)^{\sum_j n_j/2}$ and $P_{n_1+\dots+ n_k}$ are the same as in the Goeritz method.

\section{Examples of calculation of Jones polynomials}
\label{sec:Jones-ex}

In this section, we provide examples of calculation of Jones polynomials using Goeritz matrices.

\setcounter{equation}{0}

\subsection{Torus knots $T[2,p]$}

\begin{figure}[h!]
		\centering	
		\includegraphics[width =0.65\linewidth]{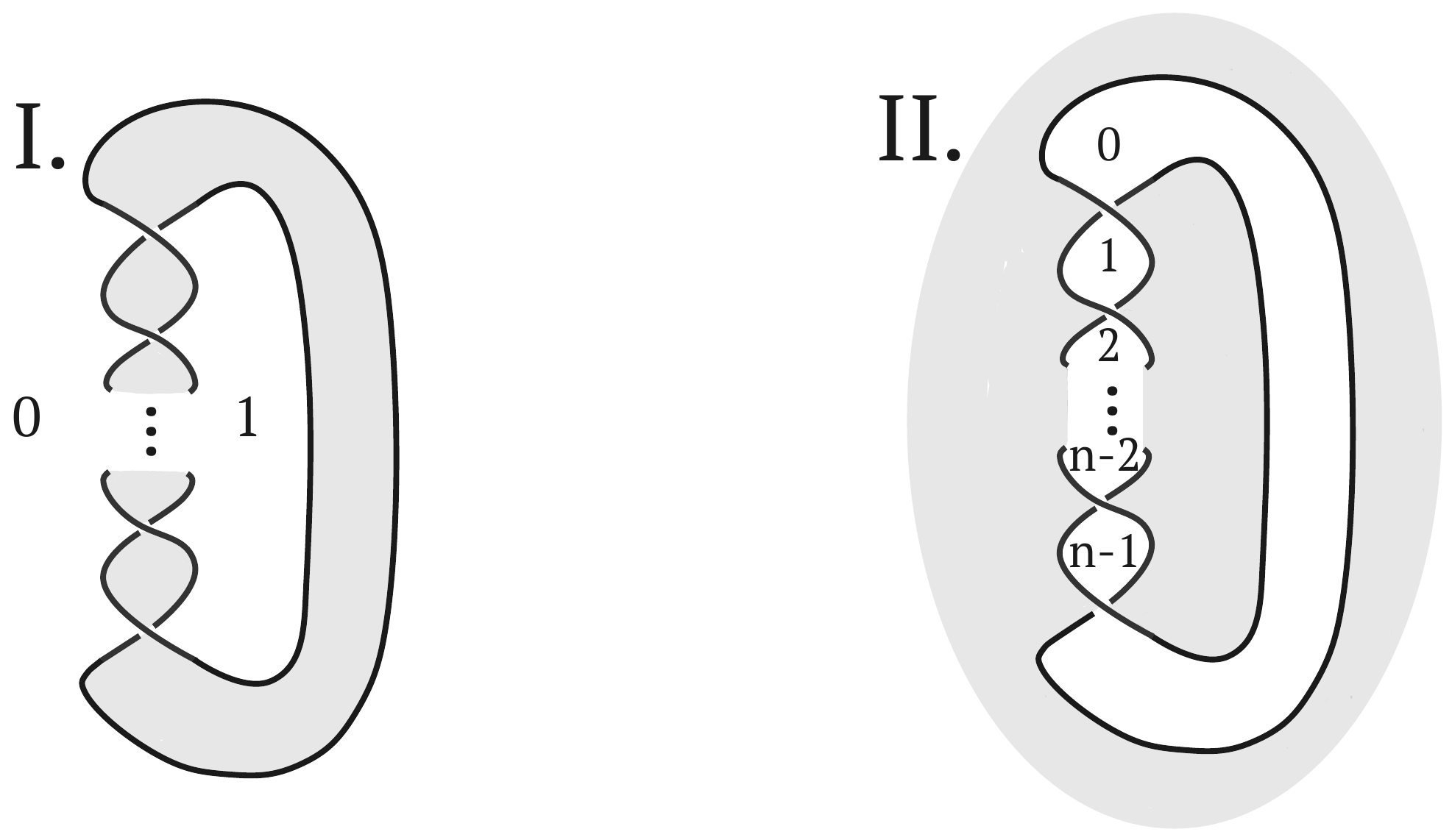}
        \caption{\footnotesize Two checkerboard coloring for two-strand links.}
        \label{fig:Tor2n}
	\end{figure}

Let us consider the family of torus knots $T[2,p]$ for an arbitrary $p$. Two variants of the checkerboard coloring lead to different Goeritz matrices, as illustrated in Fig.~\ref{fig:Tor2n}: in the first case, the outer region is white,  labeled as $0$ (there are only two white regions in total, hence, the Goeritz matrix is one-dimensional); in the second case, the outer region is black (there are $p$ white regions in total and the corresponding matrix has size is $(p-1)\times(p-1)$).

\paragraph{The first coloring variant.}

The first coloring variant is simpler, since the Goeritz matrix consists of a single number $p$ being the number of crossings between regions 0 and 1:
\begin{equation}\label{Goeritz1Tor2n}
\widetilde{G}^{T[2,p]}_I  =
\begin{pmatrix}
 \begin{array}{cc}  -p & p \\  p & -p       
 \end{array}
 \end{pmatrix}\quad \mapsto \quad G^{T[2,p]}_I  = -p\,.
\end{equation}
Then, immediately find Kauffman bracket from \eqref{Step2}\,:

\begin{equation}\label{GoeritzFunctionForT2p}
   \mu[G^{T[2,p]}_I ] = \mu[-p]=(-q)^{p/2}D_2^2 + P_p(q)D_2 = (-q)^{p/2}(q^{-2}+1+q^2+(-1)^pq^{-2p})\,.
\end{equation}
The well-known result for the Jones polynomials of two-strand knots is easily obtained using the normalization~\eqref{Jones}:

\begin{equation}\label{JonesT2pTrue}
    J^{T[2,p]}(q) = (-(-q)^{3/2})^{-p}\mu [G_{I}^{T[2,p]}]  =q^{-p} (q^{-2}+1+(-1)^pq^{-2p}+q^2)\,.
\end{equation}

\paragraph{The second coloring variant.}

In this coloring variant, the number of white regions coincides with the number of crossings $p$ in the torus knot, consequently, the Goeritz matrix $G_{II}^{T[2,p]}$ has size $(p-1)\times(p-1)\,$:

\begin{equation}\label{GoeritzMatrixTorKont2}
\widetilde{G}_{II}^{T[2,p]}  =
\begin{pmatrix}
2 & -1 & 0 & \dots  & 0 &-1\\
-1 & 2 & -1 & \ddots &\dots & 0 \\
0 & -1 & \ddots & \ddots &\ddots &\vdots\\

\vdots  & \ddots & \ddots&\ddots & -1 & 0 \\

0  &\dots &\ddots & -1 & 2 & -1 \\
-1 & 0 & \dots & 0 &-1 & 2  \\

\end{pmatrix}\quad  \mapsto \quad G_{II}^{T[2,p]}=\begin{pmatrix}
2 & -1 & 0 & \dots  & 0 \\
-1 & 2 & -1 & \ddots &\vdots  \\
0 & -1 & \ddots & \ddots &0 \\

\vdots  & \ddots & \ddots&\ddots & -1  \\

 0 & \dots & 0 &-1 & 2  \\
\end{pmatrix}. 
\end{equation}
Due to the coincidence of the off-diagonal elements and their symmetric arrangement with respect to the main diagonal, all transformations of the function $\mu[G_{II}^{T[2,p]}]$ at the first step of the Boninger algorithm are identical to each other:

\begin{equation}\label{GoeritzCalculatingTor2p}
\mu [G_{II}^{T[2,p]}] = (-q)^{-1/2} \mu \begin{bmatrix}\begin{pmatrix}
    \begin{array}{c|cc|c}
\multirow{2}{*}{\begin{tabular}{c}$G_{II}^{T[2,k]}$\end{tabular}} & \cdots  & \cdots & \multirow{2}{*}{\textbf{\Large{0}}} \\
 & -1  & \cdots &  \\
\hline
0 \ \ \cdots -1 & 1 & 0 & \cdots \\
\cdots  & 0 & 1 & -1 \ \cdots \ \ \ 0 \\
\hline
\multirow{2}{*}{\textbf{\Large{0}}} & \cdots  & -1 & \multirow{2}{*}{\begin{tabular}{c}$G_{II}^{T[2,p-k-2]}$\end{tabular}} \\
 & \cdots  & \cdots & \\
\end{array}\end{pmatrix} \end{bmatrix} + P_{-1}(q)\mu [G_{II}^{T[2,p-1]}]\,.
\end{equation}

Then, the combinatorics of the first step of the Boninger algorithm is binomial. Taking into account that zeros on the diagonals are equivalent to multiplying the function $\mu[G]$ by the quantum two $D_2 = q + q^{-1}$, we obtain the following expression for $\mu [G_{II}^{T[2,p]}]$\footnote{We present the function in a general form so that it also reproduces the result for mirror links $\overline{L}$. Obviously, the Goeritz matrix for them is equal to the negative of the matrix~\eqref{GoeritzMatrixTorKont2}: $G_{II}^{\overline{T[2,p]}} = - G_{II}^{T[2,p]}$.}:

\begin{equation}\label{ExampleTorKnot2}
    \mu [G_{II}^{T[2,p]}] = \sum_{k=0}^{|p|-3}\left\{\begin{pmatrix}
 \begin{array}{c}  |p|-2  \\  k     
 \end{array}
 \end{pmatrix}(-q)^{-{\rm sgn}(p)(|p|-2-k)/2}P_{-{\rm sgn}(p)}^{k}D_2^{|p|-3-k}\right\}\mu\begin{bmatrix}\begin{pmatrix}  1 & 0\\0 & 1     
 \end{pmatrix}
     
 \end{bmatrix} + P_{-{\rm sgn}(p)}^{|p|-2}\mu[2]\,.
\end{equation}

This expression reproduces the values of~\eqref{GoeritzFunctionForT2p} for $|p| \ge 3$, since matrices of the form~\eqref{GoeritzMatrixTorKont2} appear for two-strand knots starting from three crossings. Thus, by a somewhat longer route, we have once again reproduced the value of the Jones polynomial for two-strand knots:

\begin{equation}\label{JonesT2p2}
    J^{T[2,p]}(q)  = (-(-q)^{3/2})^{-p}\mu [G_{II}^{T[2,p]}]\,.
\end{equation}

\subsection{Twist knots Tw$_{m}$}

\begin{figure}[h!]
		\centering	
		\includegraphics[width =0.75\linewidth]{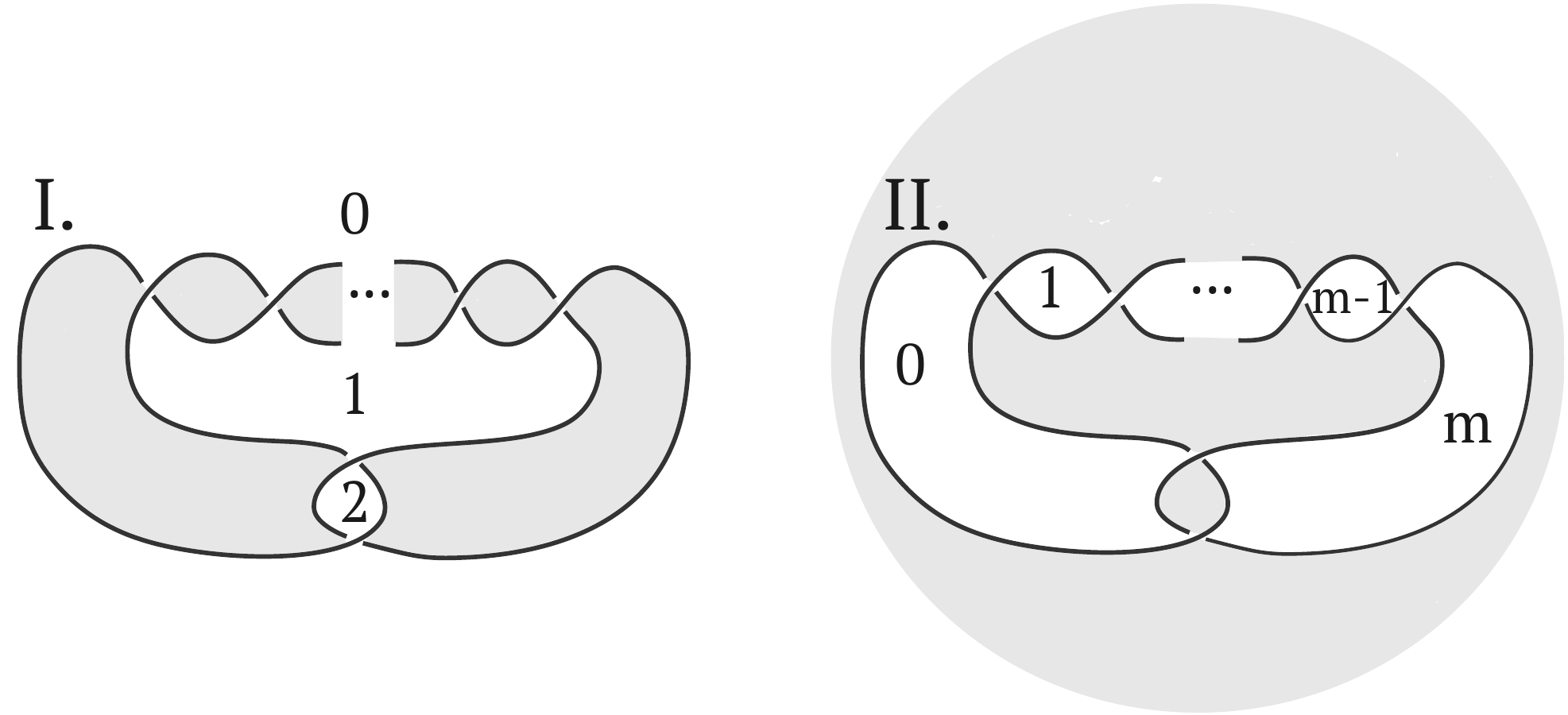}
        \caption{\footnotesize Two checkerboard coloring for twist knots ${\rm Tw}_m$ with positive $m$.}
        \label{fig:TwistKnot}
	\end{figure}

Another example convenient for demonstrating the Goeritz method is the one-parameter family of twist knots ${\rm Tw}_{m}$\footnote{By default, we assume $m$ to be positive. For the first coloring, this is not essential, since $m$ explicitly enters the corresponding Goeritz matrix~\eqref{Goeritz1Twist2m}, while the second coloring undergoes changes — among all the off-diagonal elements, only $\Tilde g_{0m}$ does not change sign. Consequently, there will be minor modifications in the computation of the Jones polynomial, but all steps of the algorithm remain completely analogous to~\eqref{ExampleTwist2}. We also note that a negative $m$ does not correspond to a mirror link.} (see Section~\ref{sec:TwistKnots}). Similar to the example with torus knots in the first coloring, we choose the outer region to be white, and in the second one, we choose it to be black (see Fig. ~\ref{fig:TwistKnot}).

\paragraph{The first coloring variant.}

With this coloring, we have three white regions and therefore the Goeritz matrix is of size $2\times2$:
\begin{equation}\label{Goeritz1Twist2m}
\widetilde{G}^{{\rm Tw}_{m}}_I  =
\begin{pmatrix}
 \begin{array}{ccc}  -m-1 & m & 1 \\ 
 m & -m-1 & 1 \\ 
 1 & 1 & -2      
 \end{array}
 \end{pmatrix} \quad \mapsto \quad G^{{\rm Tw}_{m}}_I  = \begin{pmatrix}
 \begin{array}{cc}  -m-1 & 1 \\  1 & -2       
 \end{array}
 \end{pmatrix}.
\end{equation}
Considering the only non-diagonal element, we obtain:

\begin{equation}\label{Twist1Step1}
\mu\left[G^{{\rm Tw}_{m}}_I\right] = (-q)^{\frac{1}{2}} \mu \begin{bmatrix}
 \begin{pmatrix}
     -m  & 0 \\  0 & -1  
 \end{pmatrix}     
 \end{bmatrix} + P_1(q)\,\mu [-m-1]\,,
\end{equation}

\begin{equation}\label{Twist1Step2}
\mu\left[G^{{\rm Tw}_{m}}_I\right] = (-q)^{\frac{1}{2}}D_2 \left((-q)^{-\frac{m}{2}}D_2 + P_{-m}(q)\right)\left((-q)^{-\frac{1}{2}}D_2 + P_{-1}(q)\right) + P_1(q)D_2\left((-q)^{-\frac{m+1}{2}}D_2 + P_{-m-1}(q)\right)\,,
\end{equation}

\begin{equation}\label{Twist1Step3}
\mu\left[G^{{\rm Tw}_{m}}_I\right] = (-q)^{m/2} (q^{-2}-(-q^2)^{-m}(q^{-2}+q^2)-q^4)\,.
\end{equation}
Note that depending on the parity of $m$, the sign of the crossings in the lower tangle changes. Consequently, the contribution to the normalization differs by the corresponding factor. Therefore, the values of the Jones polynomials for even $(W = -m + 2)$ and odd $(W = -m - 2)$ values of $m$ generally differ:

\begin{equation}\label{JonesTwist} 
\left\{\begin{split} 
J^{{\rm Tw}_{m}}(q)  = q^{2m-3}(q^4+q^{-2m}(q^{-2}+q^2)-q^{-2})\,,\quad \text{ even } m, \ \ \ \\
J^{{\rm Tw}_{m}}(q)  = q^{2m+3}(q^{-2}+q^{-2m}(q^{-2}+q^2)-q^{4})\,,\quad \ \ \text{ odd } m. 
\end{split}\right.
\end{equation}

\paragraph{The second coloring variant.}

There are $m+1$ white regions in this case. The relative arrangement of the regions is similar to that in the second coloring of the two-strand knot, except that region 0 is divided into two regions (0 and $m$), and there are two crossings between 	the two new regions:

\begin{equation}\label{GoeritzMatrixTwist}
\widetilde{G}^{{\rm Tw}_{m}}_{II}  =
\begin{pmatrix}
3 & -1 & 0 & \dots  & 0 &-2\\
-1 & 2 & 1 & \dots & & 0 \\
0 & -1 & \ddots &  & &\vdots\\

\vdots  & \vdots & &\ddots & -1 & 0 \\

0  & & & -1 & 2 & -1 \\
-2 & 0 & \dots & 0 & -1 & 3  \\

\end{pmatrix} \quad \mapsto \quad G^{{\rm Tw}_{m}}_{II} =\begin{pmatrix}
2 & -1 & 0 & \dots  & 0\\
-1 & 2 & -1 & \dots  \\
0 & -1 & \ddots &  &\vdots\\

\vdots  & \vdots & & \ddots & -1 \\
0 & & \dots & -1 & 3  \\

\end{pmatrix}. 
\end{equation}
It is particularly convenient that the matrix~\eqref{GoeritzMatrixTwist} has a structure similar to that of~\eqref{ExampleTorKnot2}, which makes the expression for twist knots analogous to that for two-strand knots:

\begin{equation}\label{ExampleTwist2}
    \mu [G^{{\rm Tw}_{m}}_{II}] = \sum_{k=0}^{m-2}\left\{\begin{pmatrix}
 \begin{array}{c}  m-1  \\  k     
 \end{array}
 \end{pmatrix}(-q)^{-\delta\frac{m-1-k}{2}}P_{-\delta}^{k}D_2^{m-2-k}\right\}\mu\begin{bmatrix}\begin{array}{cc}  \delta & 0\\0 & 2\delta     
 \end{array}
     
 \end{bmatrix} + P_{-\delta}^{m-1}\mu[3\delta]\,.
\end{equation}
It is easy to verify the equivalence of this formula and~\eqref{Twist1Step3}. Hence the results for the twist Jones polynomials~\eqref{JonesTwist} are correctly reproduced. For ordinary knots ${\rm Tw}_m$, $\delta = 1$, and $\delta = -1$ for their mirror counterparts $\overline{{\rm Tw}}_m$.

\section{Bipartite HOMFLY--PT polynomials from the quaternary Goeritz matrix}
\label{sec:HOMFLYfromGoeritz}

\setcounter{equation}{0}

In the works~\cite{anokhina2024planar,anokhina2025planar,anokhina2025bipartite}, the Kauffman bracket, originally introduced as a method for computing the Jones polynomial, was extended to a broader framework that allows the computation of HOMFLY--PT polynomials for a particular class of bipartite knots. These are knots that can be represented as configurations of doubled crossings with antiparallel orientation (see Fig.\,\ref{fig:pladeco}). This observation naturally suggests that Goeritz method might be generalized to the HOMFLY--PT case, at least within the class of bipartite knots. Exploring this possibility constitutes the main goal of the present work. It becomes necessary to modify both the Goeritz matrix and the Boninger algorithm to realize such a generalization. Taking into account the increased number of distinct crossings in a bipartite diagram (four in total) and the relations among them (see Appendix~\ref{sec:App1}) we derive a set of rules governing the checkerboard coloring and the filling of entries in the analogue of the Goeritz matrix that we call the quaternary Goeritz matrix $\mathcal{G}$. We also describe the corresponding transformation rules for these matrices in Subsection~\ref{sec:TransformationTripleMatrix}. Building upon the intuition developed in the analysis of the Boninger algorithm for computing the Jones polynomial (Sections~\ref{sec:muInv} and~\ref{sec:ExceptionalCases}), we introduce a function on quaternary Goeritz matrices $\mathcal{M}[\mathcal{G}]$, which taking into account the normalization condition~\eqref{HOMFLYNorm}, yields the HOMFLY--PT polynomial (see Section~\ref{sec:HOMFLYfromTripleGoeritz}).

\subsection{Quaternary Goeritz matrix}
\label{sec:quadruple-Goeritz}

In order to generalize the Goeritz matrix to the bipartite case, it is necessary to determine how the new matrix will account for the bipartite crossings within two-strand braids, which consist of four distinct types of locks (see Fig.\,\ref{fig:pladeco}). It is easy to observe that all bipartite crossings grouped within the same two–strand braid commute with one another under the Kauffman bracket resolution. In other words, a link whose diagram contains two-strand tangles formed by successively connected lock elements has the same HOMFLY--PT polynomial as any link whose diagram differs from it only by an arbitrary permutation of these lock elements within the two–strand tangles. Therefore, without loss of generality, we may restrict our consideration to diagrams composed of two-strand tangles that consist of $x$ positive vertical crossings, $\widetilde{y}$ positive horizontal crossings, $\overline{z}$ negative vertical crossings, and $\widehat{t}$ negative horizontal crossings (see Fig.\,\ref{fig:SignHOMFLYPT}). For the analogue of the Goeritz matrix $\mathcal{G}$, this implies that we must independently record the number of different types of locks appearing in the two–strand braids into which the link diagram is decomposed.

\subsubsection{Definition and examples}\label{sec:QuadrupleGoeritzParametrs}

Based on the above considerations, we introduce a new object --- the quaternary Goeritz matrix, which depends on four independent parameters (the numbers $x, \widetilde{y}, \overline{z}, \widehat{t}\,$). We define this matrix according to the following procedure, which is analogous to the algorithm described in Section~\ref{sec:Goeritz-mat}.

\begin{enumerate}

\item We assign a checkerboard coloring to a bipartite knot diagram according to the standard rules, but with bipartite vertices. All white regions are then numbered by integers.

\item  The unreduced quaternary Goeritz matrix $\widetilde{\mathcal{G}} = (\widetilde{\mathfrak{g}}_{ij})$ is defined in the same way as the ordinary Goeritz matrix \eqref{IrGoerMatr}. The only difference lies in the contribution assigned to a single crossing $c$, this contribution is now one of the four values $\sigma(c) = \{ -1,\widehat{1},-\widetilde{1}, \overline{1}\}$
\footnote{We note that the two positive bipartite locks correspond to negative values of the function $\sigma(c)$: $(-1,-\widetilde{1})$, whereas the two negative locks correspond to positive values $(\overline{1}, \widehat{1})$. As will become clear later, information about signs is not so valuable from the point of view of calculating the HOMFLY-PT polynomial, since our algorithm counts each type of bipartite crossings separately. Their purpose is to allow us to pass from bipartite diagrams to their precursor diagrams (see Section~\ref{sec:BonNorm}) and to recover the Jones polynomial of these precursors via the reduction of the Goeritz method for bipartite HOMFLY--PT polynomials $\big(\phi \rightarrow -q,\; \bphi \rightarrow -q^{-1},\; D_N \rightarrow D_2\big)$ together with the sign convention for the checkerboard coloring $(-1 \to -1;\, \widehat{1} \to 1;\, -\widetilde{1} \to -1;\, \overline{1} \to 1)$.} (see Fig.\,\ref{fig:SignHOMFLYPT}):

\begin{equation}\label{FirstGoeritzHOMFLYMatrix}
\centering
\widetilde{\mathfrak{g}}_{ij} = \left\{\begin{split}
&-\sum_{c_{ij}}\sigma(c_{ij}),  \ \  &i\ne j\,,  \ \\
&-\sum_{k \ne i} \widetilde{\mathfrak{g}}_{kj},  & i = j\,, \\ 
\end{split}\right.
\end{equation}
where $c_{ij}$ is the set of all crossings between the regions $i$ and $j$, and  
$x_{ij},\, \widetilde{y}_{ij},\, \overline{z}_{ij},\, \widehat{t}_{ij}$ denote the numbers of bipartite crossings of each type between these regions, while $\sigma(c_{ij}) = x_{ij} +\widetilde{y}_{ij}+\overline{z}_{ij}+\widehat{t}_{ij}$.

\begin{figure}[h!]
		\centering	
		\includegraphics[width =0.65\linewidth]{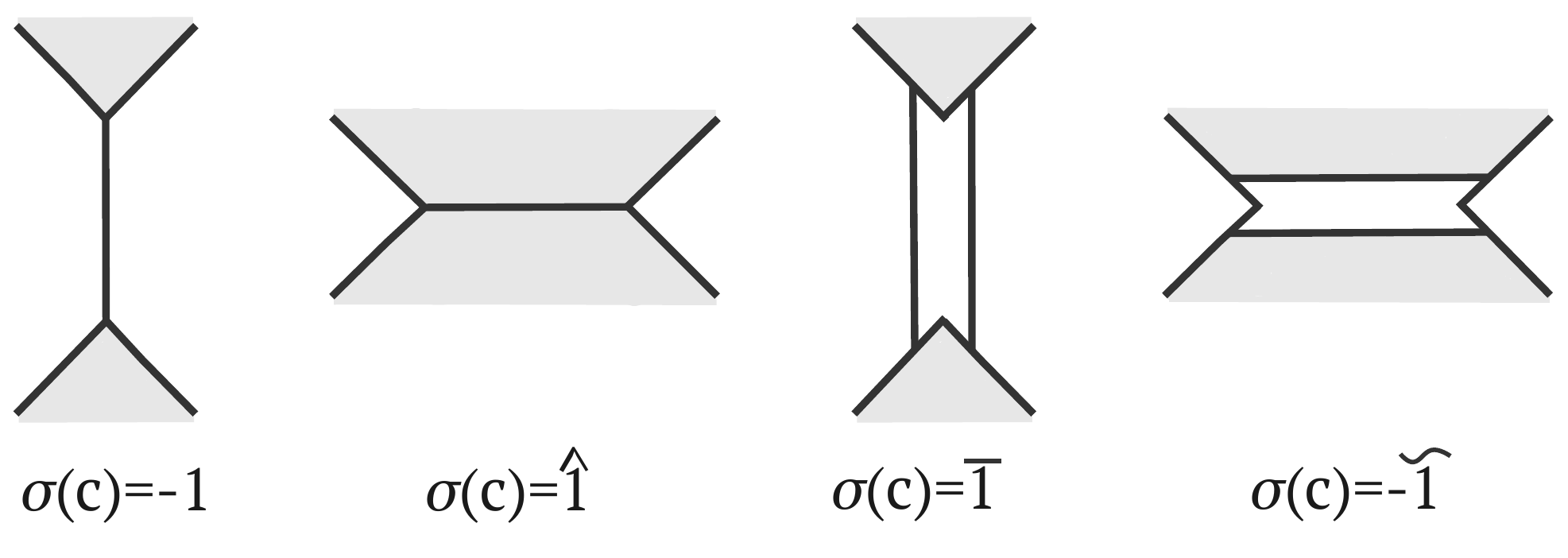}
        \caption{\footnotesize Checkerboard sign convention for bipartite diagrams. From left to right: positive vertical lock, negative horizontal lock, negative vertical lock, positive horizontal lock.}
        \label{fig:SignHOMFLYPT}
	\end{figure}

\item The (reduced) quaternary Goeritz matrix $\mathcal{G} = (\mathfrak{g}_{ij})$ is again obtained from the unreduced matrix by deleting the $k$-th row and the $k$-th column (for any choice of $k$). The resulting matrix then consists of the remaining entries $\mathfrak{g}_{ij} = \widetilde{\mathfrak{g}}_{ij},\; \forall i\ne k,\; j\ne k$. 

\end{enumerate}

The resulting quaternary Goeritz matrix evidently can be decomposed into a sum of four symmetric matrices, each depending only on one of the parameters from $(x,\, \widetilde{y},\, \overline{z},\, \widehat{t}\,)$.

\begin{equation}
    \mathcal{G}(x,\widetilde{y}, \overline{z}, \widehat{t} \ ) = \mathcal{G}(x)+\mathcal{G}(\widetilde{y} ) + \mathcal{G}(\overline{z}) + \mathcal{G}( \widehat{t} \ )\,.
\end{equation}

\subsubsection{Transformation of quaternary Goeritz matrices}\label{sec:TransformationTripleMatrix}

The quaternary Goeritz matrix inherits from the ordinary Goeritz matrix the structure of transformations acting on it, since these matrix transformations correspond to changes in diagrams arising from resolving two-strand braids located between two white regions using the Kauffman bracket (merging or splitting regions). The same correspondence carries over to the case of bipartite diagrams.
\medskip

\noindent {\bf (I)} Recall that the first transformation in~\eqref{GoeritzMovie1} corresponds to separating the regions $i$ and $j$ into regions without crossings and this property is preserved for quaternary Goeritz matrices. In matrix terms, this transformation amounts to setting the corresponding off-diagonal entry $\mathfrak{g}_{ij}$ to zero:

\begin{equation}\label{QuadroGoeritzMovie1}
\mathcal{G}=\begin{pmatrix}
\ddots & \\
& \mathfrak{g}_{ii} & \dots & \boxed{\mathfrak{g}_{ij}} \\
& \vdots & \ddots & \vdots \\
& \mathfrak{g}_{ji} & \dots & \mathfrak{g}_{jj} \\
& & & & \ddots \\
\end{pmatrix}\quad \mapsto \quad
\mathcal{G}^{(I)}=\begin{pmatrix}
\ddots & \\
& \mathfrak{g}_{ii} + \mathfrak{g}_{ij} & \dots & 0 \\
& \vdots & \ddots & \vdots \\
& 0 & \dots & \mathfrak{g}_{jj} + \mathfrak{g}_{ij} \\
& & & & \ddots \\
\end{pmatrix}.   
\end{equation}

\noindent The explicit form of the first transformation for the entry $\mathfrak{g}_{i j} = x_{ij}+\widetilde{y}_{ij}+\overline{z}_{ij}+\widehat{t}_{ij}$ is given by:

\begin{equation}
\begin{aligned}
\mathfrak{g}_{i i}   \quad &\mapsto \quad \mathfrak{g}_{i i}+\mathfrak{g}_{i j} = x_{ii} +x_{ij} +\widetilde{y}_{ii} +\widetilde{y}_{ij}+\overline{z}_{ii}  +\overline{z}_{ij}+\widehat{t}_{ii}+\widehat{t}_{ij}\,, \\
\mathfrak{g}_{j j} \quad &\mapsto \quad \mathfrak{g}_{j j}+\mathfrak{g}_{i j} = x_{jj} +x_{ij} +\widetilde{y}_{jj} +\widetilde{y}_{ij}+\overline{z}_{jj}  +\overline{z}_{ij}+\widehat{t}_{jj}+\widehat{t}_{ij}, \\
\mathfrak{g}_{i j}, \mathfrak{g}_{j i} \quad &\mapsto \quad 0\,.
\end{aligned}
\end{equation}

\noindent {\bf (II)} Similarly, the second transformation $\mathcal{G} \;\mapsto\; \mathcal{G}^{(II)}$ corresponds to merging the regions $i$ and $j$ in a bipartite diagram. As in \eqref{GoeritzMovie2}, the merging of regions reduces the size of a matrix by one, since the $j$-th row and the $j$-th column are removed:

\begin{equation}\label{QuadroGoeritzMovie2}
\begin{picture}(300,100)(90,-45)
    \put(0,0){\mbox{$\mathcal G=\begin{pmatrix}
\ddots & & & & \vdots \\
& \mathfrak{g}_{ii} & \mathfrak{g}_{ik} & \dots & \boxed{\mathfrak{g}_{ij}} \\
& \mathfrak{g}_{ki} & \ddots & & \mathfrak{g}_{kj} \\
& \vdots & & \ddots & \vdots \\
\dots & \mathfrak{g}_{ji} & \mathfrak{g}_{jk} & \dots & \mathfrak{g}_{jj} & \dots \\
& & & & \vdots &  \ddots \\
\end{pmatrix} \quad \mapsto \quad \mathcal G^{(II)} =
\begin{pmatrix}
\ddots & & & & {\vdots} \\
& \mathfrak{g}_{ii} + \mathfrak{g}_{jj} + 2\mathfrak{g}_{ij} & \mathfrak{g}_{ik} + \mathfrak{g}_{jk} & \dots & \mathfrak{g}_{ij} \\
& \mathfrak{g}_{ki} + \mathfrak{g}_{kj} & \ddots & & \mathfrak{g}_{kj} \\
& \vdots & & \ddots & \vdots \\
{\dots} & \mathfrak{g}_{ji} & \mathfrak{g}_{jk} & \dots & \mathfrak{g}_{jj} & {\dots} \\
& & & & {\vdots} &  \ddots \\
\end{pmatrix}. $}}

\thicklines

{\color{red}

\put(420,-53){\line(0,1){110}}
\put(237,-25){\line(1,0){225}}

}

\end{picture}
\end{equation}
The diagonal entry $\mathfrak{g}_{ii}$ and the off-diagonal entries in the $i$-th row and column $\mathfrak{g}_{ik}, \mathfrak{g}_{ki}$ are modified due to the merging of the crossings between regions $i$ and $k$ with those between regions $j$ and $k$, since the regions $i$ and $j$ are combined into a single region:

\begin{equation}
\begin{aligned}
& \mathfrak{g}_{i i} \quad \mapsto \quad \mathfrak{g}_{i i}+\mathfrak{g}_{j j}+2 \mathfrak{g}_{i j} = x_{ii}+x_{jj} +2x_{ij} +\widetilde{y}_{ii} +\widetilde{y}_{jj} +2\widetilde{y}_{ij}+\overline{z}_{ii}  +\overline{z}_{jj}  +2\overline{z}_{ij}+\widehat{t}_{ii}+\widehat{t}_{jj}+2\widehat{t}_{ij}\,, \\
& \mathfrak{g}_{i k} \quad \mapsto \quad \mathfrak{g}_{i k}+\mathfrak{g}_{ki} = x_{ik} + x_{ki} +\widetilde{y}_{ik}+\widetilde{y}_{ki}+\overline{z}_{ik}+\overline{z}_{ki}+\widehat{t}_{ik}++\widehat{t}_{ki}, \; \text { for all} \; k \neq i\,, \\
& \mathfrak{g}_{k i} \quad \mapsto \quad \mathfrak{g}_{k i}+\mathfrak{g}_{k j} = x_{ik} + x_{ki} +\widetilde{y}_{ik}+\widetilde{y}_{ki}+\overline{z}_{ik}+\overline{z}_{ki}+\widehat{t}_{ik}++\widehat{t}_{ki}, \; \text { for all } \; k \neq i\,.
\end{aligned}
\end{equation}

\noindent {\bf (III)} The third transformation $\mathcal{G} \;\mapsto\; \mathcal{G}^{(III)}$ for diagonal matrices again corresponds, as in the case of ordinary Goeritz matrices, to a diagram in which all regions except one intersect only with a single other region (see Fig.~\ref{fig:diag-Goeritz-ex}). The latter region corresponds to the row and column removed when reducing the quaternary Goeritz matrix. Continuing the analogy, a planar resolution of all crossings between any two regions in these diagrams produces identical diagrams up to unknots, which can be replaced by the factors $D_N$. Consequently, the third transformation must account for the removal of one region, allowing the $i$-th row and column to be deleted:

\begin{equation}
    \mathcal G = {\rm diag}(\dots,\, \mathfrak{g}_{i-1},\, \boxed{\mathfrak{g}_i}\,,\, \mathfrak{g}_{i+1},\, \dots) \quad \mapsto \quad \mathcal G^{(III)} = {\rm diag}(\dots ,\, \mathfrak{g}_{i-1},\, \cancel{\mathfrak{g}_i},\, \mathfrak{g}_{i+1},\, \dots)\,. 
\end{equation}

\subsection{Function $\mathcal{M}$ on quaternary Goeritz matrices and the HOMFLY-PT polynomial}\label{sec:HOMFLYfromTripleGoeritz}

\begin{figure}[h!]
		\centering	
		\includegraphics[width=0.8\linewidth]{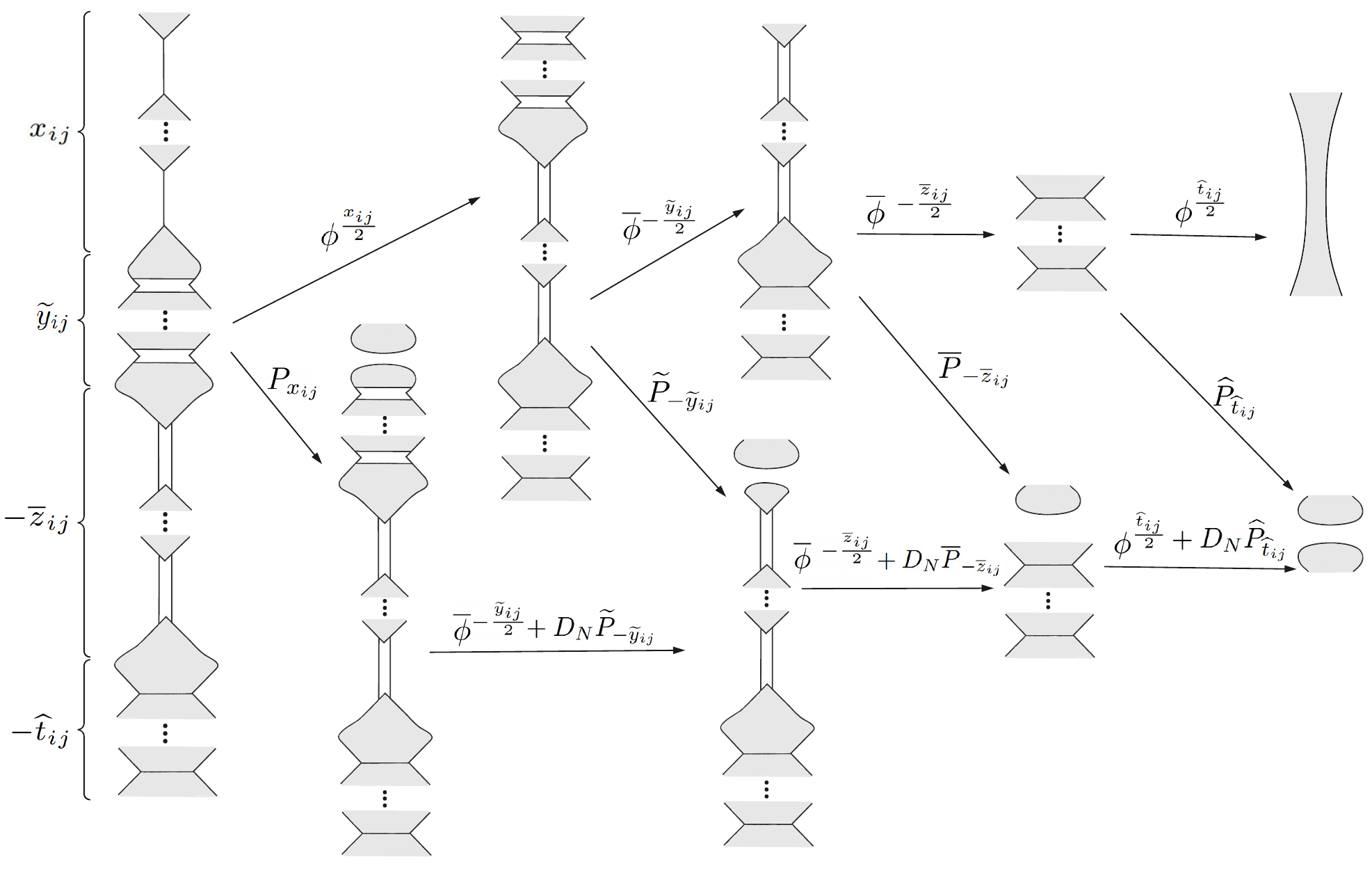}
		\caption{\footnotesize Resolution of a two–strand bipartite braid using planar decomposition. The coefficients associated with each resolution are shown above the arrows. These coefficients are consistent with the action of the function $\mathcal{M}[\mathfrak{g}_{ij}]$ on the off-diagonal entries of the quaternary Goeritz matrix. When resolving a braid corresponding to the action of $\mathcal{M}[\mathfrak{g}_{ij}]$ on diagonal entries, signs of parameters $x_{ij},\, \widetilde{y}_{ij},\, \overline{z}_{ij},\, \widehat{t}_{ij}$ are modified in transformations, as explicitly shown in~\eqref{StepHOMFLYfromGoeritz2}.} 
        \label{fig:KaufBipBraid}
	\end{figure}

Our next goal is to introduce a function $\mathcal{M}[\mathfrak{g}_{ij}]$ that performs transformations on quaternary Goeritz matrices corresponding to the planar resolution of a bipartite diagram (see Fig.\,\ref{fig:pladeco}). Using the intuition gained from resolving a two-strand braid, we have constructed the analogous function for ordinary link diagrams in Section \ref{sec:muInv}. Following the same approach, we now consider a two-strand braid composed of bipartite locks; an example of the resolution of such a braid is shown in Fig.\,\ref{fig:KaufBipBraid}. As noted earlier, bipartite diagrams can be represented as consisting of two-strand braids with sequentially arranged bipartite crossings of different types\footnote{We use the property of equality of HOMFLY-PT polynomials for diagrams that differ by an arbitrary permutation of bipartite locks within two-strand braids.}. Therefore, we focus on this general example.

As in the case of resolving an ordinary braid, the resolution of a braid composed of $n$ bipartite locks of the same type produces one vertical resolution and $n-1$ horizontal resolutions, which are identical up to the appearance of cycles corresponding to the factors $D_N$. Summing over all horizontal resolutions yields the combinatorial factors $(P_x,\, \widetilde{P}_{\widetilde{y}},\, \overline{P}_{\overline{z}},\, \widehat{P}_{\widehat{t}})$, which must also appear in a transformation of entries of the quaternary Goeritz matrix:

\begin{equation}\label{PPolynomials} 
\left\{\begin{split} 
P_{x}(\phi,\overline{\phi}, D_N) = \sum^{|x|}_{k=1} \begin{pmatrix}
|x| \\
k 
\end{pmatrix} D_N^{k-1} \phi^{|x|/2-k} \,  , \ \ \ \  \widetilde{P}_{\widetilde{y}}(\phi,\overline{\phi}, D_N) = \sum^{|\widetilde{y}|}_{k=1} \begin{pmatrix}
|\widetilde{y}| \\
k 
\end{pmatrix} D_N^{k-1} \overline{\phi}^{k-|\widetilde{y}|/2}, \\
\overline{P}_{\overline{z}}(\phi,\overline{\phi}, D_N) = \sum^{|\overline{z}|}_{k=1} \begin{pmatrix}
|\overline{z}| \\
k 
\end{pmatrix} D_N^{k-1} \overline{\phi}^{|\overline{z}|/2-k}, \ \ \ \ \widehat{P}_{\widehat{t}}(\phi,\overline{\phi}, D_N) = \sum^{|\widehat{t}|}_{k=1} \begin{pmatrix}
|\widehat{t}| \\
k 
\end{pmatrix} D_N^{k-1} \phi^{k-|\widehat{t}|/2} \,.
\end{split}\right.
\end{equation}
We want the function $\mathcal{M}[\mathfrak{g}_{ij}]$ to be expressed through combinatorial factors corresponding to the different types of bipartite crossings, since it must account for the simultaneous resolution of all bipartite locks between regions $i$ and $j$. The action of the function $\mathcal{M}[\mathfrak{g}_{ij}]$ on an off-diagonal element $\mathfrak{g}_{ij}$ becomes clear from the analysis of the resolution of a general bipartite braid (see Fig.~\ref{fig:KaufBipBraid}):

\begin{equation}\label{StepHOMFLYfromGoeritz1}
\begin{aligned}
\boxed{\mathcal{M} [\mathfrak{g}_{ij}] = \mathcal{U}_1[\mathfrak{g}_{\normalfont{ij}}] \mathcal{M}\left[\mathfrak{g}_{ij}^{(I)}\right] + \mathcal{U}_2[\mathfrak{g}_{ij}] \mathcal{M}\left[\mathfrak{g}_{ij}^{(II)}\right]\,,}
\end{aligned}
\end{equation} 
where we have introduced the following notation for brevity:
\begin{equation}
\boxed{\begin{aligned}
    \mathcal{U}_1[\mathfrak{g}_{\normalfont{ij}}] &= \phi^{(x_{ij}+\widehat{t}_{ij})/2}\overline{\phi}{}^{-(\widetilde{y}_{ij}+\overline{z}_{ij})/2}\,, \\
    \mathcal{U}_2[\mathfrak{g}_{ij}] &= \phi^{\widehat{t}_{ij}/2}\overline{\phi}{}^{-(\widetilde{y}_{ij}+\overline{z}_{ij})/2} P_{x_{ij}} + \left(\phi^{x_{ij}/2}+D_NP_{x_{ij}}\right)\cdot \\ 
    &\cdot \left[\left(\overline{\phi}{}^{-\widetilde{y}_{ij}/2}+  D_N\widetilde{P}_{-\widetilde{y}_{ij}}\right)\left(\overline{\phi}{ \ }^{-\overline{z}_{ij}/2}\widehat{P}_{\widehat{t}_{ij}}+\phi^{\widehat{t}_{ij}/2}\overline{P}_{-\overline{z}_{ij}}+D_N\overline{P}_{-\overline{z}_{ij}}\widehat{P}_{\widehat{t}_{ij}}\right)+\phi^{\widehat{t}_{ij}/2}\overline{\phi}{ \ }^{-\overline{z}_{ij}/2}\widetilde{P}_{-\widetilde{y}_{ij}} \right].
\end{aligned}}
\end{equation}
Analogously to the Boninger algorithm, we apply such transformations to the quaternary matrix until we obtain an expression involving only diagonal matrices. This means that among all results of the planar decomposition for diagrams and their checkerboard colorings, only those diagrams remain in which every region has crossings with only one other region — the one corresponding to the row and column removed in the construction of the quaternary Goeritz matrix. In other words, we arrive at a configuration analogous to Fig.\,\ref{fig:diag-Goeritz-ex}. Thus, the vertical and horizontal resolutions of crossings in the planar decomposition between two regions produce diagrams that are isotopic to each other up to isolated cycles. We also note that the closure of a braid contributes the factor of $D_N$ to the vertical resolution. From these considerations, we obtain the action on entries of diagonal quaternary matrices: 

\begin{equation}\label{StepHOMFLYfromGoeritz2}
\boxed{
\mathcal{M} [\mathfrak{g}_{i}] =\mathscr{D}[\mathfrak{g}_{i}]\mathcal{M}\left[\mathfrak{g}_{i}^{(III)}\right],}
\end{equation}
where we have introduced the following abbreviated notation:
\begin{equation}
\boxed{\begin{aligned}
    \mathscr{D}[\mathfrak{g}_{i}] &= \phi^{-(x_{i}+\widehat{t}_{i})/{2}}\overline{\phi}{}^{(\widetilde{y}_{i}+\overline{z}_{i})/{2}} D_N+\phi^{-\widehat{t}_{i}/{2}}\overline{\phi}{}^{(\widetilde{y}_{i}+\overline{z}_{i})/{2}} P_{-x_{i}} +  \\ 
    &+ \left(\phi^{-{x_{i}}/{2}}+D_NP_{-x_{i}}\right)\left[\left(\overline{\phi}{}^{{\widetilde{y}_{i}}/{2}}+  D_N\widetilde{P}_{\widetilde{y}_{i}}\right)\left(\overline{\phi}{ \ }^{{\overline{z}_{i}}/{2}}\widehat{P}_{-\widehat{t}_{i}}+\phi^{-{\widehat{t}_{i}}/{2}}\overline{P}_{\overline{z}_{i}}+D_N\overline{P}_{\overline{z}_{i}}\widehat{P}_{-\widehat{t}_{i}}\right)+\phi^{-{\widehat{t}_{i}}/{2}}\overline{\phi}{ \ }^{{\overline{z}_{i}}/{2}}\widetilde{P}_{\widetilde{y}_{i}} \right].
\end{aligned}}
\end{equation}
At the final step, just as in the Boninger algorithm, the remaining expression is proportional to the action of the function $\mathcal{M}$ on the empty set, which corresponds to a single cycle. This determines the normalization:

\begin{equation}\label{NormHOMFLYinAlgoritm}
\mathcal{M}[\varnothing]=D_N=\frac{A-A^{-1}}{q-q^{-1}}\,.
\end{equation}
Then, using the equivalence between $\mathcal{M}[\mathcal{G}]$ and the planar resolution of a bipartite diagram, we obtain the following expression for computing bipartite HOMFLY--PT polynomials by the Goeritz method:
\begin{equation}\label{HOMFLYNormBip}
    H^{\cal L} = A^{-2(N_\bullet- N_\bullet^v- N_\circ+N_\circ^v)} \cdot \phi^{(N_\bullet-N_\bullet^h-N_\bullet^v)/2}\bphi^{(N_\circ-N_\circ^h-N_\circ^v)/2}\left(\phi+D_N\right)^{N_\bullet^h}\left(\bphi+D_N\right)^{N_\circ^h}\mathcal{M}[\mathcal{G}]\,,
\end{equation}
where $N_\bullet$ and $N_\circ$ denote the numbers of the upper and lower locking elements in Fig.\,\ref{fig:pladeco}.  
The indices $v$ and $h$ indicate respectively the numbers of such vertical  
($N_\bullet^{v}$ and $N_\circ^{v}$) and horizontal ($N_\bullet^{h}$ and $N_\circ^{h}$) locks that separate the same regions, see subsection~\ref{sec:RetCalM} for details.

\subsubsection{Intersection of regions by several separate braids}\label{sec:ManyBraidsCalM}

 The resulting function $\mathcal{M}[\mathcal{G}]$ and the quaternary Goeritz matrix retain the features of the function $\mu[G]$ and the standard Goeritz matrix, as we will comment on in this and the following subsections.  Consider an element of the quaternary Goeritz matrix $\mathfrak{g}_{ij}$, corresponding to regions $i$ and $j$ that intersect in several separated bipartite two-strand tangles (analogous to Fig.~\ref{fig:sum-par}). The action of $\mathcal{M}[\mathfrak{g}_{ij}]$ on this element resolves all the two-strand tangles between regions $i$ and $j$ simultaneously. This fact becomes clear if we take into account that the contribution of each type of lock can be considered separately for all two-strand tangles at once. For each type of bipartite crossing, it is obvious that when resolving crossings in all $k$ separated bipartite two-strand tangles simultaneously, the coefficient in front of the vertical resolution calculates all crossings in the exponents of the functions $\phi$ and $\overline{\phi}$: $\Big(\phi^{\underset{k}\sum x_k/2},\overline{\phi}^{-\underset{k}\sum \widetilde{y}_k/2},\overline{\phi}^{\underset{k}\sum \overline{z}_k/2},\phi^{-\underset{k}\sum \widehat{t}_k/2}\Big)$
 which act as multiplicative factors for the function $\mathcal{M}\left[\mathfrak{g}_{ij}'^{(I)}\right]$, representing the vertical resolution of crossings $\mathfrak{g}_{ij}$ in the planar decomposition. The combinatorial coefficients corresponding to the horizontal resolutions are analogous to \eqref{PSumFormula} and satisfy the following relations:

\begin{equation}\label{QPPolynomials} 
\begin{split} 
P_{x_1+x_2} = \phi^{x_1/2}P_{x_2}+\phi^{x_2/2}P_{x_1}+P_{x_1}P_{x_2}D_N \,  ,  \ \ \ \  \widetilde{P}_{-(\widetilde{y}_1+\widetilde{y}_2)} = \overline{\phi}^{-\widetilde{y}_1/2}\widetilde{P}_{-\widetilde{y}_2}+\overline{\phi}^{-\widetilde{y}_2/2}\widetilde{P}_{-\widetilde{y}_1}+\widetilde{P}_{-\widetilde{y}_1}\widetilde{P}_{-\widetilde{y}_2}D_N, \\
\overline{P}_{\overline{z}_1+\overline{z}_2} = \overline{\phi}{ \ }^{\overline{z}_1/2}\overline{P}_{\overline{z}_2}+\overline{\phi}{ \ }^{\overline{z}_2/2}\overline{P}_{\overline{z}_1}+\overline{P}_{\overline{z}_1}\overline{P}_{\overline{z}_2}D_N, \ \ \ \ \widehat{P}_{-(\widehat{t}_1+\widehat{t}_2)}= \phi^{-\widehat{t}_1/2}\widehat{P}_{-\widehat{t}_2} + \phi^{-\widehat{t}_2/2}\widehat{P}_{-\widehat{t}_1} + \widehat{P}_{-\widehat{t}_1}\widehat{P}_{-\widehat{t}_2}D_N\,.
\end{split}
\end{equation}

In an iterative manner, the combinatorial coefficients $\Big(P_{\underset{k}\sum x_k},\widetilde{P}_{-\underset{k}\sum \widetilde{y}_k},\overline{P}_{\underset{k}\sum \overline{z}_k},\widehat{P}_{-\underset{k}\sum \widehat{t}_k}\Big)$ generalize to the case of many tangles. The expressions \eqref{QPPolynomials} can be obtained both from the definition of the combinatorial coefficients \eqref{PPolynomials} and from resolving the crossings via the planar decomposition. The latter approach makes it explicit that each term corresponds to its own diagram in the final resolution\footnote{Final here refers to resolving all types of bipartite crossings between the two regions under consideration simultaneously.}: the term with the factor $D_N$ corresponds to a disjoint sum of two certain links, while the remaining two terms correspond to the connected sum of the same links. However, the diagrams of these sums are equivalent from the point of view of the Goeritz method, since the corresponding quaternary Goeritz matrices coincide (both diagrams have no crossings between regions $i$ and $j$, and in all other aspects they are identical). Using these arguments, we can describe how the function $\mathcal{M}[\mathfrak{g}_{ij}']$ acts on an off-diagonal element $\mathfrak{g}_{ij}'$, obtained by computing locks in $k$ separated bipartite two-strand tangles between regions $i$ and $j$:

\begin{equation}\label{StepHOMFLYfromGoeritz3}
\begin{aligned}
&\mathcal{M} [\mathfrak{g}_{ij}'] = \phi^{\underset{k}\sum (x_k-\widehat{t}_k)/2}\overline{\phi}{}^{\underset{k}\sum (\overline{z}_k-\widetilde{y}_k)/2} \mathcal{M}\left[\mathfrak{g}_{ij}'^{(I)}\right]+\mathcal{M}\left[\mathfrak{g}_{ij}'^{(II)}\right] \Big(\phi^{-\underset{k}\sum \widehat{t}_k/2}\overline{\phi}{}^{\underset{k}\sum (\overline{z}_k-\widetilde{y}_k)/2} P_{\underset{k}\sum x_k} + \left(\phi^{\underset{k}\sum x_k/2}+D_NP_{\underset{k}\sum x_k}\right) \cdot\\
&\cdot \left[\left(\overline{\phi}{}^{-\underset{k}\sum \widetilde{y}_k/2}+  D_N\widetilde{P}_{-\underset{k}\sum \widetilde{y}_k}\right)\left(\overline{\phi}{ \ }^{-\underset{k}\sum \overline{z}_k/2}\widehat{P}_{-\underset{k}\sum \widehat{t}_k}+\phi^{-\underset{k}\sum \widehat{t}_k/2}\overline{P}_{\underset{k}\sum \overline{z}_k}+D_N\overline{P}_{\underset{k}\sum \overline{z}_k}\widehat{P}_{-\underset{k}\sum \widehat{t}_k}\right)+\phi^{-\underset{k}\sum \widehat{t}_k/2}\overline{\phi}{ \ }^{\underset{k}\sum \overline{z}_k/2}\widetilde{P}_{-\underset{k}\sum \widetilde{y}_k} \right]\Big), 
\end{aligned}
\end{equation}

\begin{equation}
    \mathfrak{g}_{ij}' = \underset{k}\sum x_k -\underset{k}\sum \widetilde{y}_k  +\underset{k}\sum \overline{z}_k  -\underset{k}\sum \widehat{t}_k \ .
\end{equation}

The action of the function $\mathcal{M}[\mathcal{G}]$ on a diagonal element is organized in an analogous way. Therefore, the distribution of lock elements between any two distinct regions is irrelevant, and the function $\mathcal{M}[\mathcal{G}]$ resolves such diagrams with the same result.

\subsubsection{Retouching a polynomial $\mathcal{M}$}\label{sec:RetCalM}

\begin{figure}[h!]
		\centering	
		\includegraphics[width =1.0\linewidth]{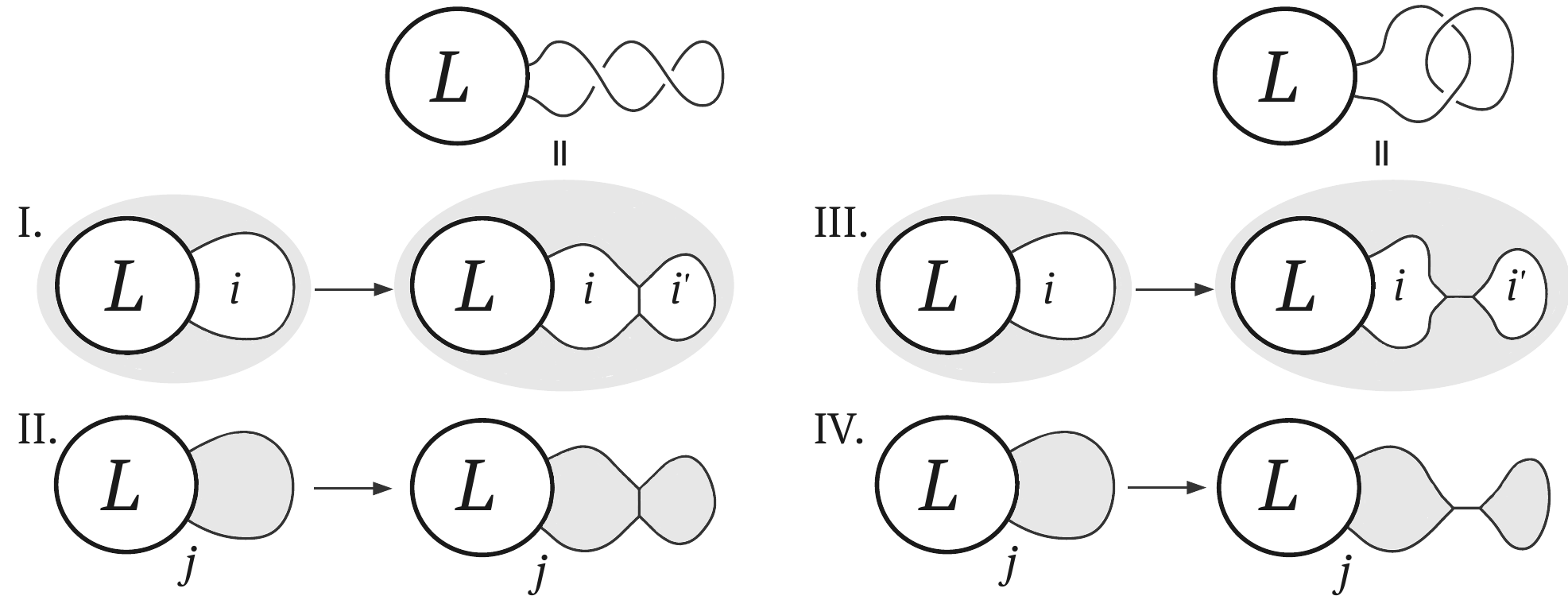}
		\caption{\footnotesize{Self-intersections of regions for positive vertical and negative horizontal bipartite crossings and the checkerboard coloring on the resulting diagrams.}}
        \label{fig:QRetoucher}
	\end{figure}

As in Section~\ref{sec:RetMu} for ordinary Goeritz matrices, problems with defining the polynomial $\mathcal{M}$ arise when considering self-intersections of a region, i.e. when a region is split by bipartite locks. Examples of such self-intersections for positive vertical and negative horizontal bipartite crossings are shown in Fig.~\ref{fig:QRetoucher} under labels II and IV, respectively. The situation is analogous for the other (opposite) types of locks. As a result, the necessity of taking self-intersections into account leads to a modification of the normalization of the polynomial $\mathcal{M}[\mathcal{G}]$ from \eqref{HOMFLYNorm}, which we refer to as retouching.

As always, there exist two types of checkerboard colorings. In the first case (I and III) a bipartite self-intersection increases the number of regions, thereby changing the size of the quaternary Goeritz matrix and a resulting expression for the polynomial $\mathcal{M}(\phi,\overline{\phi},D_N)$. This is compensated by a modification of the normalization factor \eqref{HOMFLYNorm}, due to the increased number of locks ($N_\bullet$ or $N_\circ$). The Goeritz method remains exactly as described above for such colorings.

For the second coloring (II and IV) the new region separated by the newly introduced bipartite crossing is black. Consequently, neither the number of white regions nor the number of intersections between them does not change. It follows that the quaternary Goeritz matrix also remains unchanged. At the same time, the normalization for the planar decomposition \eqref{HOMFLYNorm} changes due to the appearance of a new bipartite lock. Therefore, we need to introduce a compensating factor that retouches the polynomial $\mathcal{M}[\mathcal{G}]$ to account for bipartite crossings that are invisible to the quaternary matrix and to the Goeritz method. Moreover, horizontal and vertical bipartite crossings behave differently, which motivates introducing distinct correction factors for them.

A simpler situation arises for vertical locks (II). Resolving a self-intersection via planar decomposition produces the factors  
$\left(\phi^{-1/2}+\phi^{1/2}D_N\right) = A^{2}\phi^{-1/2}$ and $\left(\overline\phi{}^{-1/2}+\overline\phi{}^{1/2}D_N\right) = A^{-2}\overline\phi{}^{-1/2}$ for positive and negative bipartite crossings, respectively. As it is evident from Fig.\,\ref{fig:QRetoucher}, vertical self-intersecting locks do not change the link upon their removal, since they are equivalent to the first two Reidemeister moves. Therefore, the normalization factor must include the above fixed results of resolving a vertical ``lock’’ by planar decomposition, raised to the powers corresponding to the numbers of vertical positive ($N_\bullet^v$) and negative ($N_\circ^v$) bipartite self-intersections, respectively.

Horizontal bipartite locks (IV), as seen in the figure, correspond to the appearance of a new component\footnote{It is important to note that at the level of the ordinary Goeritz method (i.e. for precursor diagrams), such self-intersections could be ignored in the normalization, since one could simply remove them from the diagram (untwisting them via the first Reidemeister move). For vertical locks, which correspond to two successive first Reidemeister moves, such a simplification would still be possible. However, the appearance of a new component in the case of a horizontal self-intersection requires more careful treatment, since it changes the HOMFLY--PT polynomial.}. In terms of resolving a self-intersection via planar decomposition, this corresponds to the appearance of the factor $\left(\phi^{1/2}+\phi^{-1/2}D_N = \phi^{-1/2}\left(\phi+D_N\right)\right)$ for a negative horizontal bipartite crossing and $\left(\overline\phi{}^{1/2}+\overline\phi{}^{-1/2}D_N = \overline\phi{}^{-1/2}\left(\bphi+D_N\right)\right)$ for a positive horizontal bipartite crossing. These factors do not arise from the transformations applied to the quaternary Goeritz matrix and therefore must be included separately, raised to the powers corresponding to the numbers of positive ($N_\circ^h$) and negative ($N_\bullet^h$) horizontal locks. At the same time, this does not affect the planar normalization \eqref{HOMFLYNorm}, since horizontal self-intersecting locks actually change the link.

Nevertheless, once the retouching of the polynomial $\mathcal{M}[\mathcal{G}]$ is taken into account, the Goeritz method produces the HOMFLY--PT polynomial for any bipartite diagram. The general form of the compensating factor $\mathcal{N}_{\mathrm{com}}$ is as follows:

\begin{equation}\label{HOMFLYNormBipCom}
    \mathcal{N}_{com} = \left(A^{-2}\phi^{1/2}\right)^{-N_{\bullet}^v}\left(A^{2}\overline\phi{}^{1/2}\right)^{-N_{\circ}^v}\phi^{-N_{\bullet}^h/2}\left(\phi+D_N\right)^{N_\bullet^h}\overline\phi{}^{-N_{\circ}^h/2}\left(\bphi+D_N\right)^{N_\circ^h} .
\end{equation}

\section{Examples of calculation of HOMFLY--PT polynomials}
\label{sec:HOMFLY-ex}
\setcounter{equation}{0}

In this section, we present examples of computing HOMFLY--PT polynomials using quaternary Goeritz matrices. We consider examples of twist and two–strand torus knots as large, yet simple, knot families, each of which admits a bipartite diagram. We also provide examples of computing HOMFLY--PT polynomials for certain pretzel, rational, and Montesinos knots, for which we previously demonstrated the existence of bipartite diagrams in Section~\ref{sec:bip-diag-obtain}.

\subsection{Twist knots ${\rm Tw}_m$}\label{sec:HOMFLYTwist}

\begin{wrapfigure}[9]{r}{0.5\textwidth}
   \centering
    \vspace{-1.7cm}
    \includegraphics[width =\linewidth]{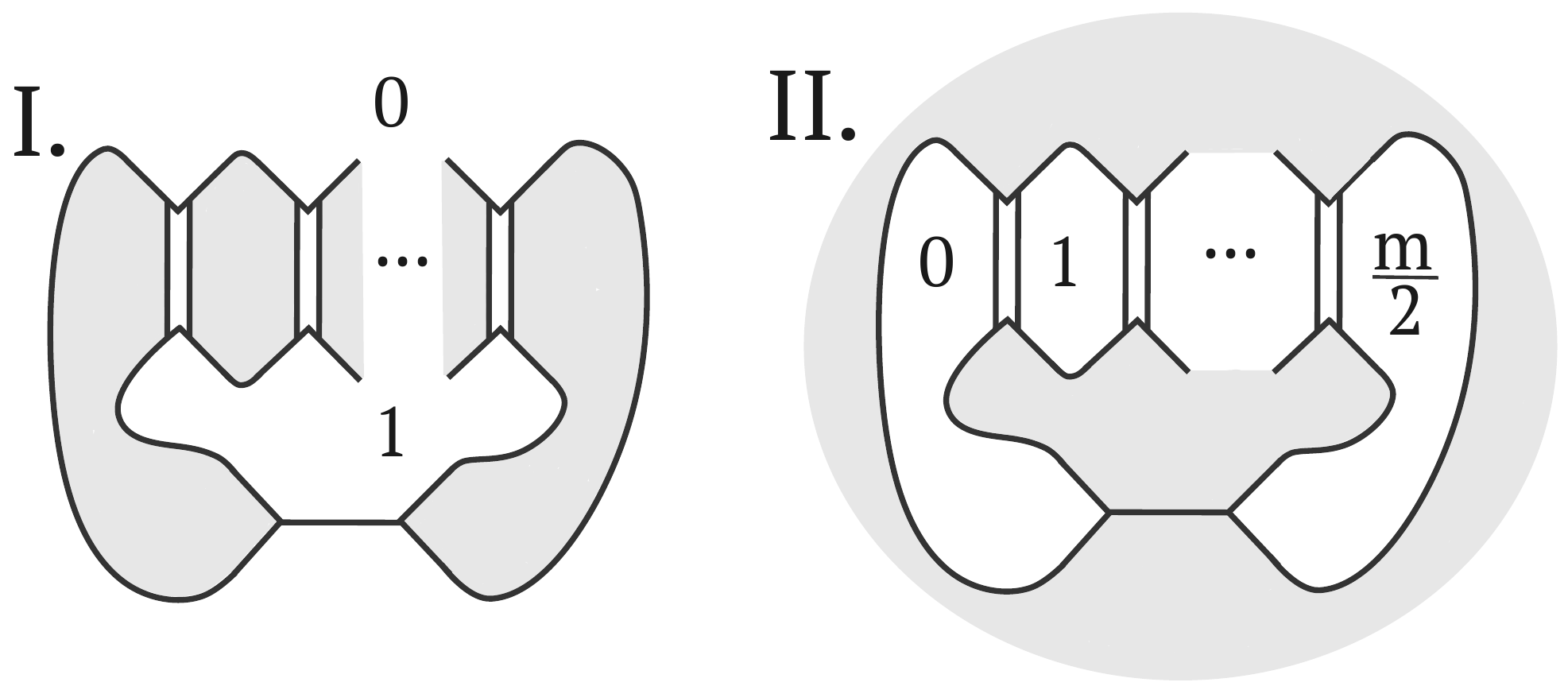}
    \caption{\footnotesize Checkerboard coloring for bipartite diagrams of twist links ${\rm Tw}_m$ for even $m$.}
        \label{fig:QuarterCheckboardForEvenTwist}
\end{wrapfigure}

One of the simplest families for computing the general form of their HOMFLY--PT polynomials by the Goeritz method is the family of twist links ${\rm Tw}_m$. In the bipartite representation, these knots and links are effectively presented as a braid of bipartite locks of only two types, which is evident for knots with even $m$ and is explicitly demonstrated in Section~\ref{sec:TwistKnots} for odd $m$. We will show the computation of the HOMFLY--PT polynomials for both types of diagrams. We begin with the case of even $m$: the quaternary Goeritz matrix for the first coloring (I) has a very simple form, since only two white regions remain:

\begin{equation}\label{QuarterGoeritzTwistEven1}
\widetilde{\mathcal G}^{\rm Tw_m}_I  =
\begin{pmatrix}
 \begin{array}{cc}  -\frac{\widetilde{m}}{2}-1 & \frac{\widetilde{m}}{2}+1 \\  \frac{\widetilde{m}}{2}+1 & -\frac{\widetilde{m}}{2}-1
 \end{array}
 \end{pmatrix}\quad \mapsto \quad \mathcal G^{\rm Tw_m}_I  = \left(-\frac{\widetilde{m}}{2}-1\right)\,.
\end{equation}
The HOMFLY--PT polynomial is obtained in a trivial manner in the general form from this matix:

\begin{equation}\label{HOMFLYFromBipTwist}
    H^{{\rm Tw}_{m}} = \left(A^{-2}\phi^{1/2}\right)\left(A^{2}\overline{\phi}^{1/2}\right)^{m/2}\mathcal{M}\left[-\frac{\widetilde{m}}{2}-1\right]\,,
\end{equation}

\begin{equation}\label{HOMFLYFromBipTwist2}
    H^{{\rm Tw}_{m}} = A^{m-2}\frac{(D_N^2-1)\phi+(D_N+\phi)(1+D_N\overline{\phi})^{m/2}}{D_N} = \frac{A^2-A^{2+m}-q^2-A^2q^2+A^mq^2+A^{4+m}q^2+A^2q^4+A^{2+m}q^4}{(A^2-1)A^2q^2}\,.
\end{equation}

There are $\frac{m}{2}$ white regions for the second coloring of a knot ${\rm Tw}_m$ (see Fig.~\ref{fig:QuarterCheckboardForEvenTwist}), which determines the size $(\tfrac{m}{2}\times\tfrac{m}{2})$ of the quaternary Goeritz matrices for the twist knots. Their structure is analogous, taking into account the different types of crossings, to the ordinary two-strand Goeritz matrices \eqref{GoeritzMatrixTorKont2}:

\begin{equation}\label{QuarterGoeritzTwistEven2}
\widetilde{\mathcal{G}}_{II}^{{\rm Tw}_{m}}  =
\begin{pmatrix}
\overline{1}+\widehat{1} & -\overline{1} & 0 & \dots  & 0 &-\widehat{1}\\
-\overline{1} & \overline{2} & -\overline{1} & \dots & & 0 \\
0 & -\overline{1} & \ddots &  & &\vdots\\

\vdots  & \vdots & &\ddots & -\overline{1} & 0 \\

0  & & & -\overline{1} & \overline{2} & -\overline{1} \\
-\widehat{1}& 0 & \dots & 0 &-\overline{1} &\overline{1}+\widehat{1}  \\

\end{pmatrix}\quad  \mapsto \quad \mathcal{G}_{II}^{{\rm Tw}_{m}}=\begin{pmatrix}
\overline{2} & -\overline{1} & 0 & \dots  & 0\\
-\overline{1} & \overline{2} & -\overline{1} & \dots  \\
0 & -\overline{1} & \ddots &  &\vdots\\

\vdots  & \vdots & & \ddots & -\overline{1} \\
0 & & \dots &-\overline{1} & \overline{1}+\widehat{1}  \\

\end{pmatrix}. 
\end{equation}

Analogous to the calculating of the Jones polynomials for knots $T[2,p]$ in \eqref{ExampleTorKnot2}, the Goeritz transformations applied to these matrices lead to a combinatorial expression for the polynomial $\mathcal{M}[\mathcal{G}]$:

\begin{equation}\label{HOMFLYFromBipTwist3}
    \mathcal{M} [\mathcal{G}_{II}^{\rm Tw_m}] = \sum_{k=0}^{m/2-2}\left\{\begin{pmatrix}
 \begin{array}{c}  m/2-1  \\  k     
 \end{array}
 \end{pmatrix}\overline{\phi}^{(m/2-1-k)/2}\overline{P}_{-{\overline{1}}}^{k}D_N^{m/2-2-k}\right\}\mathcal{M}\begin{bmatrix}\begin{pmatrix}  \overline{1} & 0\\0 & \widehat{1}     
 \end{pmatrix}
     
 \end{bmatrix} + \overline{P}_{-\overline{1}}^{m/2-1}\mathcal{M}[\overline{1}+\widehat{1}]\,.
\end{equation}
\begin{wrapfigure}[11]{r}{0.5\textwidth}
   \centering
    \vspace{-0.8cm}
    \includegraphics[width =\linewidth]{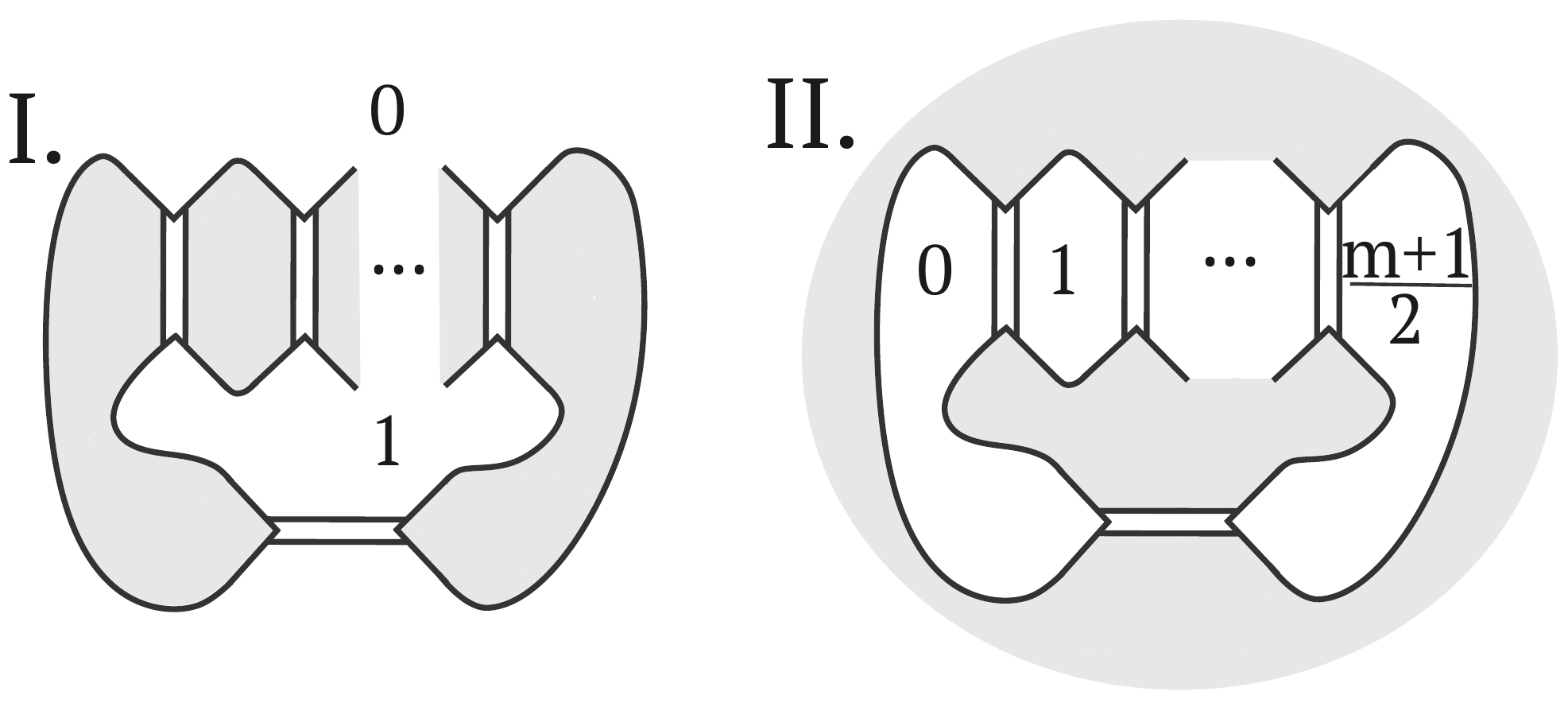}
    \caption{\footnotesize Checkerboard coloring for bipartite diagrams of twist knots ${\rm Tw}_m$ for odd $m$.}
        \label{fig:QuarterCheckboardForOddTwist}
\end{wrapfigure}
The value of the HOMFLY--PT polynomial is determined by the same normalization as in \eqref{HOMFLYFromBipTwist} (the number and types of locks are the same). The final result evidently coincides with \eqref{HOMFLYFromBipTwist2}.

For odd $m$, the actions of the Goeritz method are completely identical to the steps for twist knots with even $m$, both in terms of the number of white regions in the checkerboard coloring and the structure of the quaternary matrices, since the difference is only in a single bipartite lock. Thus, the quaternary Goeritz matrices for coloring (I) in Fig.~\ref{fig:QuarterCheckboardForOddTwist} take the following form:

\begin{equation}\label{QuarterGoeritzTwistOdd1}
\widetilde{\mathcal G}^{\rm Tw_m}_I  =
\begin{pmatrix}
 \begin{array}{cc}  -\frac{\widetilde{m}+\widetilde{1}}{2}+\overline{1} & \frac{\widetilde{m}+\widetilde{1}}{2}-\overline{1} \\  \frac{\widetilde{m}+\widetilde{1}}{2}-\overline{1} & -\frac{\widetilde{m}+\widetilde{1}}{2}+\overline{1}
 \end{array}
 \end{pmatrix}\quad \mapsto \quad \mathcal G^{\rm Tw_m}_I  = \left(-\frac{\widetilde{m}+\widetilde{1}}{2}+\overline{1}\right)\,.
\end{equation}
HOMFLY--PT polynomial takes the following value:

\begin{equation}\label{HOMFLYFromBipTwistOdd}
    H^{{\rm Tw}_{m}} = \left(A^{-2}\overline{\phi}{}^{1/2}\right)^{m/2+1}\mathcal{M}\left[-\frac{\widetilde{m}+\widetilde{1}}{2}-1\right]\,.
\end{equation}
The second coloring corresponds to matrices of size
 $\frac{m+1}{2}\times\frac{m+1}{2}\,$:

\begin{equation}\label{QuarterGoeritzTwistOdd2}
\widetilde{\mathcal{G}}_{II}^{{\rm Tw}_{m}}  =
\begin{pmatrix}
\overline{1}-\widetilde{1} & -\overline{1} & 0 & \dots  & 0 &\widetilde{1}\\
-\overline{1} & \overline{2} & -\overline{1} & \dots & & 0 \\
0 & -\overline{1} & \ddots &  & &\vdots\\

\vdots  & \vdots & &\ddots & -\overline{1} & 0 \\

0  & & & -\overline{1} & \overline{2} & -\overline{1} \\
\widetilde{1}& 0 & \dots & 0 &-\overline{1} &\overline{1}-\widetilde{1}  \\

\end{pmatrix}\quad  \mapsto \quad \mathcal{G}_{II}^{{\rm Tw}_{m}}=\begin{pmatrix}
\overline{2} & -\overline{1} & 0 & \dots  & 0\\
-\overline{1} & \overline{2} & -\overline{1} & \dots  \\
0 & -\overline{1} & \ddots &  &\vdots\\

\vdots  & \vdots & & \ddots & -\overline{1} \\
0 & & \dots &-\overline{1} & \overline{1}-\widetilde{1}  \\

\end{pmatrix}. 
\end{equation}
This matrix, like \eqref{QuarterGoeritzTwistEven2}, leads to a combinatorial expression for $\mathcal{M}[\mathcal{G}]$:

\begin{equation}\label{HOMFLYFromBipOddTwist3}
    \mathcal{M} [\mathcal{G}_{II}^{{\rm Tw}_m}] = \sum_{k=0}^{(m-3)/2}\left\{\begin{pmatrix}
 \begin{array}{c}  (m-1)/2  \\  k     
 \end{array}
 \end{pmatrix}\overline{\phi}^{((m-1)/2-k)/2}\overline{P}_{-{\overline{1}}}^{k}D_N^{(m-3)/2-k}\right\}\mathcal{M}\begin{bmatrix}\begin{array}{cc}  \overline{1} & 0\\0 & -\widetilde{1}     
 \end{array}
     
 \end{bmatrix} + \overline{P}_{-\overline{1}}^{(m-1)/2}\mathcal{M}[\overline{1}-\widetilde{1}]\,.
\end{equation}
Taking into account the normalization, this expression yields the same result as \eqref{HOMFLYFromBipTwistOdd}, namely:

\begin{equation}\label{HOMFLYFromBipTwistOdd2}
\begin{aligned}
    H^{{\rm Tw}_{m}} &= A^{m+3}\frac{(D_N^2-1)\overline{\phi}+(D_N+\overline{\phi})(1+D_N\overline{\phi})^{(m+1)/2}}{D_N} = \\ 
    &= \frac{A^2-A^{3+m}-A^2q^2-A^4q^2+A^{m+1}q^2+A^{5+m}q^2+A^2q^4-A^{3+m}q^4}{(1-A^2)q^2}
\end{aligned}
\end{equation}
Thus, the computation of the general form of the fundamental HOMFLY--PT polynomials for twist knots using the Goeritz method poses no particular difficulty.

\subsection{Two-strand links $T[2,p]$}\label{sec:HOMFLY2Torus}

\begin{figure}[h!]
		\centering	
		\includegraphics[width =\linewidth]{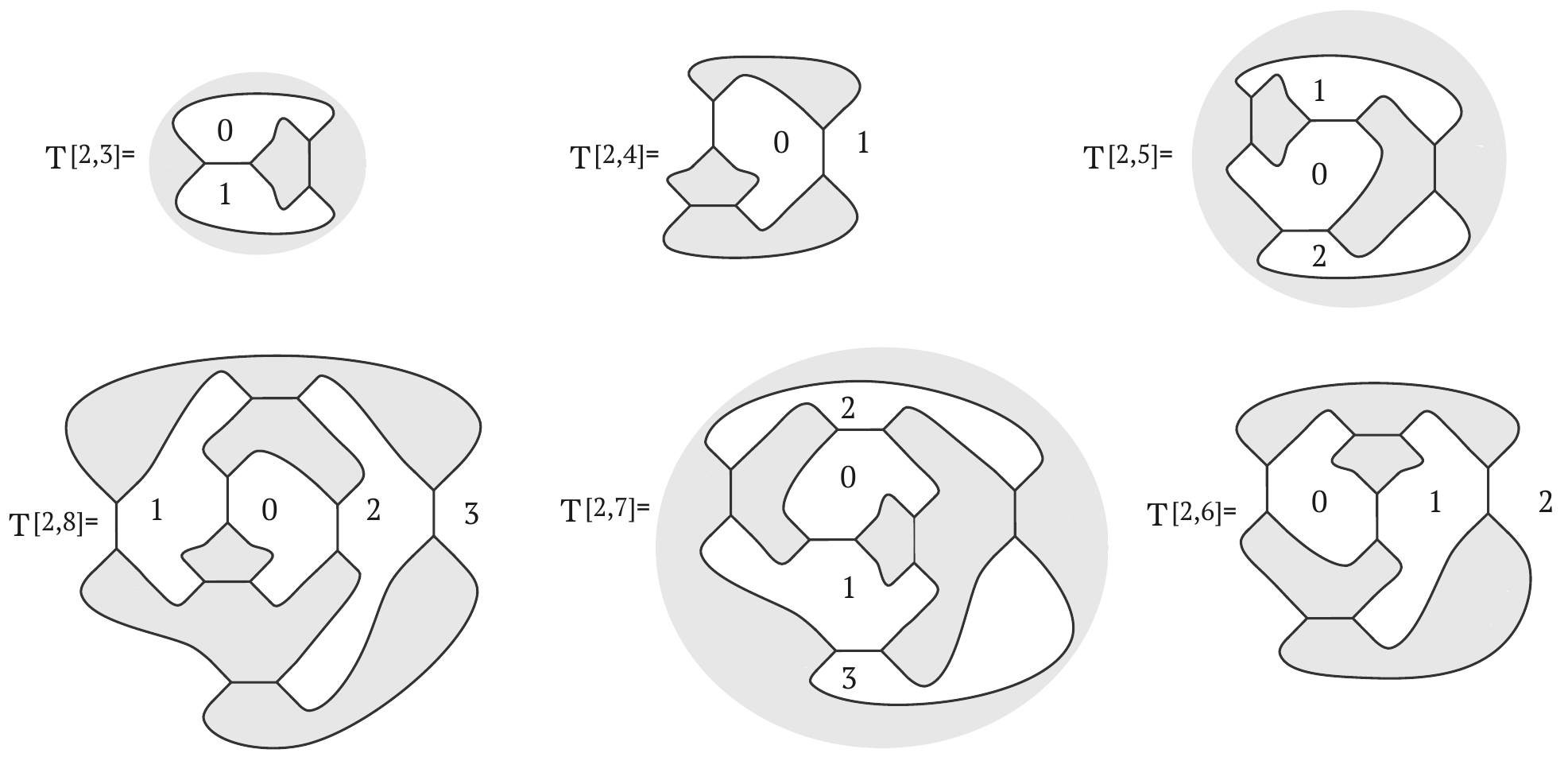}
        \caption{\footnotesize Checkerboard colorings for bipartite diagrams of the first six two-strand links, starting from the trefoil knot.}
        \label{fig:SixBipTorCheckboardDiagramm}
	\end{figure}

Two-strand links, also known as torus links $T[2,p]$, are relatively simple when represented with standard diagrams. In this regard, they are always interesting for illustrating various computations. Bipartite diagrams of two-strand knots become more complex, but they have a clear structure described in Subsection~\ref{sec:TwoBraidsKnotsAlgorithm}, namely, the new crossings are wound onto the diagrams of other two-strand links with smaller $p$ (see Fig.~\ref{fig:BipTorTangle}). This fact helps to select uniform colorings for the bipartite diagrams of $T[2,p]$, so that we can associate consistent quaternary Goeritz matrices with the entire family of two-strand links and thereby establish clear steps for computing the HOMFLY--PT polynomial of an arbitrary two-strand link. One way to define uniform colorings is to rigidly associate the white regions with the wound bipartite crossings, as explicitly illustrated in Figs.~\ref{fig:SixBipTorCheckboardDiagramm} and~\ref{fig:GeneralBipTorCheckboardDiagramm}. The smaller the number of a white region, the closer it is to the locks that appeared for smaller $p$. After several iterations the white regions are closed off from the new white regions and have no intersections with the new regions. This implies that, in general, the quaternary Goeritz matrices for two-strand diagrams will have nonzero entries only along the main diagonal and along two adjacent diagonals on each side of the main diagonal.

Next, we present the unreduced quaternary Goeritz matrices for our chosen uniform coloring and for the first six two-strand links, starting from the trefoil knot, see \eqref{BipGoeritz1Tor23}, \eqref{BipGoeritz1Tor25}, \eqref{BipGoeritz1Tor26}. The principle of associating white regions with specific locks allows the regions 0 and 1 to be closed off already in the last diagrams of Fig.~\ref{fig:SixBipTorCheckboardDiagramm}. This means that for diagrams with $p>8$, the first two rows and two columns have no nonzero entries other than those present in the matrices $\widetilde{\mathcal{G}}^{T[2,7]}$ \eqref{BipGoeritz1Tor27} and $\widetilde{\mathcal{G}}^{T[2,8]}$ \eqref{BipGoeritz1Tor28}. From these unreduced matrices, by discarding the row and column corresponding to the region 0, we can uniquely obtain the quaternary Goeritz matrices for the two-strand links. For these matrices, we provide explicit expressions for the functions $\mathcal{M}[\mathcal{G}]$ and, taking into account the normalization for the lock elements, obtain the HOMFLY--PT polynomials for the first six two-strand links. The functions $\mathcal{M}[\mathcal{G}]$ for links with larger $p$ can be expressed through $\mathcal{M}[S_3]$ and $\mathcal{M}[C_3]$, i.e. through the functions for the links $T[2,7]$ and $T[2,8]$, as we show in Appendix~\ref{sec:App2}. Therefore, we first consider the specific examples for $p \in [3,6]$ in detail and then, based on these examples, we explain how to compute the bipartite HOMFLY--PT polynomial for an arbitrary two-strand link in Appendix~\ref{sec:App2}. The expressions for the corresponding Goeritz matrices and HOMFLY–PT polynomials are given below:
\begin{equation}\label{BipGoeritz1Tor23}
\widetilde{\mathcal{G}}^{T[2,3]}  =
\begin{pmatrix}
 \begin{array}{cc}  \widehat{1}-1 & 1-\widehat{1} \\  1-\widehat{1} & \widehat{1}-1       
 \end{array}
 \end{pmatrix}\quad \mapsto \quad \mathcal{G}^{T[2,3]}  = S_1 = (\widehat{1}-1), \ \ \ \ \  \widetilde{\mathcal{G}}^{T[2,4]}  =
\begin{pmatrix}
 \begin{array}{cc}  \widehat{1}-2 & 2-\widehat{1} \\  2-\widehat{1} & \widehat{1}-2       
 \end{array}
 \end{pmatrix}\quad \mapsto \quad \mathcal{G}^{T[2,4]} =C_1 = (\widehat{1}-2), 
\end{equation}

\begin{equation}\label{BipHOMFLYPTT23}
H^{T[2,3]} = \left(A^{-2}\phi^{1/2}\right)^2\mathscr{D}[\widehat{1}-1] = \frac{(A^2 - q^2 + A^2 q^4)(A-A^{-1})}{A^4 q^2(q-q^{-1})}, 
\end{equation}

\begin{equation}\label{BipHOMFLYPTT24}
H^{T[2,4]} = \left(A^{-2}\phi^{1/2}\right)^3\mathscr{D}[\widehat{1}-2] = \frac{(A^2 - q^2 - A^2 q^2 + q^4 + A^2 q^4 - q^6 - A^2 q^6 + A^2 q^8)(A-A^{-1})}{A^5 q^4(q-q^{-1})^2}\,,
\end{equation}

\begin{equation}\label{BipGoeritz1Tor25}
\widetilde{\mathcal{G}}^{T[2,5]}  =
\begin{pmatrix}
 \begin{array}{ccc}  \widehat{1}-2 & 1-\widehat{1}  & 1  \\  1-\widehat{1} & \widehat{2}-1&-\widehat{1}  \\ 1&-\widehat{1}&\widehat{1}-1    
 \end{array}
 \end{pmatrix}\quad \mapsto \quad \mathcal{G}^{T[2,5]}  = S_2 = \begin{pmatrix}
 \begin{array}{cc}  \widehat{2}-1 & -\widehat{1} \\  -\widehat{1} & \widehat{1}-1       
 \end{array}
 \end{pmatrix}\,, 
\end{equation}

\begin{equation}\label{MPolynomialForTor25}
    \mathcal{M}\left[S_2\right] = \mathcal{U}_1[-\widehat{1}]\mathscr{D}[\widehat{1}-1]\mathscr{D}[ -1] +\mathcal{U}_2[-\widehat{1}]\mathscr{D}\left[\widehat{2}+\overline{1}-\widetilde{1}   \right] =\left( \phi^{-2}+4D_N\phi^{-1}+3+3D_N^2+4D_N\phi+\phi^2\right)D_N\,,
\end{equation}

\begin{equation}\label{HOMFLYTor25}
    H^{T[2,5]} =\left(A^{-2}\phi^{1/2}\right)^4\mathcal{M}\left[S_2\right]= \frac{(A^2 - q^2 + A^2 q^4 - q^6 + A^2 q^8)(A-A^{-1})}{A^6 q^4(q-q^{-1})}\,,
\end{equation}

\begin{equation}\label{BipGoeritz1Tor26}
\widetilde{\mathcal{G}}^{T[2,6]}  =
\begin{pmatrix}
 \begin{array}{ccc}  \widehat{1}-2 & 1-\widehat{1}  & 1  \\  1-\widehat{1} & \widehat{2}-2&1-\widehat{1}  \\ 1&1-\widehat{1}&\widehat{1}-2    
 \end{array}
 \end{pmatrix}\quad \mapsto \quad \mathcal{G}^{T[2,6]}  = C_2 = \begin{pmatrix}
 \begin{array}{cc}  \widehat{2}-2 &1 -\widehat{1} \\ 1 -\widehat{1} & \widehat{1}-2       
 \end{array}
 \end{pmatrix}\,, 
\end{equation}

\begin{equation}\label{MPolynomialForTor26}
    \mathcal{M}\left[C_2\right] = \mathcal{U}_1[1-\widehat{1}]\mathscr{D}[\widehat{1}-1]\mathscr{D}[ -1] +\mathcal{U}_2[1-\widehat{1}]\mathscr{D}\left[\widehat{2}+\overline{1}-\widetilde{1}   \right] = \small{\frac{D_N (D_N + 3 \phi + 2 D_N^2 \phi + 9 D_N \phi^2 + D_N^3 \phi^2 + 
   4 \phi^3 + 6 D_N^2 \phi^3 + 5 D_N \phi^4 + \phi^5)}{\phi^{5/2}}}\,,
\end{equation}

\begin{equation}\label{HOMFLYTor26}
    H^{T[2,6]} =\left(A^{-2}\phi^{1/2}\right)^5\mathcal{M}\left[C_2\right]= \frac{(A^2 - q^2 - A^2 q^2 + q^4 + A^2 q^4 - q^6 - A^2 q^6 + q^8 + 
 A^2 q^8 - q^{10} - A^2 q^{10} + A^2 q^{12})(A-A^{-1})}{A^7 q^5(q-q^{-1})^2}\,.
\end{equation}
The matrices for links $T[2,p]$ with larger $p$ have a common structure, which we explicitly demonstrate in Appendix~\ref{sec:App2}. There, we also show the computation of the general form of the fundamental HOMFLY--PT polynomial for an arbitrary two-strand link.

\subsection{Examples of pretzel, rational and Montesinos knots}\label{sec:HOMFLYPrRatMon}

Pretzel, rational, and Montesinos knots are large and rather complex families of knots, for which the general form of the HOMFLY--PT polynomials is unknown. Nevertheless, for each knots and even for groups of knots from these families, the computations using the Goeritz method are quite straightforward, as we demonstrate below with examples.
\paragraph{Pretzel knot $P(3,3,2)$.}

The choice of this knot to demonstrate the Goeritz method for pretzel HOMFLY--PT polynomials is motivated by the fact that in Section~\ref{sec:BipForPretzelKnots}, it served as an example of a pretzel knot admitting a bipartite diagram. The checkerboard colorings for this diagram are as follows:

\begin{figure}[h!]
		\centering	
		\includegraphics[width =0.5\linewidth]{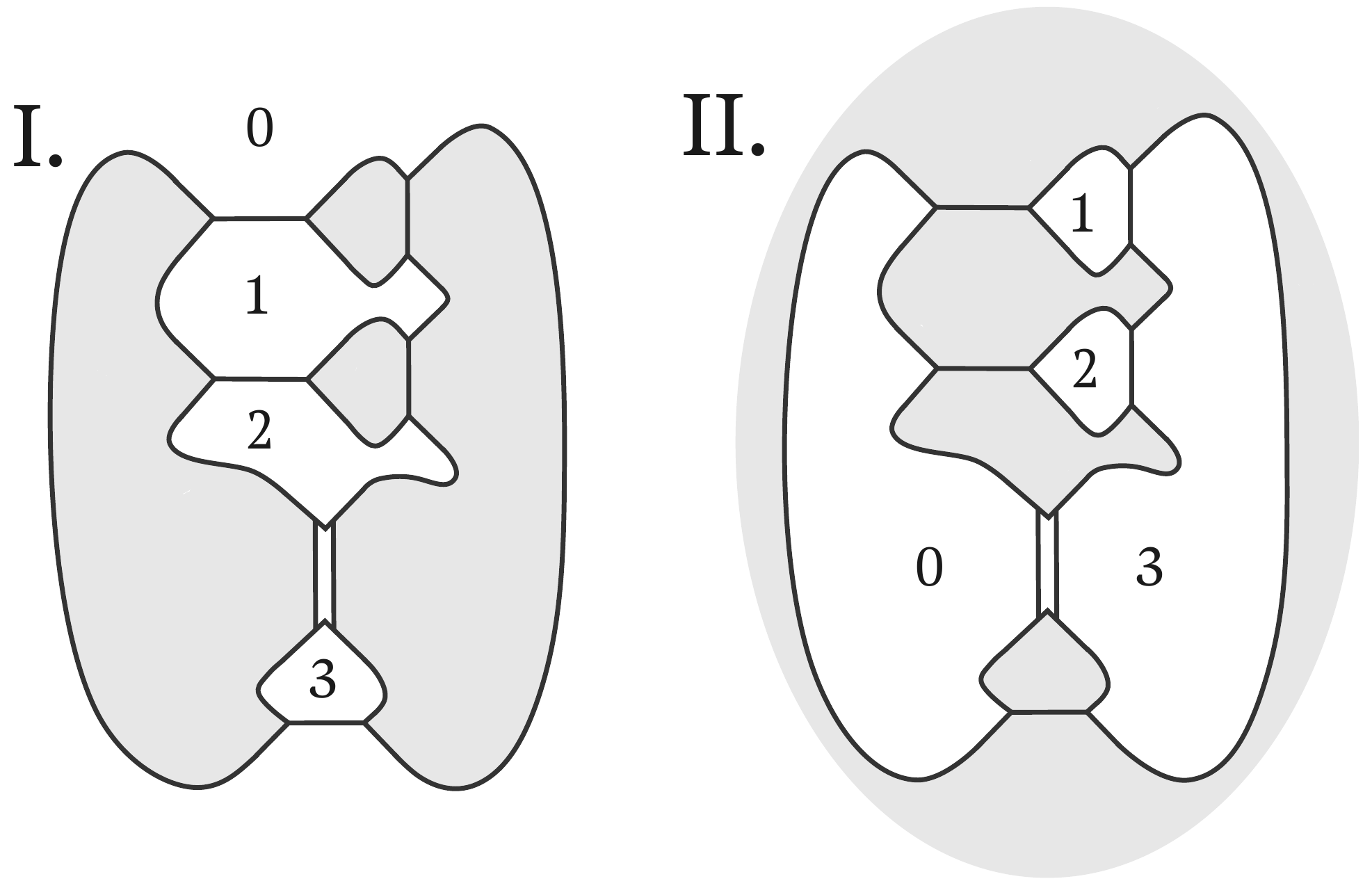}
        \caption{\footnotesize Two checkerboard colorings for the diagram of the knot $8_5$, also known as the pretzel knot $P(3,3,2)$.}
        \label{fig:P332BipTorCheckboard}
	\end{figure}

Since the number of white regions in the diagrams for both colorings is four, the quaternary Goeritz matrices for $P(3,3,2)$ have size $3 \times 3$. There are four variants for each coloring; we present one quaternary Goeritz matrix for each of them for illustration: 
    \begin{equation}\label{BipGoeritz1Pretzel332}
\widetilde{\mathcal{G}}^{P(3,3,2)}_I  =
\begin{pmatrix}
 \begin{array}{cccc} \widehat{1}-2 & 1-\widehat{1}  & 0 & 1 \\  1-\widehat{1} & \widehat{2}-2&1-\widehat{1} & 0 \\ 0&1-\widehat{1}&\widehat{1}-\widetilde{1}-1 & \widetilde{1}  \\ 1 & 0 & \widetilde{1}& -\widetilde{1} - 1    
 \end{array}
 \end{pmatrix}\quad \mapsto \quad \mathcal{G}^{P(3,3,2)}_I  =  \begin{pmatrix}
 \begin{array}{ccc}    \widehat{2}-2&1-\widehat{1} & 0 \\ 1-\widehat{1}&\widehat{1}-\widetilde{1}-1 & \widetilde{1}  \\  0 & \widetilde{1}& -\widetilde{1} - 1  
 \end{array}
 \end{pmatrix},
\end{equation}

\begin{equation}\label{BipGoeritz1Pretzel332}
\widetilde{\mathcal{G}}^{P(3,3,2)}_{II}  =
\begin{pmatrix}
 \begin{array}{cccc} \widehat{3}+\overline{1} & -\widehat{1}  & -\widehat{1} & -\widehat{1}-\overline{1} \\  -\widehat{1} & \widehat{1}-1&0 & 1 \\ -\widehat{1}&0&\widehat{1}-1 & 1  \\ -\widehat{1}-\overline{1} & 1 & 1& \widehat{1} +\overline{1}- 2    
 \end{array}
 \end{pmatrix}\quad \mapsto \quad \mathcal{G}^{P(3,3,2)}_{II}  =\begin{pmatrix}
 \begin{array}{ccc}  \widehat{3}+\overline{1} & -\widehat{1}  & -\widehat{1} \\  -\widehat{1} & \widehat{1}-1&0  \\ -\widehat{1}&0&\widehat{1}-1     
 \end{array}
 \end{pmatrix}. 
\end{equation}
For the first coloring, we present the computation of the polynomial $\mathcal{M}[\mathcal{G}]$ explicitly. For the remaining matrices, the final result is obviously the same:

\begin{equation}\label{HOMFLYPretzelP332N1}
    \mathcal{M}\left[\mathcal{G}_I^{8_5}\right] =\left( \mathcal{U}_1[1-\widehat{1}]\mathcal{U}_1[\widetilde{1}]D_N+\mathcal{U}_1[1-\widehat{1}]\mathcal{U}_2[\widetilde{1}]+\mathcal{U}_2[1-\widehat{1}]\mathcal{U}_1[\widetilde{1}]\right)\mathscr{D}\left[-1\right]\mathscr{D}\left[ \widehat{1}-1\right] D_N+\mathcal{U}_2[1-\widehat{1}]\mathcal{U}_2[\widetilde{1}]\mathscr{D}\left[\widehat{1}-2  \right]D_N \,,
\end{equation}

\begin{equation}\label{HOMFLYPretzelP332N2}
    \mathcal{M}\left[\mathcal{G}_{II}^{8_5}\right] =\mathcal{U}^2_1[-\widehat{1}]\mathscr{D}^2\left[-1\right]\mathscr{D}\left[\widehat{1}+\overline{1}\right]D_N+2\mathcal{U}_1[-\widehat{1}]\mathcal{U}_2[-\widehat{1}]\mathscr{D}[-1]\mathscr{D}\left[ \widehat{1}-1+\overline{1}\right]D_N +\mathcal{U}^2_2[-\widehat{1}]\mathscr{D}\left[\widehat{1}-2 +\overline{1} \right]D_N \,,
\end{equation}

\begin{equation}\label{HOMFLYPretzelP3320N3}
    \small{\mathcal{M}\left[\mathcal{G}_{I}^{8_{5}}\right] = \mathcal{M}\left[\mathcal{G}_{II}^{8_{5}}\right] = D_N\frac{\left(1+D_N\phi\right)\left(1+2D_N{\phi}+{\phi}^2\right)^2+\left\{D_N+\left(3+2D_N^2\right)\phi+\left(9+D_N^2\right)D_N\phi^2+2\left(2+3D_N^2\right)\phi^3+5D_N\phi^4+\phi^5\right\}\overline{\phi}}{\phi^{5/2}\overline{\phi}{}^{1/2}}} \,.
\end{equation}
Using the normalization factors and a change of variables, we arrive at the HOMFLY--PT polynomial for $P(3,3,2)$:

  \begin{equation}\label{HOMFLYFromBipForMontesinos10140}
    H^{8_5} = A^{2}\overline{\phi}{}^{1/2}\left(A^{-2}\phi^{1/2}\right)^5\mathcal{M}\left[\mathcal{G}^{8_{5}}\right]=\small{ \frac{(A-A^{-1})\left\{q^2 (1 + q^4) (1 - q^2 + q^4) + A^4 (q + q^5)^2 + A^2 (-1 + q^2 - 3 q^4 + q^6 - 3 q^8 + q^{10} - q^{12})\right\}}{A^6q^6(q-q^{-1})}}\,.
\end{equation}
    
\paragraph{Family of rational knots.}\label{sec:HOMFLYForRatioKnotsFamily}
    \begin{wrapfigure}[15]{r}{0.3\textwidth}
   \centering
    \vspace{-0.5cm}
    \includegraphics[width =\linewidth]{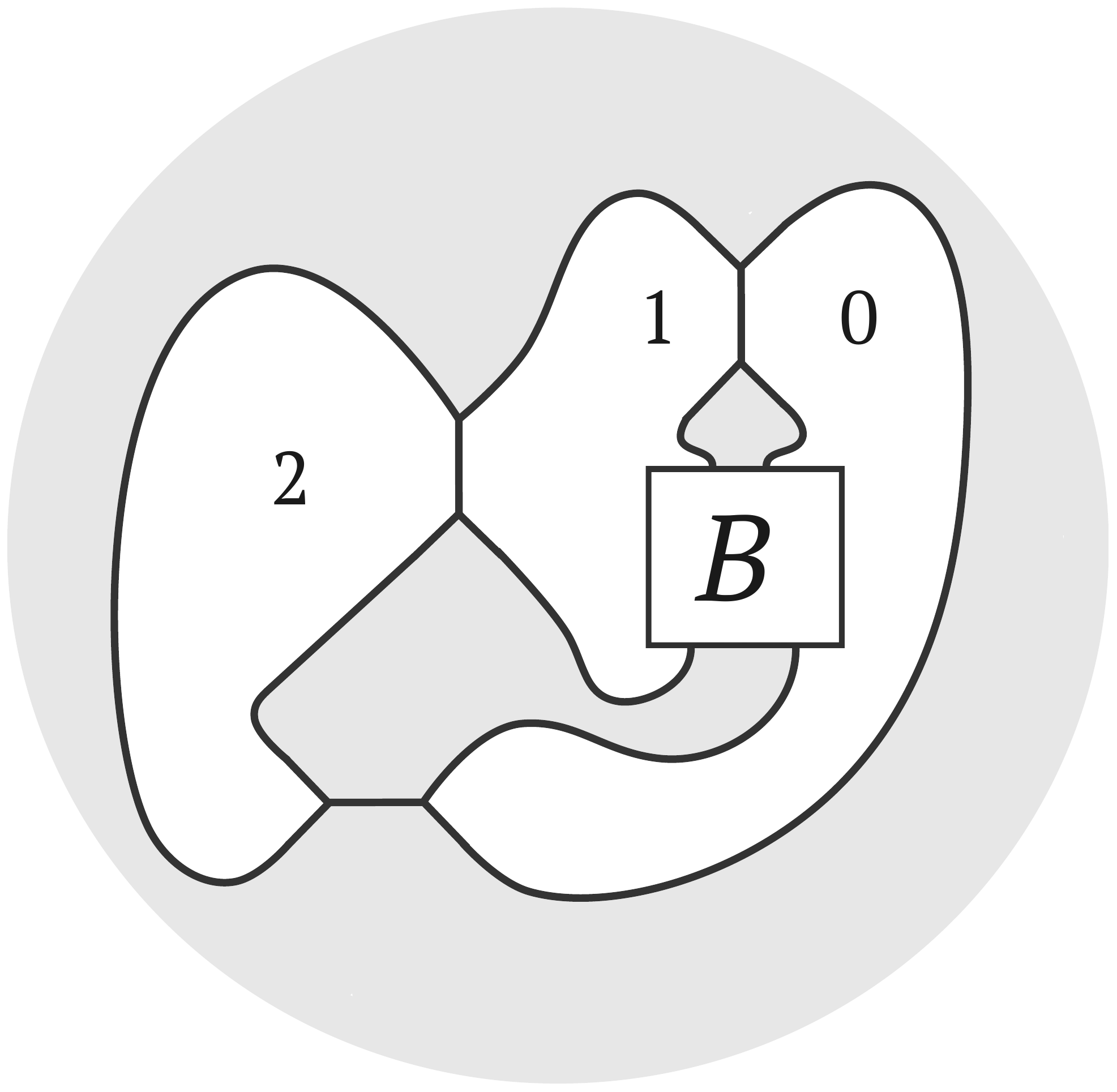}
    \caption{\footnotesize  The checkerboard coloring that allows parametrization of the family of rational knots shown in Fig.~\ref{fig:RatKnotBExample}.}
        \label{fig:BipFamilyRatio}
\end{wrapfigure}

The Goeritz method makes it possible to compute the general form of the HOMFLY--PT polynomial for parametrically defined families of knots. This becomes feasible when one parametrizes the number of bipartite crossings in the bipartite two-strand tangles between the white regions. In this case, the size of the quaternary Goeritz matrix $\mathcal{G}_{par}$ and, consequently, the number and structure of the terms in the function $\mathcal{M}[\mathcal{G}_{par}]$ (in terms of $\mathcal{U}_1, \mathcal{U}_2, \mathscr{D}$) remains fixed for any values of such parameters.

To illustrate this, let us consider a family of rational knots whose bipartite diagrams (see Fig.~\ref{fig:RatKnotBExample}) were discussed in Section~\ref{sec:ratioknots}. This class of knots is determined by the fraction~\eqref{RatTangleExample}, which in turn is specified by the value of $B$. We shall use this value to parametrize the entire family of rational knots under consideration. We choose a checkerboard coloring such that the parameter $b = |B|$ appears explicitly in the quaternary Goeritz matrices and does not affect the size of these matrices (see Fig.~\ref{fig:BipFamilyRatio}). Then, these matrices take the following form when, instead of the two-strand tangle $B$ shown in the figure, one inserts a tangle composed of positive horizontal locks $B_{+}$:

\begin{equation}\label{BipRatioPlusFamily}
\widetilde{\mathcal{G}}^{B_+}  =
\begin{pmatrix}
 \begin{array}{ccc}  \widehat{1}-1-\widetilde{b} & 1+\widetilde{b}  & -\widehat{1}  \\  1+\widetilde{b} & -2-\widetilde{b}&1 \\ -\widehat{1}&1&\widehat{1}-1   
 \end{array}
 \end{pmatrix}\quad \mapsto \quad \mathcal{G}^{B_+}  = \begin{pmatrix}
 \begin{array}{cc} -2-\widetilde{b}&1 \\  1&\widehat{1}-1  
 \end{array}
 \end{pmatrix}\, .
\end{equation}
When the parametrizing braid consists of negative horizontal locks $B_{-}$,  
the quaternary Goeritz matrix for the family under consideration changes only by the value of the parameter:

\begin{equation}\label{BipRatioMinusFamily}
\widetilde{\mathcal{G}}^{B_-}  =
\begin{pmatrix}
 \begin{array}{ccc}  \widehat{1}-1+\widehat{b} & 1-\widehat{b}  & -\widehat{1}  \\  1-\widehat{b}  & \widehat{b} -2&1 \\ -\widehat{1}&1&\widehat{1}-1   
 \end{array}
 \end{pmatrix}\quad \mapsto \quad \mathcal{G}^{B_-}  = \begin{pmatrix}
 \begin{array}{cc}  \widehat{b} -2&1 \\  1&\widehat{1}-1  
 \end{array}
 \end{pmatrix}\, .
\end{equation}
The functions $\mathcal{M}$ for the positive and negative variants transform in a simple way:
\begin{equation}\label{MPolynomialForRatioFamily} 
\left\{\begin{split} 
\mathcal{M}\left[\mathcal{G}^{B_+}\right] = \mathcal{U}_1[1]\mathscr{D}[\widehat{1}]\mathscr{D}[ -1-\widetilde{B}]D_N +\mathcal{U}_2[1]\mathscr{D}[\widehat{1}-1-\widetilde{B}]D_N , \\
\mathcal{M}\left[\mathcal{G}^{B_-}\right] = \mathcal{U}_1[1]\mathscr{D}[\widehat{1}]\mathscr{D}[ \widehat{B}-1]D_N +\mathcal{U}_2[1]\mathscr{D}[\widehat{1}+\widehat{B}-1]D_N .\ \ 
\end{split}\right.
\end{equation}
In the variables $\phi,\, \overline{\phi},\, D_N$, the functions $\mathcal{M}$ have, as expected, structurally similar expressions for both matrices. The positive variant transforms into the negative one by the simple substitution $\overline{\phi} \to \phi$, which is easy to see from the explicit form of the functions $\mathcal{M}$:

\begin{equation}\label{MPolynomialSecondForRatioFamily} 
\left\{\begin{split} 
\mathcal{M}\left[\mathcal{G}^{B_+}\right] = \overline{\phi}^{-b/2}\phi^{-3/2}\left\{\phi\left(D_N^2-1\right)\left(1+\phi\left(\phi+D_N\right)\right)+\left(\phi+D_N\right)\left(1+2D_N\phi+\phi^2\right)\left(1+D_N\overline{\phi}\right)^b\right\}, \\
\mathcal{M}\left[\mathcal{G}^{B_-}\right] = \phi^{-(b+3)/2}\left\{\phi\left(D_N^2-1\right)\left(1+\phi\left(\phi+D_N\right)\right)+\left(\phi+D_N\right)\left(1+2D_N\phi+\phi^2\right)\left(1+D_N{\phi}\right)^b\right\}.\ \ \
\end{split}\right.
\end{equation}
The expressions for the HOMFLY--PT polynomials are obtained by adding the normalization factors.  
For knots with the braid $B_{+}$, the factor is  
$\left(A^{-2}\phi^{1/2}\right)^3\left(A^{2}\overline{\phi}^{1/2}\right)^b$ 
while for knots with the braid $B_{-}$, it is  
$B_{-}$ --- $\left(A^{-2}\phi^{1/2}\right)^{3+b}$.
The final result takes the following form:

\begin{equation}\label{HOMFLYPTForRatioFamily} 
\left\{\begin{split} 
H^{B_+} = \frac{A^{2+2b}\left(A^2-q^2\right)\left(A^2q^2-1\right)(1-q^2+q^4)+\left(A^2(1+q^4)-q^2\right)\left(A^2(1-q^2+q^4)-q^2\right)}{A^7q^4(q-q^{-1})}\,, \ \ \ \ \ \ \\
H^{B_-} = A^{2b}\frac{A^{2}\left(A^2-q^2\right)\left(A^2q^2-1\right)(1-q^2+q^4)+A^{2b}\left(A^2(1+q^4)-q^2\right)\left(A^2(1-q^2+q^4)-q^2\right)}{A^7q^4(q-q^{-1})}\,. 
\end{split}\right.
\end{equation}

Similarly, one can without cumbersome calculations obtain the fundamental HOMFLY--PT polynomials for other large families of links that are parameterized by the number of locks in bipartite two-strand braids.

\paragraph{The knot $10_{140}$ in Montesinos notation.}

 \begin{wrapfigure}[11]{r}{0.3\textwidth}
   \centering
    \vspace{-0.8cm}
    \includegraphics[width =\linewidth]{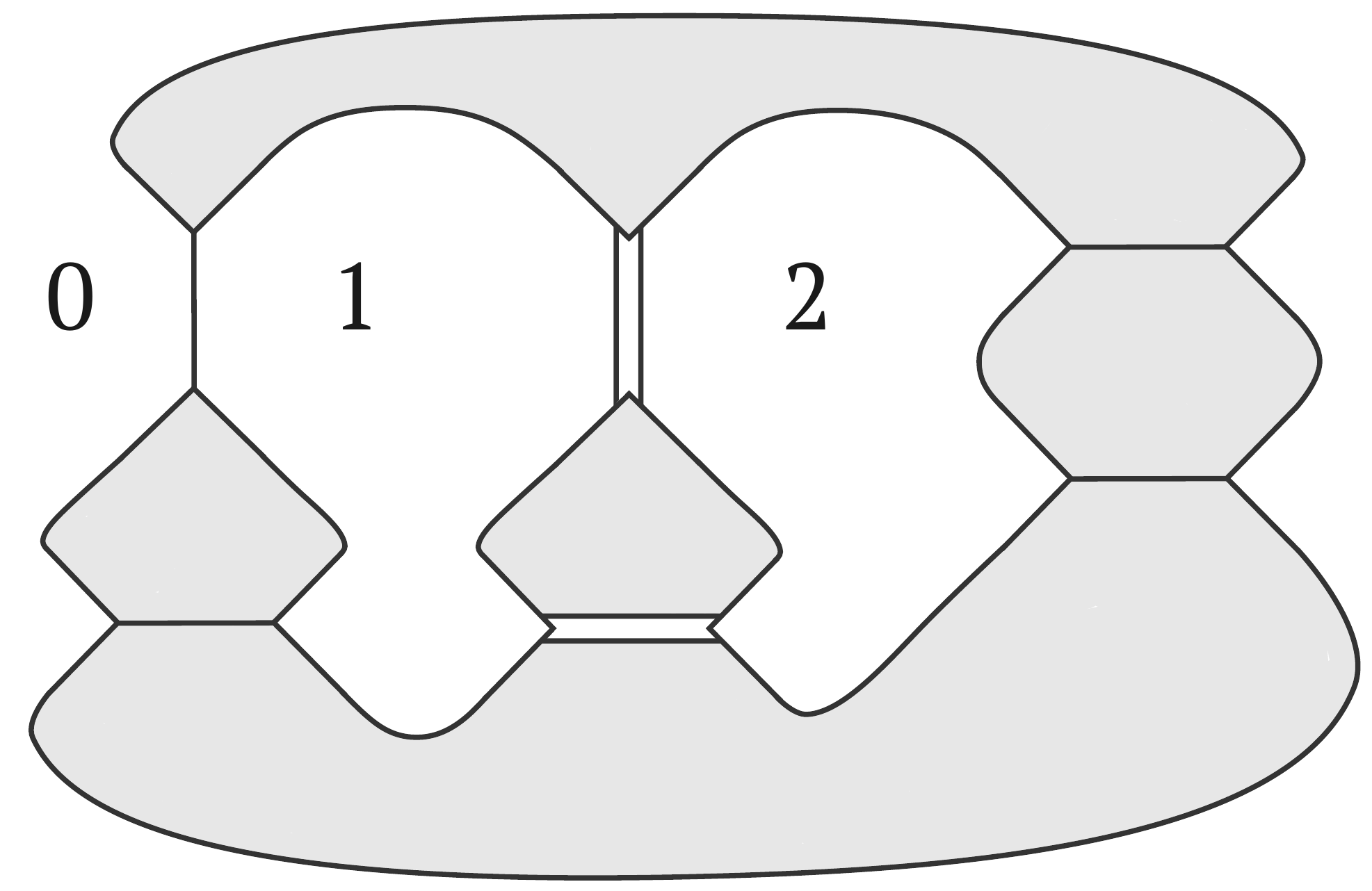}
    \caption{\footnotesize The checkerboard coloring for Montesinos knot $10_{140}=K\left(-\frac{2}{3},\,\frac{2}{3},\, -\frac{1}{4}\right)$.}
        \label{fig:Mont10140BipCheckboard}
\end{wrapfigure}

As an example of a Montesinos knot, for which explicit expressions in pretzel or rational notation are not known, we consider the knot $10_{140}$. This knot is of interest not only because it admits a bipartite diagram (see Fig.~\ref{fig:Mont10140BipCheckboard}), but also because its bipartite diagram contains all four types of locks, so that the corresponding quaternary Goeritz matrix will include all four parameters. The principles for constructing quaternary matrices have been described in detail above using various examples, so here, we limit ourselves to demonstrating the method for a single coloring, for which we provide the unreduced matrix and one of the possible quaternary Goeritz matrices:

\begin{equation}\label{BipGoeritzForMontesinos10140}
\widetilde{\mathcal{G}}^{10_{140}}  =
\begin{pmatrix}
 \begin{array}{ccc}  \widehat{3}-1 & 1-\widehat{1}  & -\widehat{2}  \\  1-\widehat{1} & \widehat{1}+\overline{1}-\widetilde{1}-1&\widetilde{1} -\overline{1} \\ -\widehat{2}&\widetilde{1}-\overline{1}&\overline{1}+\widehat{2}-\widetilde{1}    
 \end{array}
 \end{pmatrix}\quad \mapsto \quad \mathcal{G}^{10_{140}}  = \begin{pmatrix}
 \begin{array}{cc}  \widehat{3}-1 & 1-\widehat{1}   \\  1-\widehat{1} & \widehat{1}+\overline{1}-\widetilde{1}-1    
 \end{array}
 \end{pmatrix}\, .
\end{equation}
The function $\mathcal{M}$ for a given knot is expressed as follows:

\begin{equation}\label{HOMFLYMontesinos10140N1}
    \mathcal{M}\left[\mathcal{G}^{10_{140}}\right] = \mathcal{U}_1[1-\widehat{1}]\mathscr{D}[\widehat{2}]\mathscr{D}[ \overline{1}-\widetilde{1}]D_N +\mathcal{U}_2[1-\widehat{1}]\mathscr{D}\left[\widehat{2}+\overline{1}-\widetilde{1}   \right]D_N \,,
\end{equation}

\begin{equation}\label{HOMFLYMontesinos10140N2}
    \mathcal{M}\left[\mathcal{G}^{10_{140}}\right] = \small{\frac{\phi\left(D_N+2\phi+D_N\phi^2\right)\left(1+2D_N\overline{\phi}+\overline{\phi}^2\right)+\left(1+\phi\left(D_N+\phi\right)\right)\left\{D_N\overline{\phi}+\overline{\phi}{}^2+\left(1+D_N\overline{\phi}\right)\left(1+\phi+\left(2+D_N{\phi}\right)\left(D_N+\overline{\phi}\right)\right)\right\}}{\phi^2\overline{\phi}}}D_N \,.
\end{equation}
Then, applying the normalization with respect to the number of lock elements, we obtain the HOMFLY--PT polynomial for $10_{140}$:

  \begin{equation}\label{HOMFLYFromBipForMontesinos10140}
    H^{10_{140}} = \left(A^{-2}\phi^{1/2}\right)^4\left(A^{2}\overline{\phi}{}^{1/2}\right)^{2}\mathcal{M}\left[\mathcal{G}^{10_{140}}\right]= \frac{(A^6q^4-q^2(1+q^4)+A^2(1+q^4)^2-A^4q^2(1+q^4))(A-A^{-1})}{A^6q^4(q-q^{-1})}\,.
\end{equation}

\section{Conclusion}

In this work, we have explicitly shown that the method allowing one to compute the Jones polynomial from the Goeritz matrix of a checkerboard coloring of a link diagram ~\cite{boninger2023jones} can be understood as an alternative perspective on resolving a link diagram via the Kauffman bracket. Using the intuition developed from this method, we constructed a new approach for computing bipartite HOMFLY--PT polynomials. This approach, in turn, is closely related to the planar decomposition for bipartite HOMFLY--PT polynomials~\cite{anokhina2024planar,anokhina2025bipartite,anokhina2025planar}. To implement our method, we have introduced new rules for checkerboard coloring of bipartite diagrams and a new object --- the quaternary Goeritz matrix, which encodes the numbers of lock elements in a bipartite diagram separating two adjacent regions. We have also provided an algorithm that allows us to compare the polynomial in the variables of the planar decomposition for the quaternary Goeritz matrix, which gives the HOMFLY--PT polynomial, with the correct normalization. This procedure has been illustrated on numerous examples: twist knots, two-strand knots, pretzel, rational, and Montesinos knots and links. The new method turns out to be particularly convenient for computing the general form of HOMFLY--PT polynomials for classes of links distinguished by the number of bipartite locks in individual bipartite two-strand braids.

\setcounter{equation}{0}
\section*{Acknowledgements}\label{Acknowledgements}

We would like to thank A. Sleptsov, And. Morozov, A. Popolitov, and D. Galakhov for helpful discussions. This work is partially funded by the state assignment of the NRC Kurchatov Institute and by the grants of the Foundation for the Advancement of Theoretical Physics and Mathematics “BASIS”. The research of E.L. and D.K. is partly supported by the Ministry of Science and Higher Education of the Russian Federation (agreement no. 075--15--2025--343 of 29.04.2025).

\printbibliography

@article{boninger2023jones,
  title={The Jones polynomial from a Goeritz matrix},
  author={Boninger, J.},
  journal={Bulletin of the London Mathematical Society},
  volume={55},
  number={2},
  pages={732--755},
  year={2023},
  publisher={Wiley Online Library},
    eprint = {2110.03082},
    archivePrefix = {arXiv},
    primaryClass = {math.GT}
}

@article{goeritz1933knoten,
  title={Knoten und quadratische Formen},
  author={Goeritz, L.},
  journal={Mathematische Zeitschrift},
  volume={36},
  number={1},
  pages={647--654},
  year={1933},
  publisher={Springer}
}

@article{anokhina2024planar,
  title={Planar decomposition of the HOMFLY polynomial for bipartite knots and links},
  author={Anokhina, A. and Lanina, E. and Morozov, A.},
  journal={The European Physical Journal C},
  volume={84},
  number={9},
  pages={990},
  year={2024},
  publisher={Springer},
    eprint = {2407.08724},
    archivePrefix = {arXiv},
    primaryClass = {hep-th}
}

@article{anokhina2025planar,
  title={Planar decomposition of bipartite HOMFLY polynomials in symmetric representations},
  author={Anokhina, A. and Lanina, E. and Morozov, A.},
  journal={Physical Review D},
  volume={111},
  number={4},
  pages={046018},
  year={2025},
  publisher={APS},
    eprint = {2410.18525},
    archivePrefix = {arXiv},
    primaryClass = {hep-th}
}

@article{anokhina2025bipartite,
  title={Bipartite expansion beyond biparticity},
   volume={1014},
   journal={Nuclear Physics B},
   publisher={Elsevier BV},
   author={Anokhina, A. and Lanina, E. and Morozov, A.},
   year={2025},
    pages={116881},
    eprint = {2501.15467},
    archivePrefix = {arXiv},
    primaryClass = {hep-th}
}

@article{gordon1978signature,
  title={On the signature of a link},
  author={Gordon, C.Mc. and Litherland, R.A.},
  journal={Inventiones mathematicae},
  volume={47},
  pages={53--69},
  year={1978},
  publisher={Springer}
}

@article{jones2005jones,
  title={The Jones polynomial},
  author={Jones, V.F.R.},
  journal={Discrete Math},
  volume={294},
  pages={275--277},
  year={2005},
  publisher={Citeseer}
}

@article{Witten1988hf,
title={Quantum field theory and the Jones polynomial},
    author = "Witten, E.",
    editor = "Mitra, Asoke N.",
    reportNumber = "IASSNS-HEP-88-33",
    journal = "Commun. Math. Phys.",
    volume = "121",
    pages = "351--399",
    year = "1989"
}

@article{freyd1985new,
  title={A new polynomial invariant of knots and links},
  author={Freyd, P. and Yetter, D. and Hoste, J. and Lickorish, W.B.R. and Millett, K. and Ocneanu, A.},
  journal={Bulletin (new series) of the American mathematical society},
  volume={12},
  number={2},
  pages={239--246},
  year={1985},
  publisher={American Mathematical Society}
}

@article{przytycki1987kobe,
  title={Kobe J. Math.},
  author={Przytycki, J.H. and Traczyk, K.P.},
  journal={Invariants of links of Conway type},
  volume={4},
  pages={115--139},
  year={1987},
	eprint = {1610.06679},
    archivePrefix = {arXiv},
    primaryClass = {math.GT}
}

@article{mironov2012character,
  title={Character expansion for HOMFLY polynomials. II. Fundamental representation. Up to five strands in braid},
  author={Mironov, A. and Morozov, A. and Morozov, And.},
  journal={Journal of High Energy Physics},
  volume={2012},
  number={3},
  pages={1--34},
  year={2012},
  publisher={Springer},
    eprint = {1112.2654},
    archivePrefix = {arXiv},
    primaryClass = {math.QA}
}

@article{goldman1997rational,
  title={Rational tangles},
  author={Goldman, J.R. and Kauffman, L.H.},
  journal={Advances in Applied Mathematics},
  volume={18},
  number={3},
  pages={300--332},
  year={1997},
  publisher={Elsevier}
}

@article{kauffman2004classification,
  title={On the classification of rational tangles},
  author={Kauffman, L.H. and Lambropoulou, S.},
  journal={Advances in Applied Mathematics},
  volume={33},
  number={2},
  pages={199--237},
  year={2004},
  publisher={Elsevier},
    eprint = {math/0311499},
    archivePrefix = {arXiv},
    primaryClass = {math.GT}
}

@article{duzhin2010formula,
  title={A formula for the HOMFLY polynomial of rational links},
  author={Duzhin, Sergei and Shkolnikov, Mikhail},
  journal={Arnold Mathematical Journal},
  volume={1},
  number={4},
  pages={345--359},
  year={2015},
  publisher={Springer},
    eprint = {1009.1800},
    archivePrefix = {arXiv},
    primaryClass = {math.GT}
}

@article{lewark2016new,
  title={New quantum obstructions to sliceness},
  author={Lewark, L. and Lobb, A.},
  journal={Proceedings of the London Mathematical Society},
  volume={112},
  number={1},
  pages={81--114},
  year={2016},
  publisher={Oxford University Press},
    eprint = {1501.07138},
    archivePrefix = {arXiv},
    primaryClass = {math.GT}
}

@book{rolfsen2003knots,
  title={Knots and links},
  author={Rolfsen, D.},
  number={346},
  year={2003},
  publisher={American Mathematical Soc.}
}

@article{dunfield2001table,
  title={A table of boundary slopes of Montesinos knots},
  author={Dunfield, N.M.},
  journal={Topology},
  volume={40},
  number={2},
  pages={309--315},
  year={2001},
  publisher={Elsevier},
    eprint = {math/9901120},
    archivePrefix = {arXiv},
    primaryClass = {math.GT}
}

@article{Kauff,
  title={State models and the Jones polynomial},
  author={Kauffman, L.H.},
  journal={Topology},
  volume={26},
  number={3},
  pages={395--407},
  year={1987},
  publisher={Elsevier}
}

@incollection{jones1987hecke,
  title={Hecke algebra representations of braid groups and link polynomials},
  author={Jones, V.F.R.},
  booktitle={New Developments In The Theory Of Knots},
  pages={20--73},
  year={1987},
  publisher={World Scientific}
}

@article{jones1985polynomial,
  title={A polynomial invariant for knots via von Neumann algebras},
  author={Jones, V.F.R.},
  year={1985}
}

@article{reshetikhin1990ribbon,
  title={Ribbon graphs and their invaraints derived from quantum groups},
  author={Reshetikhin, N. and Turaev, V.},
  journal={Communications in Mathematical Physics},
  volume={127},
  number={1},
  pages={1--26},
  year={1990},
  publisher={Springer}
}

@article{turaev1990yang,
  title={The Yang--Baxter equation and invariants of links},
  author={Turaev, V.},
  journal={New Developments in the Theory of Knots},
  volume={11},
  pages={175},
  year={1990},
  publisher={World Scientific}
}

@article{reshetikhin1991invariants,
  title={Invariants of 3-manifolds via link polynomials and quantum groups},
  author={Reshetikhin, N. and Turaev, V.},
  journal={Inventiones mathematicae},
  volume={103},
  number={1},
  pages={547--597},
  year={1991}
}

@article{BipKnots,
  title={Bipartite knots},
  author={Duzhin, S. and Shkolnikov, M.},
  journal={Fundamenta Mathematicae},
  volume={225},
  number={1},
  pages={95--102},
  year={2014},
  publisher={Institute of Mathematics Polish Academy of Sciences},
    eprint = {1105.1264},
    archivePrefix = {arXiv},
    primaryClass = {math.GT}
}

@misc{anokhina2025khrbip,
      title={Khovanov-Rozansky cycle calculus for bipartite links}, 
      author={A. Anokhina and E. Lanina and A. Morozov},
      year={2025},
      eprint={2506.08721},
      archivePrefix={arXiv},
      primaryClass={hep-th} 
}


\appendix

\setcounter{equation}{0}
\section{Comment on the number of parameters of the generalised  \\ Goeritz Matrix}\label{sec:App1}

At first glance, one might assume that the number of distinct parameters in the generalized Goeritz matrix should be less than four, since some bipartite locks cancel each other and the planar decomposition for the HOMFLY--PT polynomial introduces only two new variables $\phi$ and $\overline{\phi}$ instead of a single variable $q$ for the Jones polynomial. The considerations in this chapter will explain why it is nevertheless necessary to choose the number of parameters in the generalized matrix to be four.

\begin{figure}[h!]
		\centering	
		\includegraphics[width =1.0\linewidth]{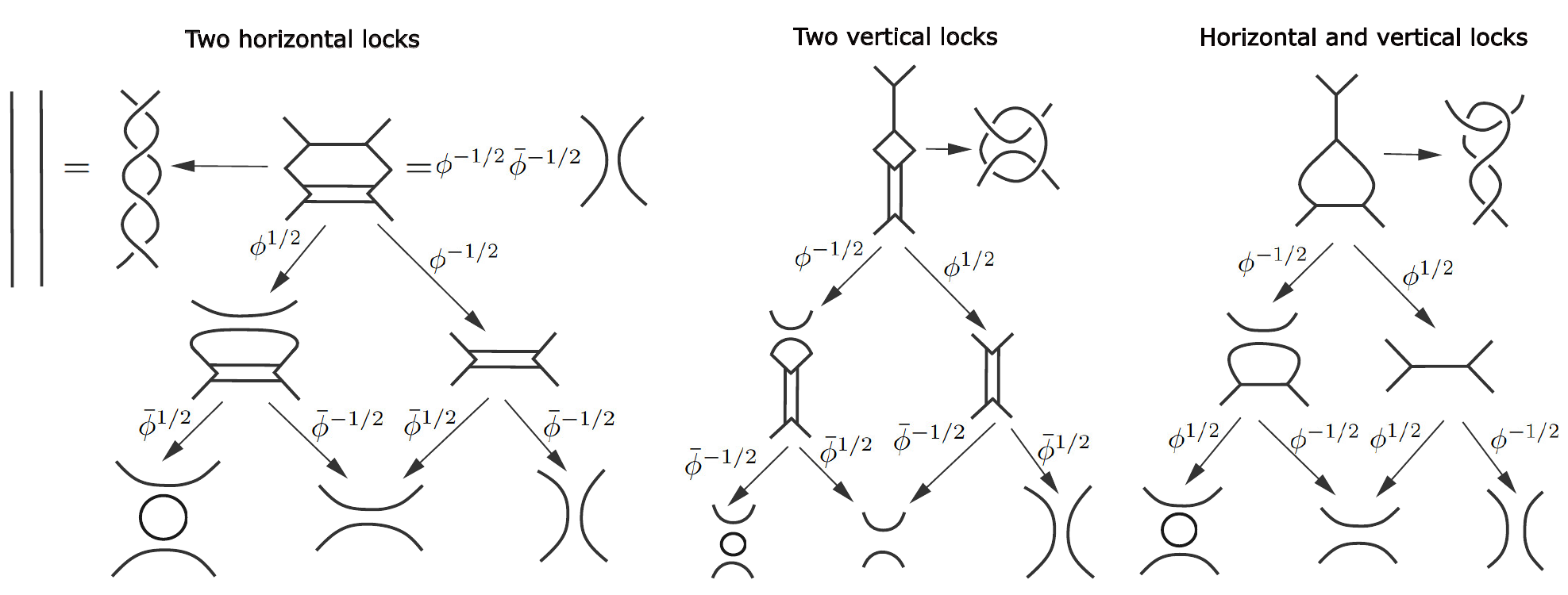}
		\caption{\footnotesize{Resolution of consecutively arranged mutually inverse horizontal and vertical locks.}}
        \label{fig:AnihilationBipLock}
	\end{figure}

A logical justification for a smaller number of parameters would be the cancellation of positive bipartite locks by negative ones\footnote{Positive horizontal and vertical bipartite crossings, like negative ones, do not cancel each other due to a clear reason. The reduction from bipartite diagrams to the precursor diagrams results in identical crossings, which are obviously not inverse of each other.}.  
Such cancellation occurs obviously for precursor crossings. However, this does not always hold for bipartite crossings and as can be seen from Fig.~\ref{fig:AnihilationBipLock}, it is true only for horizontal locks.

The cancellation of horizontal bipartite crossings can be seen explicitly at the level of ordinary diagrams: positive and negative bipartite locks joined consecutively can be untangled by two second Reidemeister moves (see Fig.~\ref{fig:AnihilationBipLock}). This is also visible from the coefficient of the planar decomposition: in our framing, the coefficient of the vertical resolution $\left(\phi^{-1/2}\overline{\phi}{}^{-1/2}\right)$
becomes unity once the normalization factor \eqref{HOMFLYNorm} is taken into account, while the coefficient of the horizontal resolution vanishes $\left(\phi^{1/2}\overline{\phi}{}^{1/2}D_N + \phi^{1/2}\overline{\phi}{}^{-1/2} + \phi^{-1/2}\overline{\phi}{}^{1/2}  = 0\right)$. 

Naively, one might expect that mutually inverse vertical bipartite crossings would also cancel. However, an explicit representation of these locks in terms of their precursor ordinary crossings shows that an additional component remains (see Fig.~\ref{fig:precursor-bip}). The same phenomenon is reflected in the planar decomposition: the coefficient of the horizontal resolution does not vanish $\left(\phi^{-1/2}\overline{\phi}{}^{-1/2}D_N + \phi^{1/2}\overline{\phi}{}^{-1/2} + \phi^{-1/2}\overline{\phi}{}^{1/2}  \ne 0\right)$.
Furthermore, mutually inverse horizontal and vertical locks cannot be governed by a single parameter, since they do not cancel each other. Equivalently, the coefficient of the resulting horizontal resolution remains nonzero $\left(D_N + \phi + {\phi}^{-1}  \ne 0\right)$\footnote{In Fig.~\ref{fig:AnihilationBipLock} and in the text, we discuss the case of upper locks from Fig.~\ref{fig:pladeco}, but the same is obviously true for the lower locks.}.

Thus, we have demonstrated that fewer than three independent parameters are insufficient. However, allowing exactly three distinct parameters leads to a further obstruction related to normalization. Indeed, along different Kauffman bracket resolution paths, the Goeritz matrix detects different numbers of crossings: for some resolution diagrams, the mutually inverse horizontal locks cancel, while for others they do not.  
Consequently, the corresponding normalization factors become path-dependent. Then, we either stop getting the correct answers for the HOMFLY-PT polynomials or we are forced to track diagrams when resolving by Kauffman brackets which would deprive the Goeritz method of its main advantage. In the next subsection, we present an explicit example illustrating how such an incorrect normalization arises.

\paragraph{Counterexample to a three-parameter generalization of the Goeritz matrix.}

Let us assume, that three parameters are sufficient for a generalized Goeritz matrix. In this setting, we let $x_{ij}$ count positive vertical bipartite crossings, $\widetilde{y}_{ij}$ count all horizontal bipartite crossings and $\overline{z}_{ij}$ count negative vertical crossings. The rules for filling the generalized Goeritz matrix according to a checkerboard coloring then take the following form:

\begin{figure}[h!]
		\centering	
		\includegraphics[width =0.6\linewidth]{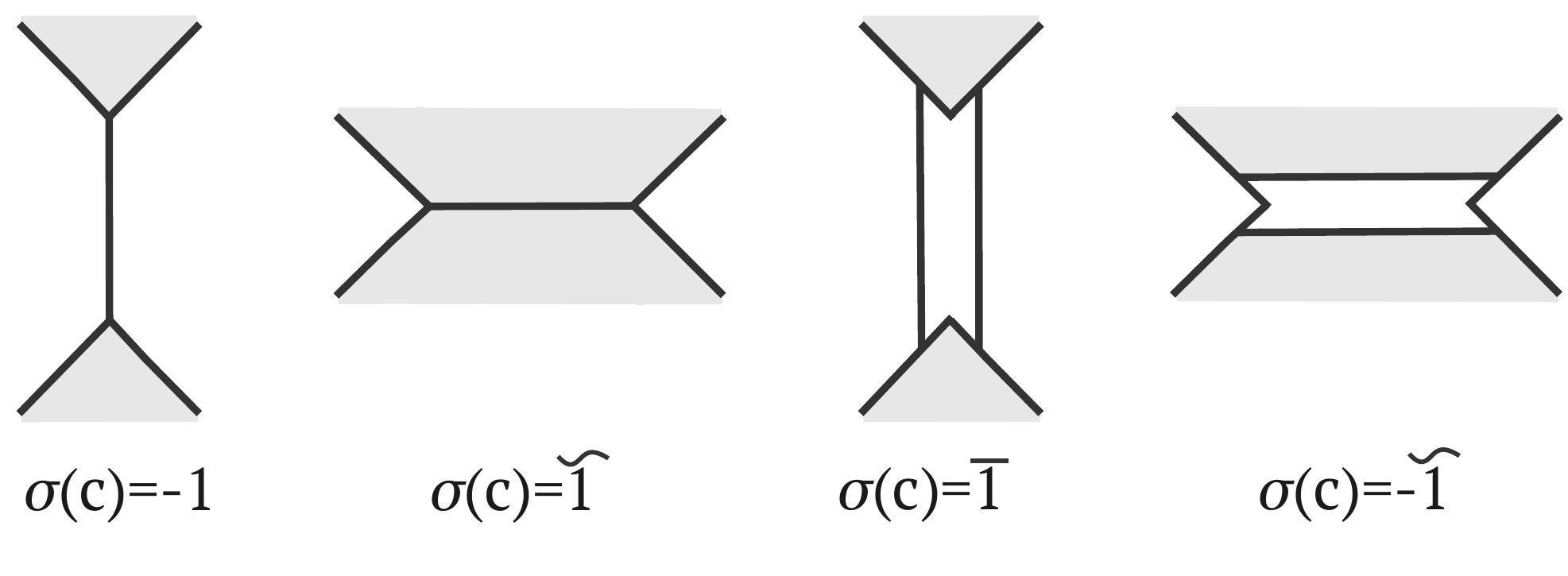}
        \caption{\footnotesize Checkerboard sign convention for bipartite diagrams assuming three independent parameters in the generalized Goeritz matrix.}
        \label{fig:CheckerboardHOMFLYPT3}
	\end{figure}

\noindent The combinatorial coefficients differ slightly from those appearing in the four-parameter case. In particular, the coefficients associated with the resolution of horizontal bipartite locks merge into a single term:
\begin{equation}\label{FalsePPolynomials} 
\begin{split} 
P_{x}(\phi,\overline{\phi}, D_N) = \sum^{|x|}_{k=1} \begin{pmatrix}
|x| \\
k 
\end{pmatrix} D_N^{k-1} \phi^{|x|/2-k} \,  , \ \ \ \  \overline{P}_{\overline{z}}(\phi,\overline{\phi}, D_N) = \sum^{|\overline{z}|}_{k=1} \begin{pmatrix}
|\overline{z}| \\
k 
\end{pmatrix} D_N^{k-1} \overline{\phi}^{|\overline{z}|/2-k}, \\
\widetilde{P}_{\widetilde{y}}(\phi,\overline{\phi}, D_N) = \sum^{|\widetilde{y}|}_{k=1} \begin{pmatrix}
|\widetilde{y}| \\
k 
\end{pmatrix} D_N^{k-1} \Big(\delta_1^{\sgn(z)}{\phi}^{k-|\widetilde{y}|/2} + \delta_{-1}^{\sgn(z)}\overline{\phi}^{k-|\widetilde{y}|/2}\Big) \ \ \ \ .
\end{split}
\end{equation}
For simplicity, let us introduce the function:

\begin{equation}
    L_{\widetilde{y}}(\phi, \overline{\phi}) = \delta_1^{\sgn(z)}{\phi^{-|\widetilde{y}|/2}} + \delta_{-1}^{\sgn(z)}\overline{\phi}^{-|\widetilde{y}|/2} + \delta^{\sgn(z)}_0.
\end{equation}
Then, the action of the Goeritz algorithm on a non-diagonal element $(i \neq j)$ of the three-parameter matrix takes the following form:

\begin{equation}\label{StepFalseHOMFLYfromGoeritz1}
\begin{aligned}
&\mathcal{M} [\mathfrak{g}_{ij}] = \phi^{x_{ij}/2}\overline{\phi}{}^{-\overline{z}_{ij}/2}L_{-\widetilde{y}_{ij}}\mathcal{M}[\mathfrak{g}_{ij}^{(I)}]+\\
&+\left\{\left( \phi^{x_{ij}/2}\overline{P}_{-\overline{z}_{ij}}+(\overline{\phi}{}^{-\overline{z}_{ij}/2}+D_N\overline{P}_{-\overline{z}_{ij}}){P}_{x_{ij}}\right)\left( L_{-\widehat{y}_{ij}}+D_N \widetilde{P}_{-\widehat{y}_{ij}}\right) + \phi^{x_{ij}/2}\overline{\phi}{}^{-\overline{z}_{ij}/2}\widetilde{P}_{-\widehat{y}_{ij}}\right\}\mathcal{M}[\mathfrak{g}_{ij}^{(II)}]\,.
\end{aligned}
\end{equation}

Similarly for diagonal elements $(\mathfrak{g}_{i} = \mathfrak{g}_{ii})$:

\begin{equation}\label{StepFalseHOMFLYfromGoeritz2} 
\mathcal{M} [\mathfrak{g}_{i}] = \left\{\phi^{-x_{i}/2}\overline{\phi}{}^{\overline{z}_{i}/2}L_{\widetilde{y}_{i}}D_N+\left( \phi^{-x_{i}/2}\overline{P}_{\overline{z}_{i}}+(\overline{\phi}{}^{\overline{z}_{i}/2}+D_N\overline{P}_{\overline{z}_{i}}){P}_{-x_{i}}\right)\left( L_{\widehat{y}_{i}}+D_N \widetilde{P}_{\widehat{y}_{i}}\right) + \phi^{-x_{i}/2}\overline{\phi}{}^{\overline{z}_{i}/2}\widetilde{P}_{-\widehat{y}_{ij}}\right\}\mathcal{M} [\mathfrak{g}_{i}^{(III)}] \,.
\end{equation}

    \begin{wrapfigure}[9]{r}{0.52\textwidth}
   \centering
    \vspace{-1.0cm}
    \includegraphics[width =0.75\linewidth]{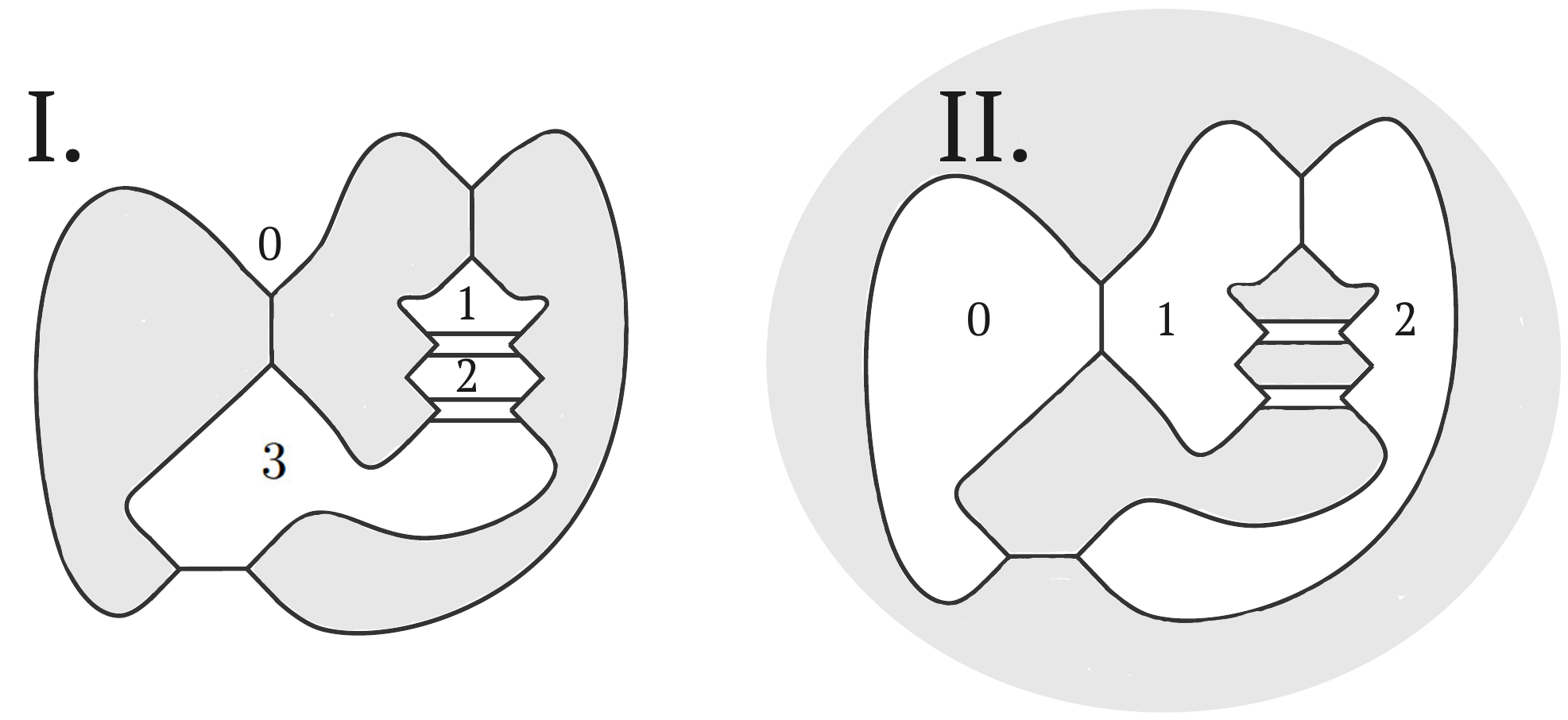}
    \caption{\footnotesize The checkerboard coloring of bipartite diagram of the knot $9_8$.}
        \label{fig:BipKnot98}
\end{wrapfigure}

As the promised example for which the HOMFLYPT polynomial is not recovered by the Goeritz algorithm generalized to three parameters, we consider the knot $9_8$. Its bipartite diagram contains various lock elements, so that the implementation of a three-parameter Goeritz matrix for one of the checkerboard colorings involves all three types of entries (see Fig.~\ref{fig:BipKnot98}). However, difficulties arise for the alternative coloring, where the matrix involves only two parameters. Some of the generalized matrices corresponding to both colorings are shown below:

\begin{equation}\label{TripleGoeritzForKnot98N1}
\widetilde{\mathcal{G}}_{I}^{9_8}  =
\begin{pmatrix}
 \begin{array}{cccc}
 \widetilde{2} - {1}&-\widetilde{1}&0&{1}-\widetilde{1}\\
 -\widetilde{1}&\overline{1}+\widetilde{1}&-\overline{1}&0\\
 0&-\overline{1}&\overline{2}&-\overline{1}\\
 {1}-\widetilde{1}&0&-\overline{1}&\overline{1}-{1}+\widetilde{1}

  \end{array}
 \end{pmatrix} \ \mapsto \
\mathcal{G}_1^{9_8} = \begin{pmatrix}
 \begin{array}{ccc}
\overline{1}+\widetilde{1}&-\overline{1}&0\\
-\overline{1}&\overline{2}&-\overline{1}\\
0&-\overline{1}&\overline{1}-{1}+\widetilde{1}

  \end{array}
 \end{pmatrix}\,, 
\end{equation}
    \begin{equation}\label{TripleGoeritzForKnot98N2}
\widetilde{\mathcal{G}}_{II}^{9_8}  =
\begin{pmatrix}
 \begin{array}{ccc}  \widetilde{1}-{1}&{1}&-\widetilde{1}  \\  {1}&-{2}-\widetilde{2}&{1}+\widetilde{2}\\-\widetilde{1}&{1}+\widetilde{2}&-{1}-\widetilde{1}       
 \end{array}
 \end{pmatrix} \ \mapsto \ \mathcal{G}_2^{9_8}  =
\begin{pmatrix}
 \begin{array}{cc}   \widetilde{1}-{1}&{1}\\{1}&-{2}-\widetilde{2}       
 \end{array}
 \end{pmatrix}, \ \ \ \mathcal{G}_3^{9_8}  =
\begin{pmatrix}
 \begin{array}{cc}  -{2}-\widetilde{2}&{1}+\widetilde{2}\\{1}+\widetilde{2}&-{1}-\widetilde{1}       
 \end{array}
 \end{pmatrix}.
\end{equation}

The HOMFLY--PT polynomial of the knot $9_8$ is correctly recovered using \eqref{StepFalseHOMFLYfromGoeritz1} and \eqref{StepFalseHOMFLYfromGoeritz2} for $\mathcal{G}_1^{9_8}$ and $\mathcal{G}_2^{9_8}$. However, for $\mathcal{G}_3^{9_8}$, this method fails to reproduce the correct polynomial, which is due to the normalization issues in certain branches of the Kauffman bracket resolution described above.

\setcounter{equation}{0}
\section{Goeritz matrix and HOMFLY--PT polynomial of \\ arbitary two-strand link}\label{sec:App2}

\begin{figure}[h!]
		\centering	
		\includegraphics[width =\linewidth]{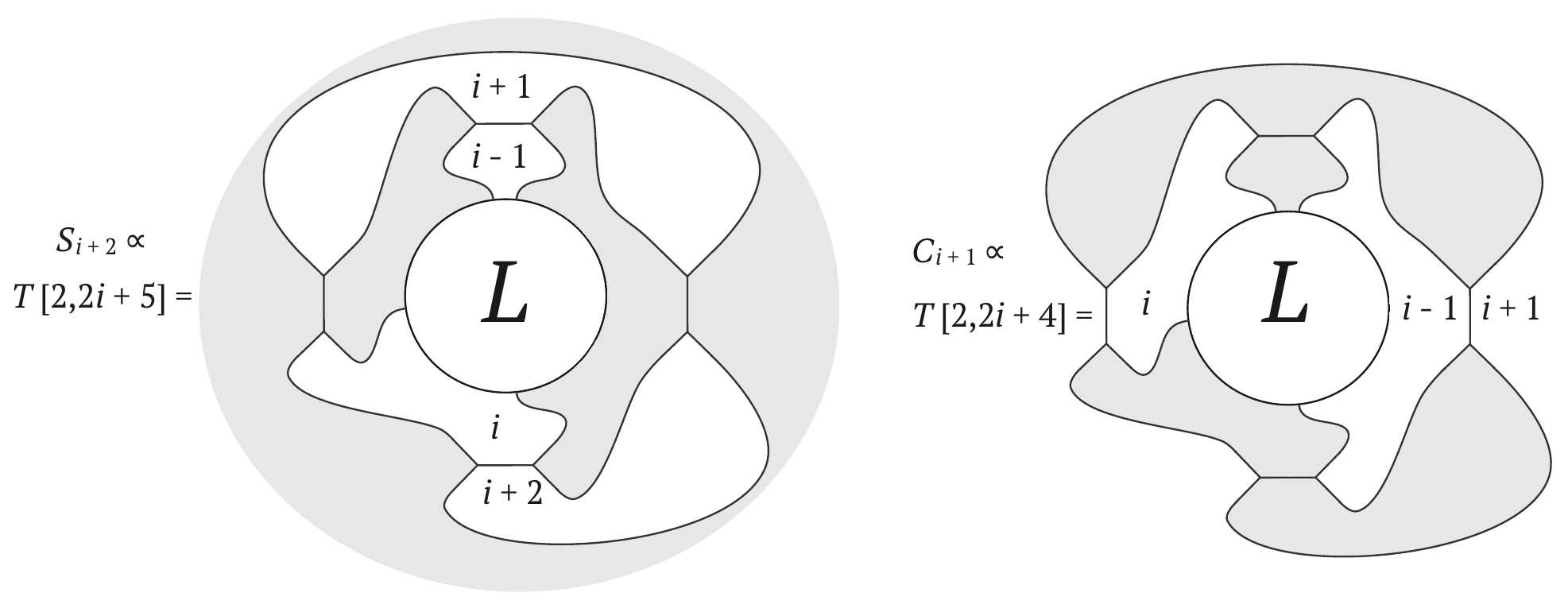}
        \caption{\footnotesize Our chosen colorings for the $S$-type and $C$-type two-strand links.}
        \label{fig:GeneralBipTorCheckboardDiagramm}
	\end{figure}

We continue the demonstration of the Goeritz method by computing the general form of the fundamental HOMFLY--PT polynomials for the torus links $T[2,p]$. Below, we obtain the HOMFLY--PT polynomials for the torus links $T[2,7]$ and $T[2,8]$:

\begin{equation}\label{BipGoeritz1Tor27}
\widetilde{\mathcal{G}}^{T[2,7]}  =
\begin{pmatrix}
 \begin{array}{cccc} \widehat{1}-2 & 1-\widehat{1}  & 1 & 0 \\  1-\widehat{1} & \widehat{2}-2&-\widehat{1} & 1 \\ 1&-\widehat{1}&\widehat{2}-1 & -\widehat{1}  \\ 0 & 1 & -\widehat{1}& \widehat{1} - 1    
 \end{array}
 \end{pmatrix}\quad \mapsto \quad \mathcal{G}^{T[2,7]}  = S_3 = \begin{pmatrix}
 \begin{array}{ccc}  \widehat{2}-2 & -\widehat{1} & 1 \\  -\widehat{1} & \widehat{2}-1 & -\widehat{1} \\ 1 & -\widehat{1} & \widehat{1}-1      
 \end{array}
 \end{pmatrix}\,, 
\end{equation}

\begin{equation}\label{MPolynomialForTor27}
\begin{aligned}
&\mathcal{M}\left[S_3\right] = \left\{\mathcal{U}^2_1[-\widehat{1}]\left(\mathcal{U}_1[1]D_N+\mathcal{U}_2[1]\right)+\mathcal{U}_2[-\widehat{1}]\mathcal{U}_1[1-\widehat{1}]\right\}D_N\mathscr{D}\left[\widehat{1}-1\right]\mathscr{D}\left[ -1\right] +\\
&+\left\{\mathcal{U}_1[-\widehat{1}]\mathcal{U}_2[-\widehat{1}]\left(\mathcal{U}_1[1]D_N+\mathcal{U}_2[1]\right)+\mathcal{U}_2[-\widehat{1}]\mathcal{U}_2[1-\widehat{1}]\right\}D_N\mathscr{D}\left[\widehat{1}-2  \right] ,
\end{aligned}
\end{equation}

\begin{equation}\label{MPolynomialForTor27Two}
\mathcal{M}\left[S_3\right] = \left(\phi^{-3}+6D_N\phi^{-2}+\left(6+9D_N^2\right)\phi^{-1}+4 D_N (4 + D_N^2)+5 (1 + 2 D_N^2)+6D_N\phi^2+\phi^3\right)D_N\,,
\end{equation}

\begin{equation}\label{HOMFLYTor27}
    H^{T[2,7]} =\left(A^{-2}\phi^{1/2}\right)^6\mathcal{M}\left[S_3\right]= \frac{(A^2 - q^2 + A^2 q^4 - q^6 + A^2 q^8 - q^{10} + A^2 q^{12})(A-A^{-1})}{A^8 q^6(q-q^{-1})}\,,
\end{equation}

\begin{equation}\label{BipGoeritz1Tor28}
\widetilde{\mathcal{G}}^{T[2,8]}  =
\begin{pmatrix}
 \begin{array}{cccc} \widehat{1}-2 & 1-\widehat{1}  & 1 & 0 \\  1-\widehat{1} & \widehat{2}-2&-\widehat{1} & 1 \\ 1&-\widehat{1}&\widehat{2}-2 & 1-\widehat{1}  \\ 0 & 1 & 1-\widehat{1}& \widehat{1} - 2    
 \end{array}
 \end{pmatrix}\quad \mapsto \quad \mathcal{G}^{T[2,8]}  = C_3 = \begin{pmatrix}
 \begin{array}{ccc}  \widehat{2}-2 & -\widehat{1} & 1 \\  -\widehat{1} & \widehat{2}-2 & 1-\widehat{1} \\ 1 & 1-\widehat{1} & \widehat{1}-2      
 \end{array}
 \end{pmatrix}\,, 
\end{equation}

\begin{equation}\label{MPolynomialForTor28}
\begin{aligned}
&\mathcal{M}\left[C_3\right] = \left\{\mathcal{U}_1[1-\widehat{1}]\mathcal{U}_1[-\widehat{1}]\left(\mathcal{U}_1[1]D_N+\mathcal{U}_2[1]\right)+\mathcal{U}_2[1-\widehat{1}]\mathcal{U}_1[1-\widehat{1}]\right\}D_N\mathscr{D}\left[\widehat{1}-1\right]\mathscr{D}\left[ -1\right] +\\
&+\left\{\mathcal{U}_1[1-\widehat{1}]\mathcal{U}_2[-\widehat{1}]\left(\mathcal{U}_1[1]D_N+\mathcal{U}_2\right)+\mathcal{U}^2_2[1-\widehat{1}]\right\}D_N\mathscr{D}\left[\widehat{1}-2  \right] ,
\end{aligned}
\end{equation}

\begin{equation}\label{MPolynomialForTor28Two}
\mathcal{M}\left[C_3\right] = \small{\frac{D_N + 4 \phi + 3 D_N^2 \phi + 18 D_N \phi^2 + 3 D_N^3 \phi^2 + 
 10 \phi^3 + 24 D_N^2 \phi^3 + D_N^4 \phi^3 + 25 D_N \phi^4 + 
 10 D_N^3 \phi^4 + 6 \phi^5 + 15 D_N^2 \phi^5 + 
 7 D_N \phi^6 + \phi^7}{\phi^{7/2}}D_N}\,,
\end{equation}

\begin{equation}\label{HOMFLYTor28}
    H^{T[2,8]} =\left(A^{-2}\phi^{1/2}\right)^7\mathcal{M}\left[C_3\right]= \frac{(-q^2 + (1 - q^2 + q^4) (q^4 - q^{10} + A^2 (1 - q^6 + q^{12})))(A-A^{-1})}{A^9 q^7(q-q^{-1})^2}\,.
\end{equation}
It is worth noting that all such matrices can be divided into two types according to the parity of $p$, which differ in the content of their entries. We have chosen the colorings so that the knot diagrams with even $p$ have a white outer region. We refer to this as the $C$-type quaternary Goeritz matrix for two-strand links. Then, for diagrams with odd $p$, which we call the $S$-type, the outer region is shaded, as illustrated in Fig.~\ref{fig:GeneralBipTorCheckboardDiagramm}.

  Now, let us fix the general form of the quaternary Goeritz matrices for two-strand knots corresponding to our chosen coloring.  
  For the $S$-type  $\left(\widetilde{\mathcal{G}}^{T[2,p]}  \quad \mapsto \quad \mathcal{G}^{T[2,p]}  = S_{(p-1)/2} \right),$ the matrices take the following form:

\begin{equation}\label{BipGoeritzMatrixTorKontSType}
\begin{pmatrix}

\widehat{1}-2 & 1-\widehat{1} & 1&0 & \dots & { \ }&\dots& 0 \\
1-\widehat{1} & \widehat{2}-2 & -\widehat{1} & \ddots &\ddots & { \ }& { \ } & \vdots \\
1 & -\widehat{1}& \widehat{2}-2  & \ddots & \ddots &\ddots & { \ }& { \ }\\
0&\ddots&\ddots&\ddots&\ddots&\ddots&\ddots&\vdots \\

\vdots  & \ddots & \ddots&\ddots & \ddots & \ddots&\ddots& 0 \\
{ \ } & { \ } & \ddots & \ddots & \ddots& \widehat{2}-2&-\widehat{1}&1
\\
\vdots & { \ } &{ \ } &\ddots & \ddots & -\widehat{1} & \widehat{2}-1&-\widehat{1} \\
0 & \dots &  &\dots & 0 &1 &-\widehat{1} & \widehat{1}-1  \\
\end{pmatrix}\quad  \mapsto \quad 
\begin{pmatrix}
\begin{array}{ccccccc}
 \widehat{2}-2 & -\widehat{1} & 1 &0 & \dots& \dots & 0 \\
 -\widehat{1}& \widehat{2}-2  &  -\widehat{1} & \ddots &\ddots & \vdots & { \ }\\
1& -\widehat{1}&\widehat{2}-2&\ddots&\ddots&\ddots&\vdots \\

 0 & \ddots&\ddots & \ddots & \ddots&\ddots& 0 \\
\vdots & \ddots & \ddots & \ddots& \widehat{2}-2&-\widehat{1}&1
\\
 \vdots & { \ } &\ddots & \ddots & -\widehat{1} & \widehat{2}-1&-\widehat{1} \\
0 & \dots &\dots & 0 &1 &-\widehat{1} & \widehat{1}-1 \\
\end{array}
\end{pmatrix}.
\end{equation}

A similar form is taken by the quaternary Goeritz matrix for the $C$-type $\left(\widetilde{\mathcal{G}}^{T[2,p]}  \quad \mapsto \quad \mathcal{G}^{T[2,p]}  = C_{(p-2)/2} \right)$:

\begin{equation}\label{BipGoeritzMatrixTorKontCType}
\begin{pmatrix}
\widehat{1}-2 & 1-\widehat{1} & 1&0 & \dots &{ \ } &\dots& 0 \\
1-\widehat{1} & \widehat{2}-2 & -\widehat{1} & \ddots &\ddots &{ \ } &{ \ } & \vdots \\
1 & -\widehat{1}& \widehat{2}-2  & \ddots & \ddots &\ddots &{ \ } &{ \ }\\
0&\ddots&\ddots&\ddots&\ddots&\ddots&\ddots&\vdots \\

\vdots  & \ddots & \ddots&\ddots & \ddots & \ddots&\ddots& 0 \\
{ \ } &{ \ } & \ddots & \ddots & \ddots& \widehat{2}-2&-\widehat{1}&1
\\
\vdots &{ \ } & 
{ \ } &\ddots & \ddots & -\widehat{1} & \widehat{2}-2&1-\widehat{1} \\
0 & \dots & { \ }  &\dots & 0 &1 &1-\widehat{1} & \widehat{1}-2  

\end{pmatrix} \quad  \mapsto \quad 
\begin{pmatrix}
\begin{array}{ccccccc}
 \widehat{2}-2 & -\widehat{1} & 1 &0 & \dots & \dots & 0\\
 -\widehat{1}& \widehat{2}-2  & -\widehat{1}& \ddots &\ddots & { \ } & \vdots\\
1&-\widehat{1}&\widehat{2}-2&\ddots&\ddots&\ddots&\vdots \\

 0 & \ddots&\ddots & \ddots & \ddots&\ddots& 0 \\
 \vdots& \ddots & \ddots & \ddots& \widehat{2}-2&-\widehat{1}&1
\\
\vdots &{ \ } &\ddots & \ddots & -\widehat{1} & \widehat{2}-2&1-\widehat{1} \\
0 & \dots &\dots & 0 &1 &1-\widehat{1} & \widehat{1}-2 \\
\end{array}
\end{pmatrix}
\end{equation}

Note that the action of the function $\mathcal{M}$ on the off–diagonal elements of the matrix $S_n$ allows one to express $\mathcal{M}[S_n]$ in terms of $S_{n-1}$ or $C_{n-1}$.  
Similarly, the action of $\mathcal{M}$ on the off–diagonal elements of $C_n$ yields:

\begin{equation}\label{MPolynomialForSnFromCandS}
    \mathcal{M}\left[S_n\right] = \mathcal{U}_1[-\widehat{1}]\left(\mathcal{U}_1[1]D_N+\mathcal{U}_2[1]\right)\mathcal{M}\left[S_{n-1}\right] +\mathcal{U}_2[-\widehat{1}]\mathcal{M}\left[C_{n-1}\right] = \frac{1+D_N\phi}{\phi}\mathcal{M}\left[S_{n-1}\right] +\phi^{1/2}\mathcal{M}\left[C_{n-1}\right] \,,
\end{equation}

\begin{equation}\label{MPolynomialForCnFromCandS}
    \mathcal{M}\left[C_n\right] = \mathcal{U}_1[1-\widehat{1}]\left(\mathcal{U}_1[1]D_N+\mathcal{U}_2[1]\right)\mathcal{M}\left[S_{n-1}\right] +\mathcal{U}_2[1-\widehat{1}]\mathcal{M}\left[C_{n-1}\right] = \frac{1+D_N\phi}{\phi^{1/2}}\mathcal{M}\left[S_{n-1}\right] +\frac{1+D_N\phi+\phi^2}{\phi}\mathcal{M}\left[C_{n-1}\right] \,.
\end{equation}

Thus, $\mathcal{M}\left[S_n\right]$ and $ \mathcal{M}\left[C_n\right]$ are expressed iteratively in terms of $\mathcal{M}\left[S_3\right]$ and $ \mathcal{M}\left[C_3\right]$, which have been already computed above in \eqref{MPolynomialForTor27Two} and \eqref{MPolynomialForTor28Two}:

\begin{equation}\label{MPolynomialForCnAndSnFromC3AndS3}
\begin{pmatrix}
 \begin{array}{c} \mathcal{M}\left[S_n\right]   \\   \mathcal{M}\left[C_n\right]  
 \end{array}
 \end{pmatrix} = \begin{pmatrix}
 \begin{array}{cc}  \phi^{-1}+D_N & \phi^{1/2} \\ \phi^{-1/2}+D_N\phi^{1/2} & {\phi^{-1}+D_N+\phi}      
 \end{array}
 \end{pmatrix}^{n-3}\begin{pmatrix}
 \begin{array}{c} \mathcal{M}\left[S_{3}\right]   \\   \mathcal{M}\left[C_{3}\right]  
 \end{array}
 \end{pmatrix}\,. 
\end{equation}
Taking into account the normalization $\left(A^{-2}\phi^{1/2}\right)^{2n}$ for diagrams of $S$-type and $\left(A^{-2}\phi^{1/2}\right)^{2n+1}$ for diagrams of $C$-type, one easily arrives at the HOMFLY--PT polynomials for arbitrary two-strand links:

\begin{equation}\label{HOMFLYPTForCnAndSnFromC3AndS3}
\begin{pmatrix}
 \begin{array}{c}  H^{T[2,2n+1]}    \\  H^{T[2,2n+2]}  
 \end{array}
 \end{pmatrix} = \begin{pmatrix}
 \begin{array}{cc}  \phi^{-1}+D_N & \phi^{1/2} \\ \phi^{-1/2}+D_N\phi^{1/2} & {\phi^{-1}+D_N+\phi}      
 \end{array}
 \end{pmatrix}^{n-3}\begin{pmatrix}
 \begin{array}{c} \left(A^{-2}\phi^{1/2}\right)^{2n} \mathcal{M}\left[S_{3}\right]  \\   \left(A^{-2}\phi^{1/2}\right)^{2n+1} \mathcal{M}\left[C_{3}\right]
 \end{array}
 \end{pmatrix}\, 
\end{equation}
In explicit form, they are rewritten as follows:

\begin{equation}\label{HOMFLYTor2np1}
    H^{T[2,2n+1]} = \frac{\left\{q^{2+2n}A^{-2n}(A^2q^2-1)+(A^2-q^2)A^{-2n}q^{-2n}\right\}(A-A^{-1})}{A^2 (q^4-1)(q-q^{-1})}\,,
\end{equation}

\begin{equation}\label{HOMFLYTor2np2}
    H^{T[2,2n+2]} = \frac{\left\{q^{4+2n}A^{-2n}(A^2q^2-1)+(A^2-q^2)A^{-2n}q^{-2n}\right\}(A-A^{-1})}{A^3q (q^4-1)(q-q^{-1})}\,.
\end{equation}

It is easy to see that these general expressions are also valid for all six two-strand torus links considered separately in \eqref{BipHOMFLYPTT23}, \eqref{BipHOMFLYPTT24}, \eqref{HOMFLYTor25}, \eqref{HOMFLYTor26}, \eqref{HOMFLYTor27}, and \eqref{HOMFLYTor28}. Thus, Goeritz method turns out to be a convenient tool that allows one to derive the HOMFLY--PT polynomials of two-strand torus links, in a straightforward and uniform way.

\end{document}